\documentclass[11pt,draft]{amsart}
\input epsf.tex

\usepackage[all]{xy}
\usepackage{amsthm,array,amssymb,amscd,amsfonts,latexsym}

\DeclareFontFamily{OT1}{wncyr}{\hyphenchar\font45}
\DeclareFontShape{OT1}{wncyr}{m}{n}{%
   <5> <6> <7> <8> <9> gen * wncyr
   <10> <10.95> <12> <14.4> <17.28> <20.74>  <24.88>wncyr10}{}
\DeclareFontShape{OT1}{wncyr}{m}{it}{%
   <5> <6> <7> <8> <9> gen * wncyi
   <10> <10.95> <12> <14.4> <17.28> <20.74> <24.88> wncyi10}{}
\DeclareFontShape{OT1}{wncyr}{m}{sc}{%
   <5> <6> <7> <8> <9> <10> <10.95> <12> <14.4>
   <17.28> <20.74> <24.88>wncysc10}{}
\DeclareFontShape{OT1}{wncyr}{b}{n}{%
   <5> <6> <7> <8> <9> gen * wncyb
   <10> <10.95> <12> <14.4> <17.28> <20.74> <24.88>wncyb10}{}
\input cyracc.def
\def\rus{\usefont{OT1}{wncyr}{m}{n}\cyracc\fontsize{9}{11pt}\selectfont}

\DeclareMathSizes{9}{9}{7}{5} % This only is different.

\newtheorem{thm}[equation]{Theorem}
\newtheorem{prop}[equation]{Proposition}
\newtheorem{lem}[equation]{Lemma}

\newtheorem{question}[equation]{Question}
\newtheorem{problem}[equation]{Problem}

\theoremstyle{definition}
\newtheorem{defn}[equation]{Definition}
\newtheorem{remark}[equation]{Remark}
\newtheorem{example}[equation]{Example}

\numberwithin{equation}{section}

\date{\today}
\begin{document}
%%%%%%%%%%%%%%%%%%%%%%%%%%%%%%%%%%%%%%%%%%

\newcommand{\Span}{\operatorname{Span}}
\newcommand{\Symp}{\mbox{\boldmath$\rm Sp$}}
\newcommand{\g}{\mathfrak{g}}
\newcommand{\el}{\mathfrak{l}}
\newcommand{\lt}{\mathfrak{t}}
\newcommand{\lc}{\mathfrak{c}}
\newcommand{\lu}{\mathfrak{u}}
\newcommand{\lr}{\mathfrak{r}}
\newcommand{\Id}{\operatorname{id}}
\newcommand{\id}{\operatorname{id}}
\newcommand{\pr}{\operatorname{pr}}
\newcommand{\Hom}{\operatorname{Hom}}
\newcommand{\sign}{\operatorname{sign}}

\newcommand{\bbA}{{\mathbb A}}
\newcommand{\bbC}{{\mathbb C}}
\newcommand{\bbZ}{{\mathbb Z}}
\newcommand{\bbP}{{\mathbb P}}
\newcommand{\bbQ}{{\mathbb Q}}
\newcommand{\bbG}{{\mathbb G}}
\newcommand{\bbN}{{\mathbb N}}
\newcommand{\bbF}{{\mathbb F}}

\newcommand{\e}{{\lambda}}

\newcommand{\ve}{{\varepsilon}}
\newcommand{\vp}{{\varpi}}

\newcommand{\kbar}{\overline k}

\newcommand{\mo}{\mathopen<}
\newcommand{\mc}{\mathclose>}

\newcommand{\Ker}{\operatorname{Ker}}
\newcommand{\Aut}{\operatorname{Aut}}
\newcommand{\sdp}{\mathbin{{>}\!{\triangleleft}}} % <--- semidirect product
\newcommand{\Alt}{\operatorname{A}}   % < ---- the alternating group
\newcommand{\GL}{\mbox{\boldmath$\rm GL$}}
\newcommand{\PGL}{\mbox{\boldmath$\rm PGL$}}
\newcommand{\SL}{\mbox{\boldmath$\rm SL$}}
\newcommand{\U}{\mbox{\boldmath$\rm U$}}

\newcommand{\rank}{\operatorname{rank}}
\newcommand{\aut}{\operatorname{Aut}}
\newcommand{\Char}{\operatorname{\rm char\,}} %% \char is already a command
\newcommand{\diag}{\operatorname{\rm diag}}
\newcommand{\Gal}{\operatorname{Gal}}
\newcommand{\galois}{\Gal}
\newcommand{\lra}{\longrightarrow}
\newcommand{\SO}{\mbox{\boldmath$\rm SO$}}
\newcommand{\M}{\operatorname{M}}        % <--- matrix algebra
\newcommand{\ord}{\mathop{\rm ord}\nolimits}
\newcommand{\Sym}{{\operatorname{S}}}    % <--- the symmetric group
\newcommand{\tr}{\operatorname{\rm tr}}
\newcommand{\trace}{\tr}
\newcommand{\ad}{\operatorname{Ad}}

\newcommand{\Res}{\operatorname{Res}}
\newcommand{\Sha}{\mbox{\rus{\fontsize{11}{11pt}\selectfont{SH}}}}
\newcommand{\G}{\mathcal{G}}
\renewcommand{\H}{\mathcal{H}}
\newcommand{\gen}[1]{\langle{#1}\rangle}
\renewcommand{\O}{\mathcal{O}}
\newcommand{\C}{\mathcal{C}}
\newcommand{\Ind}{\operatorname{Ind}}
\newcommand{\End}{\operatorname{End}}
\newcommand{\Spin}{\mbox{\boldmath$\rm Spin$}}
\newcommand{\T}{\mathbf G}
\newcommand{\GT}{\mbox{\boldmath$\rm T$}}
\newcommand{\Inf}{\operatorname{Inf}}
\newcommand{\Tor}{\operatorname{Tor}}
\newcommand{\m}{\mbox{\boldmath$\mu$}}
\newcommand{\Cr}{\operatorname{Cr}}

\newcommand{\A}{{\sf A}}
%%{\mbox{\sf \fontsize{11pt}{0mm}\selectfont
%%A}}

\newcommand{\D}{{\sf D}}
%%{\mbox{\sf \fontsize{11pt}{0mm}\selectfont
%%D}}

\newcommand{\Lbd}{{\sf \Lambda}}
%%{\mbox{\sf \fontsize{11pt}{0mm}
%%\selectfont \Lambda}}

%%%%%%%%%%%%%%%%%%%%%%%%%%%%%%%%%%%%%%

\title[Cayley groups]{Cayley groups}
\author{Nicole Lemire}
\address{Department of Mathematics, University of Western Ontario,
London, Ontario N6A 5B7, Canada}
\email{nlemire@uwo.ca}

\author{ Vladimir L. Popov}
\address{Steklov Mathematical Institute,
Russian Academy of Sciences, Gubkina 8, Moscow
119991, Russia} \email{popovvl@orc.ru}

\author{Zinovy Reichstein}
\address{Department of Mathematics, University of British Columbia,
       Vancouver, BC V6T 1Z2, Canada}
\email{reichst@math.ubc.ca}
\urladdr{www.math.ubc.ca/$\stackrel{\sim}{\phantom{.}}$reichst}
\thanks{The first and last authors
were supported in
 part by  NSERC research grants. \\
${\ }$\quad The second author
%% was supported in part by ETH, Z\"urich,
%% Switzerland and Russian grant {\rus
%%N{SH}--123.2003.01}.
was supported in part by ETH, Z\"urich,
 Switzerland, Russian grants {\rus RFFI
05--01--00455}, {\rus N{SH}--123.2003.01}, and
a (granting) program of the Mathematics Branch of the Russian
Academy of Sciences.} 


\begin{abstract} The classical Cayley map,
$X \mapsto (I_n-X)(I_n+X)^{-1}$, is a
birational isomorphism between the special
orthogonal group $\SO_n$ and its Lie algebra
${\mathfrak so}_n$, which is
$\SO_n$-equivariant with respect to the
conjugating and adjoint actions, respectively.
We ask whether or not maps with these
properties can be constructed for other
algebraic groups.
%% $G$.
We show that the answer is usually ``no", with
a few exceptions. In particular, we show that
a Cayley map for the group $\SL_n$ exists if
and only if $n \leqslant 3$, answering an old
question of {\sc Luna}.
\end{abstract}

\maketitle

\begin{quote}{\fontsize{8pt}{3.7mm}
\selectfont\tableofcontents}
\end{quote}

\section {\bf Introduction}
\label{sect.intro}

The exponential map is a fundamental
instrument of Lie theory that yields local
li\-nea\-rization of various problems
involving Lie groups and their actions;
see\;\cite{Bourbaki}.  Let $L$ be a real Lie
group with Lie algebra $\el$. As the
differential at $0$ of the exponential
$\operatorname{\rm exp}\!:\el\rightarrow L$ is
bijective, $\operatorname{\rm exp}$ yields a
diffeomorphism of an open neighborhood of $0$
in $\el$ onto an open neighborhood $U$ of the
identity element $e$ in $L$. The inverse
diffeomorphism $\lambda$ (logarithm) is
equivariant with respect to the action of $L$
on $\el$ via the adjoint representation
$\ad_L\!: L\rightarrow \Aut\el$ and on $L$ by
conjugation, i.e., $\lambda(gug^{-1})=
\ad^{}_L \hskip -.4mm g \hskip .7mm
(\lambda(u))$ if $g\in L$, $u\in U$ and
$gug^{-1}\in U$. This shows that the
conjugating action of $L$ on its underlying
manifold is linearizable in a neighborhood
of~$e$.

In this paper we study what happens if $L$ is
replaced with a connected linear algebraic
group $G$ over an algebraically closed
field~$k$: what is a natural algebraic
counterpart of $\lambda$ for such $G$ and for
which $G$ does it exist?

In what follows we assume that
$\operatorname{\rm char} k=0$ (in fact in many
places this assumption is either redundant or
can be bypassed by modifying the relevant
proof).

\subsection { The classical Cayley map.}
\label{caley-example} %Subsection 1.1
Let $\g$ be the Lie algebra of $G$. One way to
look at the problem is to replace the
Hausdorff topology in the Lie group setting by
the \'etale topology, i.e., to define the
algebraic counterpart of $\lambda$ as a
$G$-equivariant morphism $G\rightarrow \g$
\'etale at~$e$. Then, at least for reductive
groups, there is no existence problem: such
morphisms always exist; see the Corollary to
Lemma \ref{fixed} below. Properties of some of
them have been studied by {\sc Kostant} and
{\sc Michor} in \cite{komi}; see
Example~\ref{km} below. Note also that a
$G$-equivariant dominant morphism $G \lra \g$
exists for every linear algebraic group $G$;
see Theorem~\ref{thm.gen-morph} below.

In the present paper we look at the problem
differently. Our point of view stems from a
discovery made by {\sc Cayley} in 1846,
\cite{cayley};  cf.\,\cite{weyl},
\cite{postnikov}. It suggests that the most
direct approach, i.e., replacing the Hausdorff
topology by the Zariski topology, leads to
something really interesting.
%%a really interesting notion.
 Namely, let $G$ be the special
orthogonal group,
$$G=\SO_n:=\{X\in
{\rm Mat}_{n\times n}\mid X^{\sf T}X=I_n\},$$
where $I_n$ is the identity $n \times
n$-matrix. Then
%%In this case,
$$\g=\mathfrak{o}_n:=\{Y\in {\rm Mat}_{n\times n}\mid
Y^{\sf T}=-Y\},$$ and the adjoint
representation $\ad_G: G\rightarrow \Aut\g$ is
given by
\begin{equation} \label{adj}
\ad^{}_G \hskip -.25mm g \hskip .7mm
(Y)=gY\hskip -.6mmg^{-1}, \quad g\in G,\
Y\!\in\g.
\end{equation}
{\sc Cayley} discovered that there exists a
birational isomorphism
\begin{equation} \label{e.cayley-map1}
\e\colon G \overset{\simeq}
%%\\
%%\mbox{}}
{\dashrightarrow} \g
\end{equation}
%%that is
equivariant with respect to the conjugating
and adjoint actions of $G$ on the underlying
varieties of $G$ and
%%itself
%%by conjugating  and on
$\g$,
%%by $\ad$
respectively, i.e., such that
\begin{equation} \label{e.cayley-map2}
\e(g X \hskip -.2mm g^{-1}) = \ad^{}_G\hskip
-.25mm g \hskip .7mm (\e (X))
\end{equation}
if $g$ and $X \in G$ and both sides of
\eqref{e.cayley-map2} are defined.
%%More
%%precisely, he explicitly constructed
%%gave
His proof is given by the explicit formula
%%for
defining such $\e$:
\begin{equation} \label{e.classical}
%%\text{$
\e\colon  X \mapsto (I_n-X)(I_n+X)^{-1}
%%$ \quad and \quad
%%$\e^{-1} \colon  Y \mapsto
%%(I_n-Y)(I_n+Y)^{-1}$. }
\end{equation}
(one immediately deduces from
\eqref{e.classical} that $
%%\e^{-1} \colon
Y \mapsto (I_n-Y)(I_n+Y)^{-1}$ is the inverse
of $\e$, and from \eqref{adj} that
\eqref{e.cayley-map2} holds).

\subsection { Basic definitions,
main problem, and examples.} \label{subsect1.2}
Inspired by this example, we introduce the
following definition for an arbitrary
connected linear algebraic group $G$.

\begin{defn}\label{Cayley}
A {\em Cayley map} for $G$ is a birational
isomorphism \eqref{e.cayley-map1}
satisfying~\eqref{e.cayley-map2}. A group $G$
is called a {\em Cayley group} if it admits a
Cayley map. If $G$ is defined over a subfield
$K$ of $k$, then a Cayley map defined over $K$
is called a {\it Cayley $K$-map}. If $G$
admits a Cayley $K$-map, $G$ is called a {\it
Cayley $K$-group}.
\end{defn}

Our starting point was a question, posed in
1975 to the second-named author by {\sc Luna},
\cite{lu2}. Using Definition \ref{Cayley}, it
can be reformulated as follows:

\begin{question}\label{q.luna} For
what $n$ is the special linear group $\SL_n$ a
Cayley group?
\end{question}

It is easy to show  (see Example \ref{ex.sp}
below) that $\SL_2$ is a Cayley group. {\sc
Popov} in \cite{popov-luna} has proved that,
contrary to what was expected,
\cite[Remarque, p.\,14]{lu1}, $\SL_3$ is a
Cayley group as well.

More generally, given Definition \ref{Cayley},
it is natural to pose the following problem:

\begin{problem} \label{c.problem}
{\it Which connected linear algebraic groups
are Cayley groups?}
\end{problem}

Before stating our main results, we will
discuss several examples. Set
\[
\m_d:=\{a\in \T_m\mid a^d=1\}.
\]
This is a cyclic subgroup of order $d$ of the
multiplicative group $\T_m$. Below we use the
same notation $\m_d$ for the central cyclic
subgroup $\{aI_n\mid a\in \m_d\}$ of $\GL_n$.

\begin{example}\label{product}
If $G_1, \ldots, G_n$ are Cayley, then
$G:=G_1\times\ldots\times G_n$ is Cayley (the
converse is false; see
Subsection~\ref{semisimple-to-simple}).
Indeed, if $\g_i$ is the Lie algebra of $G_i$
and $\lambda_i\!:
G_i\stackrel{\simeq}{\dasharrow} \g_i$ a
Cayley map, then
$\g=\g_1\oplus\ldots\oplus\g_n$ and
$\lambda_1\times\ldots\times\lambda_n\!:
G\stackrel{\simeq}{\dasharrow}\g$ is a Cayley
map.\quad $\square$
\end{example}

\begin{example}
\label{generalization1} \label{ex.gl_n}
Consider a finite-dimensional associative
algebra $A$ over $k$ with
%%the
%%unit element
identity element~$1$. Let $\mathcal L_A$
%%${\rm Lie}\,A$
be the Lie algebra whose underlying vector
space is that of $A$ and
%%the
whose Lie bracket is given by
%%the commutator:
\begin{equation}\label{bracket}
[X_1, X_2]:=X_1X_2-X_2X_1.
\end{equation}
The group $$G:=A^*$$ of invertible elements of
$A$ is a connected linear algebraic group
whose underlying variety is an open subset of
that of $A$.
%%containing $1$.
This implies that $\g$ is naturally identified
with ${\mathcal L}_A$, and the adjoint action
is given by formula \eqref{adj}. Hence the
natural embedding $\lambda: A^*\hookrightarrow
{\mathcal L}_A$, $X\mapsto X$, is a Cayley
map. Therefore $G$ is a Cayley group.

Taking $A={\rm Mat}_{n\times n}$, we obtain
that $G:=\GL_n$ is Cayley for every $n
\geqslant 1$. \quad $\square$
\end{example}

\begin{example}
\label{generalization2} \label{ex.pgl_n}
Maintain the notation of Example
\ref{generalization1}. For any element $a\in
A$, denote by $\tr a$ the trace of the
operator $L_a$ of left multiplication of $A$
by~$a$. Since the algebra $A$ is associative,
$a\mapsto L_a$ is a homomorphism of $A$ to the
algebra of linear operators on the underlying
vector space of $A$.  From this and
\eqref{bracket}, we deduce that $k\!\cdot\!1$
is an ideal of ${\mathcal L}_A$, the map
$$
\tau: {\mathcal L}_A\rightarrow
k\!\cdot\!1,\quad a\mapsto {\textstyle
%%\frac{\tr a}{\tr 1}
\tr a}\!\cdot\!1,
$$
is a surjective homomorphism of Lie algebras,
and
\begin{equation}\label{direct}{\mathcal
L}_A=\Ker\tau\oplus k\!\cdot\!1.
\end{equation}

The subgroup $k^{*}\!\cdot\!1$ of $A^*$ is
normal; set
\begin{equation}\label{pgl}
G:=A^*/k^{*}\!\cdot\!1.
\end{equation}
As the Lie algebras of $A^*$ and
$k^{*}\!\cdot\!1$ are, respectively, ${\mathcal
L}_A$ and $k\!\cdot\!1$, it follows from
\eqref{direct} that one can identify $\g$ with
$\Ker\tau$. Let $A^*\rightarrow G$,
$a\mapsto[a]$, be the natural projection. Then
the formula
\begin{equation}\label{trace}
[a]\mapsto { \textstyle \frac{\tr 1}{\tr
a}}a-1
\end{equation}
defines a rational map $\lambda:
G\dasharrow\g=\Ker\tau$. Since $\tr
xax^{-1}=\tr a$ for any $a\in A$, $x\in
A^{*}$, it follows from \eqref{trace} that
\eqref{e.cayley-map2} holds. On the other
hand, \eqref{trace} clearly implies that
\begin{equation}\label{+1}
a\mapsto[a+1]
\end{equation}
is the inverse of $\lambda$. Thus $G$ is a
Cayley group.

If $A$ is defined over a subfield $K$ of $k$,
then the group $G$ and birational isomorphisms
\eqref{trace}, \eqref{+1} are defined over $K$
as well. Hence $G$ is a Cayley $K$-group.

For $A={\rm Mat}_{n\times n}$ this shows that
$\PGL_n$ is a Cayley group for every
$n\geqslant 1$. Note that in this case,
$\frac{\tr 1}{\tr a}=\frac{n}{{\rm Tr}\, a}$,
where ${\rm Tr}\,a$ is the trace of matrix
$a$. Let $K$ be a subfield of $k$. Since every
inner $K$-form $G$ of $\PGL_n$ is given by
\eqref{pgl} for $A=D\otimes_K k$, where $D$ is
an $n^2$-dimensional central simple algebra
over $K$ and the $K$-structure of $A$ is
defined by $D$, cf.\,\cite{kneser}, all inner
$K$-forms of $\PGL_n$ are Cayley $K$-groups.

Setting $A=\bigoplus_{i=1}^s{\rm
Mat}_{n_i\times n_i}$, we conclude that
$\prod_{i=1}^{s}\GL_{n_i}
/k^{*}I_{n_1+\ldots+n_s}$ is a Cayley group.
Here $\prod_{i=1}^{s}\GL_{n_i}$ is
block-diagonally embedded in
${\GL}_{n_1+\ldots+n_s}$. \quad $\square$
\end{example}

\begin{example}
\label{next-example} \label{ex.sp} \label{so}
The following construction was noticed by {\sc
Weil} in \cite[p.~599]{weil}. Namely, maintain
the notation of Example \ref{generalization1}
({\sc Weil} assumed that $A$ is semisimple,
but his construction, presented below, does
not use this assumption). Let $\iota$ be an
involution (i.e., an involutory
$k$-anti\-automorphism) of the algebra $A$.
Set
\begin{equation}\label{algebra}
G:=\{a\in A^*\mid a^\iota
a=1\}\!^{\circ}\end{equation} (as usual,
$S^\circ$ denotes the identity component of an
algebraic group $S$). The Lie algebra of $G$
is the subalgebra of odd elements of
${\mathcal L}_A$ for $\iota$,
$$
\g=\{a\in {\mathcal L}_A\mid a^\iota=-a\}.
$$
The formula
\begin{equation}\label{1-a}
a \mapsto (1-a)(1+a)^{-1}
\end{equation}
defines an equivariant rational map
$\lambda\colon G {\dashrightarrow} \g$, and
the formula
\begin{equation}\label{1-b}
b \mapsto (1-b)(1+b)^{-1}
\end{equation}
defines its inverse, $\lambda^{-1}\colon \g
{\dashrightarrow} G$. Thus $\lambda$
 is a Cayley map and $G$ is a
Cayley group.

If $A$ and $\iota$ are defined over a subfield
$K$ of $k$, then the group $G$ and birational
isomorphisms \eqref{1-a}, \eqref{1-b} are
defined over $K$ as well. Hence $G$ is a
Cayley $K$-group.

For $A={\rm Mat}_{n\times n}$ and the
involution $X\stackrel{\iota}{\mapsto} X^{\sf
T}$, this turns into the classical Cayley
construction for $G={\bf SO}_n$, proving that
this group is Cayley for every $n\geqslant 1$.
In particular, the following groups are
Cayley: $\T_m \simeq \SO_2$ (see Examples
\ref{ex.gl_n} and \ref{ex.abelian}),
$\PGL_2\simeq\SL_2/\m_2\simeq \SO_3$  (see
Example \ref{ex.pgl_n}),
$(\SL_2\times\SL_2)/\m_2\simeq \SO_4$ (here
$\SL_2\times\SL_2$ is block-diagonally
embedded in $\SL_4$), $\Symp_4/\m_2\simeq
\SO_5$, and $\SL_4/\m_2\simeq \SO_6$.

For $A\!=\!{\rm Mat}_{2n\times 2n}$ and the
involution $X\!\stackrel{\iota}{\mapsto}\!
J_{2n}\!\!{}^{-1}X^{\sf T}J_{2n}$, where
\begin{math}J_{2n}\!:=
\!\Bigl[\begin{smallmatrix}0&
I_{n}\\-I_{n}&0\end{smallmatrix}\Bigr]
\end{math},
we have
\begin{equation*}
\begin{gathered}
G = \Symp_{2n}:=\{X\in {\rm Mat}_{2n\times
2n}\mid X^{\sf T}J_{2n}X=J_{2n}\},\\
\g=\mathfrak{sp}_{2n}:=\{Y\in {\rm
Mat}_{2n\times 2n}\mid Y^{\sf
T}J_{2n}=-J_{2n}Y\},
\end{gathered}
\end{equation*}
\vskip 1.5mm \noindent so the contruction
shows that \eqref{e.classical} is a Cayley map
for $\Symp_{2n}$; cf.\;\cite[Examples 6,
7]{postnikov}. Thus $\Symp_{2n}$ is Cayley
for every $n \geqslant 1$. In particular,
$\SL_2\simeq \Spin_{3}\simeq \Symp_2$,
$\Spin_{4}\simeq \SL_2\times\SL_2$, and
$\Spin_{5}\simeq \Symp_4$ are Cayley. Below we
shall prove that $\Spin_n$  is not Cayley for
$n\geqslant 6$.

Let $K$ be a subfield of $k$. Since every
$K$-form $G$ of ${\bf SO}_n$ or ${\bf
Sp}_{2n}$ is given by \eqref{algebra} for some
algebra $A$ and its involution $\iota$, both
defined over $K$, see \cite{weil},
\cite{kneser}, all $K$-forms of ${\bf SO}_n$
and ${\bf Sp}_{2n}$ are Cayley $K$-groups.
\quad $\square$
\end{example}

\begin{example}
\label{ex.abelian} Every connected commutative
linear algebraic group $G$ is Cayley. In fact,
in this case, condition~\eqref{e.cayley-map2}
is vacuous, so the existence of
\eqref{e.cayley-map1} is equivalent to the
property that the underlying variety of $G$ is
rational. {\sc Chevalley} in \cite{c1} proved
that over an algebraically closed field of
characteristic zero this property holds for
any connected linear algebraic group (not
necessarily commutative). In particular, the
al\-geb\-raic torus $\T_m^d$, where $$\T_m^d:=
\underbrace{\T_m\times\ldots\times\T_m}_d
\hskip2.3mm \mbox{if $d\geqslant 1$,} \quad
{\bf G}_m^0=e,
$$
\vskip -1mm \noindent is a Cayley group for
every $d\geqslant 0$ (as ${\bf
G}_m\!=\!\GL_1$, this also follows from
Examples~\ref{product} and \ref{ex.gl_n}).
\end{example}

\begin{example}\label{uni}
Every unipotent linear algebraic group $G$ is
Cayley ($G$ is automatically connected because
$\operatorname{char}k=0$). Indeed, we may
assume without loss of generality that $G
\subset \GL_n$, so that elements of $G$ are
unipotent $n \times n$-matrices, elements of
$\g$ are nilpotent $n \times n$-matrices, and
$\operatorname{Ad}^{}_G$ is given by
\eqref{adj}. So we have $(I_n-X)^n=Y^n=0$ for
any $X\in G$, $Y\in \g$. Hence the exponential
map is given by
 \begin{equation*}\label{unip-exp}
 \textstyle
 \exp: \g\longrightarrow G,\quad
 Y\mapsto \sum_{i = 0}^{n-1}
 \frac{1}{i!} Y^i.
 \end{equation*}
Therefore $\exp$ is a $G$-equivariant morphism
of algebraic varieties. Moreover, it is an
isomorphism since the formula
 \begin{equation} \label{e.ln}
 \textstyle
 \lambda:=\ln \colon G \lra \g,\quad
 X\mapsto - \sum_{i = 1}^{n-1} \frac{1}{i} (I_n-X)^i,
 \end{equation}
defines its inverse.

More generally, by the Corollary of
Proposition~\ref{prop.levi} below, every
connected solvable linear algebraic group is
Cayley. \quad $\square$
\end{example}

\subsection { Notational
conventions.}
\label{not} % Subsection 1.3
In order to formulate our main results, we need
some notation and definitions.

For any algebraic torus $T$, we denote by
$\widehat{T}$ its character group,
\[
\widehat{T}:=\Hom_{\rm alg}(T, \T_m),
\]
written additively. It is a lattice (i.e., a
free abelian group of finite rank).

Let $T$ be a maximal torus of $G$ and let
\begin{gather}\label{torus-notations}
\begin{cases}
\begin{gathered}
N =N_{G, T}
:=\{g\in G\mid gTg^{-1}=T\},\\[1.4mm]
C =C_{G,T} :=\{g\in G\mid gtg^{-1}=t \mbox{
\rm
for all } t\in T\},\\[1.4mm]
W=W_{G} =W_{G, T}:=N/C
\end{gathered}
\end{cases}
\end{gather}
be, respectively, its normalizer, centralizer
(which is the Cartan subgroup of $G$), and the
Weyl group. The group $C$ is the identity
component of $N$, and if $G$ is reductive,
then $C=T$; see \cite[12.1, 13.17]{borel}. The
finite group $W$ naturally acts by
automorphisms of $\widehat{T}$. Since all
maximal tori in $G$ are conjugate, $W$ and the
$W$-lattice $\widehat{T}$ do not depend, up to
isomorphism, on the choice of $T$.

\begin{defn} The $W$-lattice $\widehat{T}$
is called the {\it character lattice} of $G$
and is denoted by~${\mathcal X}_G$.
\end{defn}

\begin{remark} The reader should be
careful about this terminology: the elements
of the character lattice of $G$ are the
characters of $T$, not of $G$.
\end{remark}

\begin{defn}\label{stably-cayley}
A group $G$ is called {\em stably Cayley} if
$G\times\T_m^d$ is Cayley for some $d
\geqslant 0$. If $G$ is defined over a
subfield $K$ of $k$ and $G\times {\bf G}_m^d$
is a Cayley $K$-group for some $d \geqslant
0$, then $G$ is called a {\it stably Cayley
$K$-group}.
\end{defn}

It is easy to see that $G$ is stably Cayley if
and only if $G \times A$ is Cayley for some
connected abelian algebraic group $A$. In what follows
 we will sometimes use
Definition~\ref{stably-cayley} in this form.

\subsection { Main results.}
\label{subsection 1.4} Now we are ready to
state our main results. In what follows we shall
denote the generic torus of $G$ by ${\GT}_G$.
(For the definition of ${\GT}_G$, see
\cite{voskresenskii}, \cite{ck} or
Definition~\ref{def.generic-torus} in
Subsection~\ref{sect.generic-tori}.)

\begin{thm} \label{thm1}
Let $G$ be a connected reductive algebraic
group. Then the following implications hold:
\[
\xymatrix@=16pt@R=-2pt{ \txt{$\mathcal{X}_G$
is
sign- \\
permutation}\ar@{=>}[r]^{\ \hskip 3mm\rm(a)}&
\txt{$G$ is \\
Cayley} \ar@{=>}[r]^{\ \hskip -2mm\rm(b)}&
\txt{${\GT}_G$
is \\
rational} \ar@{=>}[r]^{\ \hskip -6mm\rm(c)}&
\txt{${\GT}_G$
is stably \\
rational}\ar@{<=>}[r]^{\hskip 0mm\rm(d)}&
\txt{
$\mathcal{X}_G$ is quasi- \\
permutation} \ar@{<=>} [r]^{\ \hskip
0.6mm\rm(e)}& \txt{$G$ is stably
\\ Cayley } }.  \]
\vskip 1mm
 \noindent Moreover, the
implications {\rm (a)} and {\rm (b)} cannot be
reversed. In particular, not every stably
Cayley group is Cayley.
\end{thm}

For the definitions of sign-permutation and
quasi-permutation lattices, see
Subsection~\ref{sect.lattices}. Note that it is a
long-standing open question whether or not
every stably rational torus is rational;
see~\cite[p.\;52]{voskresenskii}. In
particular, we do not know whether or not
implication (c) can be reversed. We also
remark that (d) is well known;
see,~e.g.,~\cite[Theorem\;4.7.2]{voskresenskii}.

A proof of Theorem~\ref{thm1} will be given in
Subsection~\ref{sect.pr-of-thm1}. In
Section~\ref{sect.levi} we will partially
reduce Problem~\ref{c.problem} to the case
where $G$ is a simple group.

We will then use Theorem~\ref{thm1} to
translate results about stable rationality of
generic tori into statements about the
existence (and, more often, the non-existence)
of Cayley maps for various simple algebraic
groups (i.e., groups having no proper
connected normal subgroups). In particular,
{\sc Lemire} and {\sc Lorenz} in \cite{ll} and
{\sc Cortella} and {\sc Kunyavski\v \i} in
\cite{ck} have recently proved that the
character lattice of $\SL_n$ is
quasi-permutation if and only if $n \leqslant
3$. (This result had been previously
conjectured and proved for prime $n$ by {\sc
Le Bruyn} in \cite{lebruyn1},
\cite{lebruyn2}.) Theorem~\ref{thm1} now tells
us that $\SL_n$ is not stably Cayley (and thus
not Cayley) for any $n \geqslant 4$. On the
other hand, Example~\ref{ex.sp} shows that
$\SL_2$ is Cayley, and {\sc Popov} in
\cite{popov-luna} has proved that $\SL_3$ is
Cayley as well (an outline of the arguments
from \cite{popov-luna} is reproduced in the
Appendix; see also an explicit construction in
Section~\ref{sect.sl_3}). This settles Luna's
original Question \ref{q.luna} about~$\SL_n$.

In a similar manner, we proceed to classify
the connected simple groups $G$ with
quasi-permutation character lattices
${\mathcal X}_G$. For simply connected and
adjoint groups this was done by {\sc Cortella}
and {\sc Kunyavski\v\i}~in \cite{ck}. In
Sections~\ref{sect.intA_n1}
and~\ref{sect.intD_n} we extend their results
to all other connected simple groups.
Combining this classification with
Theorem~\ref{thm1}, we obtain the following
result.

\begin{thm} \label{thm2}
Let $G$ be a connected simple algebraic group.
Then the following conditions are equivalent:

\begin{enumerate}
\item[\rm (a)] $G$ is stably Cayley; \item[\rm
(b)] $G$ is one of the following groups:
\begin{equation} \label{list_st_Cal}
\text{ $\SL_n$ for $n\leqslant 3$,\hskip 2.5mm
$\SO_{n}$ for $n\neq 2, 4$,\hskip 2.5mm
$\Symp_{2n}$, \hskip 1mm $\PGL_n$, \hskip
1mm${\bf G}_2$.}
\end{equation}
\end{enumerate}
\end{thm}
\begin{remark} The groups ${\bf SO}_2$
and ${\bf SO}_4$ are stably Cayley (and even
Cayley; see Example~\ref{so}) but they are
excluded because they are not simple. Note
also that, due to exceptional isomorphisms,
some groups are listed twice
in~\eqref{list_st_Cal}. (For example, $\Symp_2
\simeq \SL_2$.)
\end{remark}

It is now natural to ask which of the stably
Cayley simple groups listed in
Theorem~\ref{thm2}(b) are in fact Cayley. Here
is the answer:

\begin{thm} \label{thm3}
Let $G$ be a connected simple algebraic group.
\begin{enumerate}
\item[{\rm (a)}] The following conditions are
equivalent:
\begin{enumerate}
\item[\rm (i)] $G$ is Cayley; \item[\rm (ii)]
$G$ is one of the following groups:
\begin{equation}\label{list-Cay}
\text{ $\SL_n$ for $n\leqslant 3$,\hskip 2.5mm
$\SO_{n}$ for $n\neq 2, 4$,\hskip 2.5mm
$\Symp_{2n}$, \hskip 1mm $\PGL_n$. }
\end{equation}
\end{enumerate}
\item[{\rm (b)}] The group
$\mbox{\boldmath${\rm G}$}_2$ is not Cayley
but the group $\mbox{\boldmath$\rm G$}_2
\times \mbox{\boldmath$\rm G$}_m^2$ is Cayley.
\end{enumerate}
\end{thm}

The first assertion of part (b) is based on
the recent work of {\sc
Iskovskikh}~\cite{iskovskih1}. The groups
$\SO_n$, $\Symp_{2n}$, and $\PGL_n$ were shown
to be Cayley in Examples~\ref{ex.sp}
and~\ref{ex.pgl_n}. The groups $\SL_3$ and
$\mbox{\boldmath$\rm G$}_2$ will be discussed
in Section~\ref{sect.sl_3}.

\begin{remark}
Question \ref{q.luna} was inspired by {\sc
Luna}'s interest in the existence (for
reductive~$G$) of ``algebraic linearization''
of the conjugating action in a Zariski
neighborhood of the identity element $e\in G$,
i.e., in the existence of $G$-isomorphic
neighborhoods of $e$ and $0$ in $G$ and~$\g$,
respectively; cf.\,\cite{lu1}. In our
terminology this is equivalent to the
existence of a Cayley map
\eqref{e.cayley-map1} such that $\lambda$ and
$\lambda^{-1}$ are defined at $e$ and $0$
respectively, and $\lambda(e) = 0$. Not all
Cayley maps have this property. However, note
that our proof of Theorem~\ref{thm3} (in
combination with \cite[p.\,13,
Proposition]{lu1}) shows that each of the
simple groups listed in~\eqref{list-Cay}
admits a Cayley map with this property (and so
does any direct product of these groups); see
Examples \ref{product}, \ref{ex.gl_n}, \ref{ex.pgl_n}, \ref{ex.sp}, \ref{ex.abelian}, and \ref{uni}, Subsections
\ref{subssl3} and \ref{sect.g_2}, and the
Appendix.
\end{remark}

\smallskip

Let $K$ be a subfield of $k$. It follows from
Theorems~\ref{thm2} and \ref{thm3} and
Examples~\ref{ex.pgl_n} and \ref{ex.sp} that
classifying simple Cayley (respectively,
stably Cayley) $K$-groups is reduced to
classifying outer $K$-forms of $\PGL_n$ for
$n\geqslant 3$ and $K$-forms of $\SL_3$
(respectively, outer $K$-forms of $\PGL_n$ for
$n\geqslant 3$ and $K$-forms of $\SL_3$ and
${\bf G}_2$) that are Cayley (respectively,
stably Cayley) $K$-groups. Note that not all
of these $K$-forms are Cayley (respectively,
stably Cayley) $K$-groups. Indeed,
Definitions~\ref{Cayley} and ~\ref{stably-cayley}
imply the following special property of Cayley
(respectively, stably Cayley) $K$-groups:
their underlying varieties are rational
(respectively, stably rational) over~$K$. For
some of the specified $K$-forms this property
does not hold:
\begin{example}\label{kpgl}
{\sc Berhuy}, {\sc  Monsurr\`o}, and {\sc
Tignol} in~\cite{be_mo_ti} have shown that for
every $n\equiv 0\,{\rm mod}\,4$, the group
${\bf PGL}_n$ has a $K$-form $G$ of outer type
whose underlying variety is not stably
rational over $K$. Hence $G$ is not a stably
Cayley $K$-group. \quad $\square$
\end{example}
\begin{remark}
The underlying varieties of all outer
$K$-forms of ${\bf PGL}_n$ with odd $n$ are
rational over $K$;
see~\cite{voskresenskii-klyachko}.  Note also
that the underlying variety of any $K$-form of
a linear algebraic group of rank at most $2$
is rational over $K$; e.g.,\,see
\cite[p.\,189]{merkurjev} and  \cite[4.1,
4.9]{voskresenskii}.
\end{remark}

\subsection { Application to Cremona groups.}
\label{subsection 1.5} The Cremona group
$\Cr_d$, i.e., the group of birational
automorphisms of the affine space $\bbA^d$, is
a classical object in algebraic geometry; see
\cite{iskovskih-enc} and the references
therein. Classifying the subgroups of $\Cr_d$
up to conjugacy is an important research
direction originating in the works of {\sc
Bertini}, {\sc Enriques}, {\sc Fano}, and {\sc
Wiman}. Most of the currently known results on
Cremona groups relate to ${\Cr}_2$ and $\Cr_3$
(the case $d=1$ is trivial because
${\Cr}_1=\PGL_2$). For $d \ge 4$ the groups
$\Cr_d$ are poorly understood, and any results
that shed light on their structure are prized
by the experts.

Our results provide some information about
subgroups of ${\Cr}_d$ by means of the
following simple construction. Consider an
action of an algebraic group $G$ on a rational
variety $X$ of dimension $d$. Let $G_0$ be the
kernel of this action. Any birational
isomorphism between $X$ and $\bbA^d$ gives
rise to an embedding $\iota_X \colon G/G_0
\hookrightarrow \Cr_d$. A different birational
isomorphism between $X$ and $\bbA^d$ gives
rise to a conjugate embedding, so $i_X$  is
uniquely determined by $X$ (as a $G$-variety)
up to conjugacy in $\Cr_d$. If $Y$ is another
rational variety on which $G$ acts, then the
embeddings $\iota_X$ and $\iota_Y$ are
conjugate if and only if $X$ and $Y$ are
birationally isomorphic as $G$-varieties.

Now consider a special case of this
construction, where $G$ is a connected linear
algebraic group, $X$ is the underlying variety
of $G$ (with the conjugating $G$-action), $Y =
\g$ (with the adjoint $G$-action), and the
kernel $G_0$ (for both actions) is the center
of $G$; see \cite[3.15]{borel}.
Definition~\ref{Cayley} can now be rephrased
as follows: a connected algebraic group $G$ is
Cayley if and only if the embeddings $\iota_G$
and $\iota_{\g} \colon G/G_0 = {\rm Ad}^{}_GG
\hookrightarrow \Cr_{\dim G}$ are conjugate in
$\Cr_{\dim G}$. In this paper we show that
many connected algebraic groups are not
Cayley; each non-Cayley group $G$ gives rise
to a pair of non-conjugate embeddings of the
form $\iota_G$, $\iota_{\g} \colon {\rm
Ad}^{}_GG \hookrightarrow \Cr_{\dim G}$.

Definition~\ref{stably-cayley} can be
interpreted in a similar manner. For every $d
\ge 1$ consider the embedding $\Cr_d
\hookrightarrow \Cr_{d+1}$ given by writing $
\bbA^{d + 1}$ as $\bbA^d \times \bbA^1$ and
sending an element $g \in \Cr_d$ to $g \times
\id_{\bbA^1} \in \Cr_{d+1}$. Denote the direct
limit for the tower of groups $\Cr_1
\hookrightarrow \Cr_2 \hookrightarrow \ldots$
obtained in this way by $\Cr_{\infty}$.
Suppose $G$ is a group acting on rational
varieties $X$ and~$Y$ (possibly of different
dimensions) with the same kernel $G_0$. Then
it is easy to see that the embeddings $\iota_X
\colon G/G_0 \hookrightarrow \Cr_{\dim X}$ and
$\iota_Y \colon G/G_0 \hookrightarrow
\Cr_{\dim Y}$ are conjugate in $\Cr_{\infty}$
(or equivalently, in $\Cr_m$ for some $m \ge
\max \{ \dim \,X, \,  \dim\,Y \}$) if and only
if $X$ and $Y$ are stably isomorphic as
$G$-varieties.

If $V_1$ and $V_2$ are vector spaces with
faithful linear $G$-actions, then
$\iota_{V_1}$ and $\iota_{V_2}$ are conjugate
in $\Cr_{\infty}$ by the ``no-name lemma'';
cf.  Subsection~\ref{speiserlemma}. We call an
embedding $G \hookrightarrow \Cr_d$ {\em
stably linearizable} if it is conjugate, in
$\Cr_{\infty}$, to $i_V$ for some faithful
linear $G$-action on a vector space $V$.
Definition~\ref{stably-cayley} and the
``no-name lemma'' now  tell us that the
following conditions are equivalent: (a) $G$
is stably Cayley, (b) the embeddings $\iota_G$
and $\iota_{\g} \colon {\rm Ad}^{}_GG
\!\hookrightarrow\! \Cr_{\dim G}$ are
conjugate in $\Cr_{\infty}$, and (c) $\iota_G$
is stably linearizable. Once again, the
results of this paper (and, in particular,
Theorem~\ref{thm2}) can be used to produce
many examples of pairs of embeddings of the
form ${\rm Ad}^{}_GG\! \hookrightarrow\!
\Cr_{\dim G}$ that are not conjugate
in~$\Cr_{\infty}$.

Now suppose that $\Gamma$ is a finite group
and $L$ and $M$ are faithful
$\Gamma$-lattices; see
Subsection~\ref{sect.lattices}. Then $\Gamma$
acts on their dual tori, which we will denote
by $X$ and $Y$.  It now follows from
Lemma~\ref{stabiso} that the embeddings
$\iota_X \colon \Gamma \hookrightarrow
\Cr_{\rank L}$ and $\iota_Y \colon \Gamma
\hookrightarrow \Cr_{\rank M}$ are conjugate
in $\Cr_{\infty}$ if and only if $L$ and $M$
are equivalent in the sense of
Definition~\ref{equivalent}. Taking $M$ to be
a faithful permutation lattice, we conclude
that the embedding $\iota_X \colon \Gamma
\hookrightarrow \Cr_{\rank X}$ is stably
linearizable if and only if $L$ is
quasi-permutation (cf. Definition~\ref{qp} and
the Corollary to Lemma~\ref{stabiso}).

In the special case where $L = \mathcal{X}_G$
is the character lattice of the algebraic group
$G$, $\Gamma = W_G$ is the Weyl group, and $X
= T$ is a maximal torus with Lie algebra
$\mathfrak t$, we see that the following
conditions are equivalent: (a) $G$ is stably
Cayley, (b) $\mathcal{X}_G$ is
quasi-permutation, (c) the embeddings
$\iota_{\lt}$ and $\iota_T \colon W
\hookrightarrow \Cr_{\dim T}$ are conjugate in
$\Cr_{\infty}$, and (d) $\iota_T$ is stably
linearizable. (Note that (a) and (b) are
equivalent by Theorem~\ref{thm1}, and (c) and
(d) are equivalent because the $W$-action on
$\lt$ is linear.) Consequently, every
reductive non-Cayley group $G$ gives rise to a
pair of embeddings $i_T$, $i_{\mathfrak t}
\colon W \hookrightarrow \Cr_{\rank G}$ which
are not conjugate~in~$\Cr_{\infty}$.

\begin{example} Let $G$ be
a simple group of type ${\sf A}_{n-1}$ which
is not stably Cayley, i.e., $G=\SL_n/\m_d$,
where $d \, | \, n$, $d < n$, $n \ge 4$, and
$(n, d) \ne (4, 2)$. Then the embeddings
$\iota_T$ and $\iota_{\lt} \colon \Sym_n
\hookrightarrow \Cr_{n-1}$ are not conjugate
in $\Cr_{\infty}$.

Assume further that $n \ne 6$. Then by
H\"older's theorem (see~\cite{holder}),
$\Sym_n$ has no outer automorphisms. Thus the
images $\iota_T(\Sym_n)$ and
$\iota_{\lt}(\Sym_n)$ are isomorphic finite
subgroups of $\Cr_{n-1}$ which are not
conjugate in $\Cr_{\infty}$. \quad $\square$
\end{example}

\section{\bf Preliminaries}
\label{sect2}

In this section we collect certain preliminary
facts for subsequent use. Some of them are
known and some are new. Throughout this
section $\Gamma$ will denote a group; starting
from Subsection \ref{sect.lattices}, it is
assumed to be finite.

\subsection { \boldmath$\Gamma$-fields and
\boldmath$\Gamma$-varieties.}
\label{f-fields} % Subsection 2.1

In what follows we will use the following
terminology. A $\Gamma$-{\it field} is a field
$K$ together with an action of $\Gamma$ by
automorphisms of $K$. Let $K_1$ and $K_2$ be
$\Gamma$-fields containing a common
$\Gamma$-subfield $K_0$. We say that $K_1$ and
$K_2$ are {\it isomorphic as $\Gamma$-fields
{\rm(or} $\Gamma$-isomorphic\!{\rm)} over}
$K_0$ if there is a $\Gamma$-equivariant field
isomorphism $K_1\rightarrow K_2$ which is the
identity on $K_0$. We say that $K_1$ and $K_2$
are {\it stably isomorphic as $\Gamma$-fields
{\rm(or} stably $\Gamma$-isomorphic{\rm)}
over} $K_0$ if, for suitable $n$ and $m$,
$K_1(x_1,\ldots,x_n)$ and
$K_2(y_1,\ldots,y_m)$ are isomorphic as
$\Gamma$-fields over $K_0$. Here, $x_1, \dots,
x_n$ and $y_1, \dots,y_m$ are algebraically
independent variables over $K_1$ and $K_2$,
respectively; these variables are assumed to
be fixed by the $\Gamma$-action.

If $\Gamma$ is an algebraic group, a {\it
$\Gamma$-variety} is an algebraic variety $X$
endowed with an algebraic (morphic) action of
$\Gamma$. A $\Gamma$-equivariant morphism
(respectively, rational map) of
$\Gamma$-varieties is a {\it
$\Gamma$-morphism} (respectively, {\it
rational $\Gamma$-map}). If $X_1$ and $X_2$
are irreducible $\Gamma$-varieties, then
$k(X_1)$ and $k(X_2)$ are $\Gamma$-fields with
respect to the natural actions of $\Gamma$.
These fields are stably $\Gamma$-isomorphic
over $k$ if and only if there is a birational
$\Gamma$-isomorphism $X_1\times \bbA^r
\dasharrow X_2\times \bbA^s$ for some $r$ and
$s$, where $\Gamma$ acts on $X_1\times \bbA^r$
and $X_2\times \bbA^s$ via the first factors.
In this case, $X_1$ and $X_2$ are called {\it
stably birationally $\Gamma$-isomorphic}.

\subsection{ \boldmath$\Gamma$-lattices.}
\label{sect.lattices} % Subsection 2.2.
 From now on we
assume that $\Gamma$ is a finite group.

A {\it lattice} $L$ of rank $r$ is a free
abelian group of rank $r$. A $\Gamma$-{\it
lattice} is a lattice equipped with an action
of $\Gamma$ by automorphisms. It is called
{\it faithful} (respectively, {\it trivial}) if
the homomorphism $\Gamma \rightarrow {\rm
Aut}_{\bbZ} L$ defining the action is
injective (respectively, trivial). If $H$ is a
subgroup of $\Gamma$, then $L$ considered as
an $H$-lattice is denoted by $L|_H$.

Given a group $H$ and a ring $R$, we denote by
$R[H]$ the group ring of $H$ over $R$. If $K$
is a field and $L$ is a $\Gamma$-lattice, we
denote by $K(L)$ the fraction field of $K[L]$;
both $K[L]$ and $K(L)$ inherit a
$\Gamma$-action from $L$. We usually think of
these objects multiplicatively, i.e., we
consider the set of symbols $\{x^a\}_{a\in L}$
as a basis of the $K$-vector space $K[L]$, and
the multiplication being defined by
$x^{a}x^{b}=x^{a+b}$. So $\sigma\cdot
x^a=x^{\sigma\cdot a}$ for any $\sigma\in
\Gamma$. If $a^{}_1,\ldots, a^{}_r$ is a basis
of $L$ and $x^{}_i:=x^{a^{}_i}$, then
$K[L]=K[x^{}_1, {x^{}_1},\hskip -1.5mm^{-
1}\ldots, x^{}_r, {x^{}_r}\hskip -1.5mm^{-
1}]$ and $K(L) = K(x^{}_1, \dots, x^{}_r)$.
Note that any group isomorphism $L\rightarrow
\widehat{{\bf G}_m^r}$ induces the
$K$-isomorphisms of algebras $K[L]\rightarrow
K[{\bf G}_m^r]$ and fields $K(L)\rightarrow
K({\bf G}_m^r)$, and therefore it induces a
$K$-defined algebraic action of $\Gamma$ on
the torus $\T_m^r$ by its automorphisms. Any
such action is obtained in this way.

An important example is $L={\mathcal X}_G$,
the character lattice of a connected algebraic
group $G$, and $\Gamma=W$, the Weyl group of
$G$. In this case, $k({\mathcal X}_G)$ is the
field of rational functions on a maximal torus
of $G$.

\begin{defn} \label{def.permutation}
A $\Gamma$-lattice $L$ is called {\em
permutation} (respectively, {\em
sign-permutation}) if it has a basis
$\ve^{}_1, \dots, \ve^{}_r$ such that the set
$\{\ve^{}_1, \dots, \ve^{}_r \}$
(respectively, $\{\ve^{}_1, -\ve^{}_1, \dots,\break
\ve^{}_r, -\ve^{}_r\}$) is $\Gamma$-stable.
\end{defn}

If $X$ is a finite set endowed with an action
of $\Gamma$, we denote by $\bbZ[X]$ the free
abelian group generated by $X$ and endowed
with the natural action of $\Gamma$.
Permutation lattices may be, alternatively,
defined as those of the form $\bbZ[X]$. Since
$X$ is the union of $\Gamma$-orbits, any
permutation lattice is isomorphic to some
$\bigoplus_{i=1}^s\bbZ[\Gamma/\Gamma_i]$, where
each $\Gamma_i$ is a subgroup of~$\Gamma$.

\begin{defn}\label{equivalent}
(\cite{cs}). Two $\Gamma$-lattices $M$ and $N$
are called \emph{equivalent}, written $M\sim
N$, if they become $\Gamma$-isomorphic after
extending by permutation lattices, i.e., if
there are exact sequences of $\Gamma$-lattices
\begin{equation} \label{E:equiv}
0\longrightarrow M\longrightarrow
E\longrightarrow P\longrightarrow 0
\qquad\text{and}\qquad 0\longrightarrow
N\longrightarrow E\longrightarrow
Q\longrightarrow 0
\end{equation}
where $P$ and $Q$ are permutation lattices.
\end{defn}

For a direct proof that this does indeed
define an equi\-va\-lence relation and for
further background see \cite[Lemma 8]{cs} or
\cite{swan}.

\begin{defn}\label{qp}
A $\Gamma$-lattice $L$ is called {\em
quasi-permutation} if $L \sim 0$ under this
equivalence relation, i.e., $L$ becomes
permutation after extending by a permutation
lattice. In other words, $L$ is
quasi-permutation if and only if there is an
exact sequence of $\Gamma$-lattices
$$0\longrightarrow L
\longrightarrow P\longrightarrow Q
\longrightarrow  0,
$$
where $P$ and $Q$ are permutation lattices.
\end{defn}
%
% It is easily seen that the properties of
% being permutation, sign-permutation and
% quasi-permuta\-tion are preserved under
% passing to $\Gamma$-isomorphic
% $\Gamma$-lattices and that replacing
% equivalent $\Gamma$-lattices with
% $\Gamma$-isomorphic ones yields
% equivalent lattices as well.

\begin{lem}\label{stabiso}
Let $M$ and $N$ be faithful $\Gamma$-lattices
and let $K$ be a field. Then the following
properties are equivalent:
\begin{enumerate}
\item[\rm(i)] $K(M)$ and $K(N)$ are stably
isomorphic as $\Gamma$-fields over $K$;
\item[\rm(ii)] $M\sim N$.
\end{enumerate}
\end{lem}

\begin{proof} See \cite[Proposition
1.4]{ll}; this assertion is also implicit
in~\cite{swan}, \cite{cs}
and~\cite[4.7]{voskresenskii}. \quad $\square$
\renewcommand{\qed}{} \end{proof}

Lemma \ref{stabiso} and Definition \ref{qp}
immediately imply the following.

\medskip
\noindent {\bf Corollary.} {\it Let $L$ be a
faithful $\Gamma$-lattice and let $K$ be a
field. Then the following properties are
equivalent:
\begin{enumerate}
\item[\rm(i)] $K(L)$ is stably isomorphic to
$K(P)$ {\rm(}as a $\Gamma$-field over
$K${\rm)} for some faithful permutation
$\Gamma$-lattice $P$; \item[\rm(ii)] $L$ is
quasi-permutation.
\end{enumerate} }

\smallskip

\subsection { Stable equivalence and
flasque resolutions.}
\label{stableeq} % Subsection 2.3.
In addition to the equivalence relation $\sim$
on $\Gamma$-lattices, we will also consider a
stronger equivalence relation $\approx$ of
stable equivalence. Two $\Gamma$-lattices
$L_1$ and $L_2$ are called {\em stably
equivalent} if $L_1\oplus P_1 \simeq L_2\oplus
P_2$ for suitable permutation
$\Gamma$-lattices $P_1$ and $P_2$.

A $\Gamma$-lattice $L$ is called
\emph{flasque} if ${\rm H}^{-1}(S,L)=0$ for
all subgroups $S$ of $\Gamma$. Every
$\Gamma$-lattice $L$ has a \emph{flasque
resolution}
\begin{equation} \label{E:flasqueres}
0\lra L\lra P\lra Q\lra  0
\end{equation}
with $P$ a permutation $\Gamma$-lattice and
$Q$ a flasque $\Gamma$-lattice. Moreover, $Q$
is determined by $L$ up to stable equivalence:
If $0\to L\to P'\to Q'\to  0$ is another
flasque resolution of $L$, then $Q \approx
Q'$. Following \cite{cs}, we will denote the
stable equivalence class of $Q$ in the flasque
resolution \eqref{E:flasqueres} by
$$\rho(L).$$ Note that by~\cite[Lemme
8]{cs}, for $\Gamma$-lattices $M$, $N$,
\begin{equation}\label{E:flreseq}
M\sim N\hskip 1.5mm\Longleftrightarrow\hskip
1.5mm\rho(M) = \rho(N).
\end{equation}

Dually, every $\Gamma$-lattice $L$ has a
\emph{coflasque resolution}
\begin{equation} \label{E:coflasqueres}
0\lra R\lra P\lra L\lra  0
\end{equation}
with $P$ a permutation $\Gamma$-lattice and
$R$ a \emph{coflasque} $\Gamma$-lattice; that
is, ${\rm H}^{1}(S,R)=0$ holds for all
subgroups $S$ of $\Gamma$. Similarly, $R$ is
determined by $L$ up to stable equivalence.
Note that the dual of a flasque resolution for
$L$ is a coflasque resolution for $L^*$ since
 the finite abelian group
${\rm H}^1(S, L)$ is dual to  ${\rm H}^{-1}(S,
L^*)$. For details, see~\cite[Lemme 5]{cs}.
Note that since ${\rm H}^{\pm 1}$ is trivial
for permutation modules, ${\rm H}^{\pm
1}(\Gamma, L)$ depends only on the stable
equivalence class $[L]$ of $L$ and therefore
is denoted by ${\rm H}^{\pm 1}(\Gamma, [L])$.

\smallskip

Following {\sc Colliot--Th{\'e}l{\`e}ne} and
{\sc Sansuc}, \cite{cs,cs2}, we define
$$
\textstyle \Sha^i(\Gamma ,M)=\bigcap_{a\in
\Gamma}
\Ker\bigl(\Res^{\Gamma}_{\gen{a}}\!\colon {\rm
H}^i(\Gamma, M)\longrightarrow {\rm
H}^i(\gen{a}, M)\bigr)
$$
for any $\bbZ[\Gamma]$-module $M$.  Of
particular interest to us will be the case
where $M$ is a $\Gamma$-lattice $L$ and $i=1$
or $2$.

The following lemma is extracted from
\cite[pp.\;199--202]{cs2}. For a proof, see
also \cite[Lemma 4.2]{ll}.

\begin{lem}\label{shalem}
{\rm (a)} For any exact sequence of
$\bbZ[\Gamma]$-modules
$$0 \longrightarrow M\longrightarrow
P\longrightarrow N\longrightarrow 0$$ with $P$
a permutation projective $\Gamma$-lattice,
${\Sha}^2(\Gamma, M) \simeq {\Sha}^1(\Gamma,
N).$

{\rm (b)} ${\rm H}^1(\Gamma,\rho(L))
\simeq{\Sha}^2(\Gamma, L)$ for any
$\Gamma$-lattice $L$.

{\rm (c)} If $L$ is equivalent to a direct
summand of a quasi-permutation
$\Gamma$-lattice, then \linebreak $\Sha^2(S,
L)=0$ holds for all subgroups $S$ of $\Gamma$.
\end{lem}

In particular, ${\Sha}^2(\Gamma,\,{\cdot}\,)$
is constant on $\sim$-classes.

The following technical proposition will help
us show that certain $\Gamma$-lattices are
equivalent.

\begin{prop}\label{equivlat} Let $X$ and $Y$ be
$\Gamma$-lattices satisfying the exact
sequence
$$0\lra X\lra Y\lra \bbZ/d\bbZ\lra 0$$
where $\Gamma$ acts trivially on $\bbZ/d\bbZ$.

{\rm (a)} If $(d,|\Gamma|)=1$, then $X\oplus
\bbZ\simeq Y\oplus \bbZ$ so that $X\approx Y$
and $X^*\approx Y^*$.

{\rm (b)} If the fixed point sequence
$$0\lra X^S\lra Y^S\lra \bbZ/d\bbZ \lra 0$$
is exact for all subgroups $S$ of $\Gamma$,
then $X^*\sim Y^*$ as $\Gamma$-lattices.
\end{prop}

\begin{proof}
(a) This follows directly from Roiter's form
of Schanuel's Lemma~\cite[31.8]{cr1} applied
to the sequence of the proposition and
$$0\lra \bbZ\stackrel{\times d}
{\lra} \bbZ\lra \bbZ/d\bbZ\lra 0.$$

\smallskip

(b) We claim that any coflasque resolution
$$0\lra C_1\lra P\lra X\lra 0$$
for $X$ can be extended to a coflasque
resolution
$$0\lra C_2 \lra P\oplus Q\lra Y\lra 0$$
for $Y$ so that the following diagram commutes
and has exact rows and columns:
\begin{equation}\label{diagr}
\begin{matrix}
\xymatrix@C=5mm@R=4mm{
&0\ar[d]&0\ar[d]&0\ar[d]& \\
 0\ar[r]
& C_1\ar[r]\ar[d] &P \ar[r]\ar[d]& X \ar[r]
\ar[d] &0
 \\
0\ar[r] & C_2\ar[r]\ar[d]&P\oplus Q
\ar[d]\ar[r] & Y \ar[r] \ar[d] & 0
\\
0\ar[r] & U\ar[r]\ar[d]
&Q\ar[r]\ar[d]&\bbZ/d\bbZ \ar[r]\ar[d]&0
\\
&0&0&0& }
\end{matrix}\hskip 6mm 
\end{equation}
Here $C_1,C_2$ are $\Gamma$-coflasque and $P,
Q$ are $\Gamma$-permutation. Indeed, as is
described in~\cite[Lemme 3]{cs}, given a
surjective homomorphism $\pi_0$ from a
permutation $\Gamma$-lattice $P_0$ to a given
$\Gamma$-lattice $X$, we form a coflasque
resolution of $X$ by  defining a new
permutation $\Gamma$-lattice $P$ containing
$P_0$ as a $\Gamma$-sublattice and a new
surjective homomomorphism $\pi:P\to X$ which
extends
 $\pi_0$
and such that $\Ker\pi$ is coflasque.
Explicitly, take $P=P_0\oplus \oplus_{S
}\bbZ[\Gamma/S]\otimes X^S$ where the sum is
taken over all subgroups $S$ of $\Gamma$ for
which $\pi:P^S\to X^S$ is not a surjection and
such that $\Gamma$ acts on
$\bbZ[\Gamma/S]\otimes X^S$ via the first
factor. Then we take $\pi:P\to X$ to be the
unique  $\Gamma$-map
 such that $\pi\vert_{P_0}=\pi_0$
and such that for each $S$, $\pi(gS\otimes
x)=x$ for $g\in \Gamma$ and $x\in X^S$. Then
$\pi:P\to X$ is a surjective $\Gamma$-map
which maps $P^S$ surjectively onto $X^S$ for
all subgroups $S$ of $\Gamma$ so that ${\rm
H}^1(S,\Ker \pi)=0$ as required. To obtain a
compatible coflasque resolution for $Y$,
extend the surjection from the permutation
lattice $P$ onto $X$ to a surjection from the
permutation lattice $P\oplus Q_0$ onto $Y$ and
then adjust this surjection $P\oplus Q_0\to Y$
to one with a coflasque kernel $P\oplus Q\to
Y$ as above. Then the top two rows are exact
and commutative. The bottom row is obtained
via the Snake Lemma.

Let $S$ be a subgroup of $\Gamma$. Taking
$S$-fixed points in \eqref{diagr}, we obtain
\begin{equation*}
\begin{matrix}
\xymatrix@C=5mm@R=4mm{
&0\ar[d]&0\ar[d]&0\ar[d]& \\
0\ar[r] & C_1^{{}\,S}\ar[r]\ar[d] &P^S
\ar[r]\ar[d]& X^S \ar[r]\ar[d] &0
 \\
0\ar[r] & C_2^{{}\,S}\ar[r]\ar[d] &P^S\oplus
Q^S \ar[d]\ar[r] &
Y^S\ar[r]\ar[d] &0\\
0\ar[r] &U^S\ar[r]\ar[d] &Q^S\ar[r]\ar[d]&
\bbZ/d\bbZ \ar[r]\ar[d] &0
\\
&0&0&0& }
\end{matrix}\hskip 6mm 
\end{equation*}
Since $C_1,C_2$ are coflasque, we find that
the first two rows and columns are exact. By
hypothesis, the third column is exact. Then a
diagram chase shows that the bottom row is
exact. But then this means that $U$ is
coflasque since
$$0\lra U^S\lra Q^S\lra \bbZ/d\bbZ \lra {\rm H}^1(S,U)
\lra {\rm H}^1(S,Q)=0$$ is exact.
Applying~\cite[Lemma 1.1]{ll}  to
$$0\lra U\lra Q\lra \bbZ/d\bbZ\lra 0,$$
we find that $U$ is also quasi-permutation as
it satisfies
$$0\lra U\lra Q\oplus \bbZ \lra \bbZ\lra 0.$$
So as $U$ is coflasque, this sequence splits
and $U$ is in fact stably permutation with
$U\oplus \bbZ\simeq Q\oplus \bbZ$. Combining
this isomorphism with the sequence from the
first column of the first commutative diagram
gives us an exact sequence
$$0\lra C_1\lra C_2\oplus \bbZ\lra  Q\oplus \bbZ\lra 0.$$
Since $C_1$ is coflasque and $Q\oplus \bbZ$ is
permutation, this new sequence splits so that
$C_1\oplus Q\oplus \bbZ \simeq C_2\oplus
\bbZ$. Since
$$0\lra X^*\lra P\lra C_1^*\lra 0,\qquad
0\lra Y^*\lra P\oplus Q\lra C_2^*\lra 0$$ are
flasque resolutions of $X^*$ and $Y^*$, this
implies $\rho(X^*)=\rho(Y^*)$ (i.e., that the
corresponding flasque lattices are stably
equivalent). By~\cite[Lemme~8]{cs}, we
conclude that $X^*\sim Y^*$. \quad $\square$
\renewcommand{\qed}{} \end{proof}

\subsection { Speiser's Lemma.}
\label{speiserlemma} % Subsection 2.4
Let $\pi: Y\to X$ be an algebraic vector
bundle. We call it an algebraic {\it vector
$\Gamma$-bundle} if $\Gamma$ acts on $X$ and
$Y$, the morphism $\pi$ is
$\Gamma$-equivariant, and $g \colon \pi^{-1}(x)
\to \pi^{-1}(g(x))$ is a linear map for every
$x \in X$ and $g \in \Gamma$.

The first of the following related rationality
results is an immediate consequence of the
classical  Speiser's Lemma; the others follow
from the first. In a broader context, when
$\Gamma$ is any algebraic group, results of
this type appear in the literature under the
names of ``no-name method" (\cite{dolgachev})
and ``no-name lemma'' (see \cite{colliot}).

\begin{lem} \label{lem.no-name}
{\rm (a)} Suppose $E$ is a $\Gamma$-field and
$K$ is a $\Gamma$-subfield of $E$ such that
$\Gamma$ acts on $K$ faithfully,
$E=K(x_1,\ldots,x_m)$, and $Kx_1+\ldots+Kx_m$
is $\Gamma$-stable. Then $E=K(t_1,\ldots,
t_m)$, where $t_1,\ldots,t_m$ are
$\Gamma$-invariant elements of
$Kx_1+\ldots+Kx_m$. %%

{\rm(b)}  Let $\pi \colon Y \to X$ be an
algebraic vector $\Gamma$-bundle. Suppose that
$X$ is irreducible and the action of $\Gamma$
on $X$ is faithful. Then $\pi$ is birationally
$\Gamma$-trivial; i.e., there exists a
birational $\Gamma$-isomorphism $\varphi:
Y\stackrel{\simeq}{\dasharrow}X\times k^m$,
where $\Gamma$ acts on $X\times k^m$ via the
first factor, such that the diagram
\begin{equation*}
\begin{matrix}
\xymatrix{Y\ar@{-->}
[rr]^{\varphi}\ar_{\pi}[rd]&&X\times
k^m\ar^{\pi_1}[dl]\\
&X&}
\end{matrix}
%%\quad \text{where}\ \pi_1(x,
%%a)=x,
\end{equation*}
is commutative {\rm(}$\pi_1$ denotes
projection to the first factor{\rm)}.

{\rm (c)} Let $V_1$ and $V_2$ be finite-dimensional 
vector spaces over $k$ endowed
with faithful linear actions of $\Gamma$. Then
$V_1$ and $V_2$ are stably
$\Gamma$-isomorphic.

{\rm(d)} Suppose $L$ is a field and
\[ 0 \lra S \stackrel{\iota}{\lra} N
\stackrel{\tau}{\lra} P \lra 0 \] is an exact
sequence of $\Gamma$-lattices, where $S$ is
faithful and $P$ is permutation. Then the
$\Gamma$-field $L(N)$ is $\Gamma$-isomorphic
over $L$ to the $\Gamma$-field $L(S)(t_1,
\dots, t_r)$, where the elements $t_1, \dots,
t_r$ are $\Gamma$-invariant and algebraically
independent over $L(S)$.
\end{lem}

\begin{proof} Part (a) follows from
Speiser's Lemma, \cite{speiser};
cf.~\cite[Theorem 1]{hk} or~\cite[Appendix
3]{shafarevich}.

(b) Recall that, by definition, algebraic
bundles are locally trivial in the \'etale
topo\-logy, but algebraic vector bundles are
automatically locally trivial in the Zariski
topology; see \cite{serre}. This implies that
after replacing $X$ by a $\Gamma$-stable dense
open subset $U$ and $Y$ by $\pi^{-1}(U)$, we
may assume that $Y=X\times k^m$ (but we do not
claim that $\Gamma$ acts via the first
factor!) and $\pi$ is the projection to the first
factor.

Using the projections $Y\to X$ and $Y\to k^m$,
we shall view $k(X)$ and $k(k^m)$ as subfields
of $k(Y)$. Put $E:=k(Y)$, $K:=k(X)$ and let
$x_1,\ldots,x_m$ be the standard coordinate
functions on $k^m$. If $g\in \Gamma$ and $b\in
X$, then the definition of $\Gamma$-bundle
implies that $g(x_i)\hskip
-.3mm|_{\pi^{-1}(b)}\in k\,x_1\!|\hskip
-.2mm_{\pi^{-1}(b)}+\ldots+ k\,x_m\hskip
-.3mm|\hskip -.2mm_{\pi^{-1}(b)}$. In turn,
this implies %%
that the assumptions of (a) hold. Part (b) now
follows from part (a).

(c) Applying part (b) to the projections
$V_1\leftarrow V_1 \times V_2 \to V_2$, we see
that both $V_1$ and $V_2$ are stably
$\Gamma$-isomorphic to $V_1 \times V_2$.

%%\smallskip
(d) Identify $S$ with $\iota(S)$; then
$K:=L(S)$ is a $\Gamma$-subfield of
$E:=L(N)$. Put $x_1=1\in E$ and choose
$x_2\ldots, x_m\in N\subset E$ such that
$\tau(x_2),\ldots, \tau(x_m)$ is a basis of
$P$ permuted by $\Gamma$. The elements
$x_2,\ldots, x_m$ are algebraically
independent over $K$. If $g\in \Gamma$, then
for every $i$ there is a $j$ such that
$a_{ij}:=g(x_i)-x_j\in \Ker\tau=S\subset K$;
so $g(x_i)=a_{ij}x_1+x_j$. This shows that the
assumptions of (a) hold. The claim (with
$r=m-1$) now follows from part~(a). \quad
$\square$
\renewcommand{\qed}{} \end{proof}

\subsection { Homogeneous fiber spaces.} \label{fiber-spaces} % Subsection 2.5
Let $H$ be an algebraic group and let $S$ be a
closed subgroup of $H$. Consider an algebraic
variety $X$ endowed with an algebraic
(morphic) action of $S$ and the algebraic
action of $S$ on $H\times X$ defined by
\begin{equation}\label{fs}
s(h, x)=(hs^{-1}, s(x)), \quad s\in S,\ (h,
x)\in H\times X.
\end{equation}
 Assume that
there exists a geometric quotient,
\cite{mumford}, \cite[4.2]{popov-vinberg},
\begin{equation}\label{GXH}
H\times X\lra (H\times X)/S.
\end{equation}
This is always the case if every finite subset
of $X$ is contained in an affine open subset
of $X$ (note that this property holds if the
variety $X$ is quasi-projective) 
(\cite[3.2]{serre}; 
cf.\;\cite[4.8]{popov-vinberg}). The variety
$(H\times X)/S$, called a {\it homogeneous
fiber space over $H/S$ with fiber $X$}, is
denoted by $H\times^S\hskip -.7mm X$.
%%Clearly if
If  $H$ is connected and $X$ is irreducible,
then $H\times^S\!\!X$ is irreducible. We
denote by $[h, x]$ the image of a point $(h,
x)\in H\times X$ under the morphism
\eqref{GXH}.

The group $H$ acts on $H\times X$ by left
translations of the first factor.
%%Since
As this action commutes with the $S$-action
\eqref{fs}, the universal property of
geometric quotients implies that the
corresponding $H$-action on $H\times^S\!\!X$,
\begin{equation*}\label{left-action}
h'[h, x]= [h'h, x],\quad h', h\in H,\ x\in X,
\end{equation*}
is algebraic. It also implies that since the
composition of the projection $H\times X\to H$
with the canonical morphism $H\to H/S$ is
constant on $S$-orbits of the action
\eqref{fs},
%%it
this composition induces a morphism
\begin{equation}\label{morphi}
\pi=\pi_{{}_{H, S, X}}: H\times^S\hskip -.6mm
X \lra H/S, \quad [h, x]\mapsto hS.
\end{equation}
This morphism
%%$\pi$
is $H$-equivariant and its fiber over the
point  $o\in H/S$ corresponding to $S$ is
$S$-stable and $S$-isomorphic to $X$; in what follows
 we identify $X$ with this fiber. Since
$H$ acts transitively on $H/S$ and $\pi$ is
$H$-equivariant, the $H$-orbit of any point of
$H\times^S\hskip -.6mm X$ intersects $X$. If
$Z$ is an open (respectively, closed)
$H$-stable subset of $X$, and $\iota:
Z\hookrightarrow X$ is the identity
%%identical
embedding, then $H\times^S\hskip -.6mm Z\to
H\times^S\hskip -.6mm X$, $[h, z]\mapsto [h,
\iota(z)]$, is the embedding of algebraic
varieties whose image is an $H$-stable closed
(respectively, open) subset of $H\times^S\hskip
-.6mm X$. Every $H$-stable closed
(respectively, open) subset of $H\times^S\hskip
-.6mm X$ is obtained in this way.

If the action of $S$ on $X$ is trivial, then
$H\times^S\hskip -.6mm X=H/S\times X$ and
$\pi$ is the projection to the first factor.

The morphism $\pi$ is a locally trivial
fibration in the \'etale topology; i.e., each
point of $H/S$ has an open neighborhood $U$
such that the pull back of
$\pi^{-1}(U)\stackrel{\pi}{\to} U$ over a
suitable \'etale covering $\widetilde U\to U$
is isomorphic to the trivial fibration
$\widetilde U\times X\to \widetilde U$, $(y,
x)\mapsto x$; see \cite[\S2]{serre}, \cite[4.8]{popov-vinberg}. If $X$ is a
$k$-vector space and the action of $S$ on $X$
is linear, then \eqref{morphi} is an algebraic
vector $H$-bundle, so
%%the
%%morphism
$\pi$ is locally trivial in the Zariski
topology; i.e.,
$\pi^{-1}(U)\stackrel{\pi}{\to} U$ is
isomorphic to $U\times X\to U$, $(u, x)\mapsto
x$, for a suitable $U$; see \cite{serre}.

If $\psi$ is a (not necessarily
$H$-equivariant) morphism
(respectively, rational map) of
$H\times^S\hskip -.7mm X$ to $H\times^S\hskip
-.6mm Y$ such that
\begin{equation}\label{commut}
\pi_{{}_{H, S, X}}=\pi_{{}_{H, S,
Y}}\circ\psi,
\end{equation}
 then we say that $\psi$ is a
 %%{\it
morphism {\rm(respectively}, rational
map\,{\rm)} {\it over} $H/S$.

\begin{lem} \label{properties} {\rm (a)} If $\psi:
H\times^S\hskip -.7mm X\to H\times^S\hskip
-.6mm Y$ is an $H$-morphism over $H/S$, then
$\psi|_X$ is an $S$-morphism $X\to Y$. The map
$\psi\mapsto \psi|_X$ is a bijection between
$H$-morphisms $H\times^S\hskip -.7mm X\to
H\times^S\hskip -.6mm Y$ over $H/S$ and
$S$-morphisms $X\to Y$. Moreover, $\psi$ is
dominant {\rm (}\hskip -.5mm re\-spectively,
an isomorphism{\rm)} if and only if $\psi|_X$
is dominant {\rm(}\hskip -.5mm respectively,
an isomorphism{\rm)}.

{\rm (b)} Let $H$ be connected and let $X$ and
$Y$ be irreducible. Then the statements in
{\rm(a)} hold with ``morphism\!'' and
``isomorphism\!'' replaced by, respectively,
``rational map\!'' and ``birational
isomorphism\!''.
\end{lem}
\begin{proof} (a) Since $X=
\pi_{ {H, S, X}}^{-1}(o)$, $Y=\pi_{ {H, S,
Y}}^{-1}(o)$, the first statement follows from
\eqref{commut}. As every $H$-orbit in
$H\times^S\hskip -.7mm X$ intersects $X$ and
$\psi$ is $H$-equivariant, $\psi$ is uniquely
determined by $\psi|_X$.
%%Consider action
%%\eqref{fs} of $H$ on $G\times X$ and the
%%analogous action on $G\times Y$.
If $\varphi: X\to Y$ is an $S$-morphism, then
$H\times X\to H\times Y$, $(h, x)\mapsto (h,
\varphi(x)),$ is a morphism commuting with the
actions of $S$ (defined for $H\times X$ by
\eqref{fs} and analogously for $H\times Y$)
and $H$.
%% for $G\times X$ and
%%similarly for $G\times Y$.
By the universal property of geometric
quotients,
%% factors,
the $H$-map $\psi: H\times^S\hskip -.7mm X\to
H\times^S\hskip -.6mm Y$, $[h, x]\mapsto[h,
\varphi(x)]$, is a morphism over $H/S$. We
have $\psi|_X=\varphi$. The
%%the
same argument proves the last statement.
%% is proved using the same
%%argument.

(b) Since $\psi$ is $H$-equivariant, its
indeterminacy locus is $H$-stable. As every
$H$-orbit in $H\times^S\hskip -.7mm X$
intersects $X$, this locus cannot contain $X$.
Consequently, $\psi|_X: X\dashrightarrow
H\times^S\hskip -.6mm Y$ is a well-defined
rational $S$-map. In view of \eqref{commut},
its image lies in $Y$. Now (b) follows from
(a) because rational maps are the equivalence
classes of morphisms of dense open subsets,
and all $H$-stable open subsets in
$H\times^S\hskip -.7mm X$ are of the form
$H\times^S\hskip -.7mm Z$ where $Z$ is an
$S$-stable open subset of $X$.\quad
 $\square$
\renewcommand{\qed}{}
\end{proof}

\section{\bf Cayley maps, generic tori, and lattices}
\label{sect3}

\subsection { Restricting Cayley
maps to Cartan subgroups.}
\label{sect.cartan} % Subsection 3.1
Let $G$ be a connected linear algebraic group
and let $T$ be its maximal torus. Consider the
Cartan subgroup $C$, its normalizer $N$, and
the Weyl group $W$ defined by
\eqref{torus-notations}. Let $\g$, $\lt$, and
$\lc$ be the Lie algebras of $G$, $T$, and $C$,
respectively.

Since $C$ is the identity component of $N$ and
the Cartan subgroups of $G$ are all conjugate
to each other, \cite[12.1]{borel}, assigning
to a point of $G/N$ the identity component of
its $G$-stabilizer (respectively, the Lie
algebra of this $G$-stabilizer) yields a
bijection between $G/N$ and the set of all
Cartan subgroups in $G$ (respectively, all
Cartan subalgebras in $\g$). So $G/N$ can be
considered as the {\it variety of all Cartan
subgroups in} $G$ (respectively, the {\it
variety of all Cartan subalgebras in} $\g$).

Moreover the Cartan subgroups in $G$
(respectively, the Cartan subalgebras in $\g$)
paramet\-ri\-zed in this way by the points of
$G/N$ naturally ``merge" to form a homogeneous
fiber space over $G/N$ with
%%the
fiber $C$ (respectively,~$\lc$).
%%Namely
More precisely, consider the homogeneous fiber
space $G\times^N\hskip -.7mm C$ over $G/N$
defined by the conjugating action of $N$ on
$C$
%%by conjugation
(respectively, the homogeneous fiber space
$G\times^N\hskip -.7mm \lc$ over $G/N$ defined
by the adjoint action of $N$ on $\lc$). Then
for any $g\in G$, the map $\pi_{ {G, N,
C}}^{-1}(g(o))\to gCg^{-1}$, $[g, c]\mapsto
gcg^{-1}$ (respectively, the map $\pi_{ {G, N,
\lc}}^{-1}(g(o))\to {\rm Ad}^{}_G\,g(\lc)$,
$[g, x]\mapsto {\rm Ad}^{}_G\,g(x)$), is a
well-defined isomorphism (we use the notation
of Subsection \ref{fiber-spaces} for $H=G$,
$S=N$).

Consider the conjugating  and adjoint actions, 
respectively, of $G$ on $G$ and $\g$. Then the
definition of homogeneous fiber space implies
that
\begin{equation}\label{cartan-morphisms}
\gamma_{{}_C}: G\times^N\hskip -.7mm C\lra G,
\ [g, c]\mapsto gcg^{-1}, \qquad \gamma_{\lc}:
G\times^N\hskip -.7mm \lc\lra \g, \
[g,x]\mapsto {\rm Ad}^{}_G\,g(x),
\end{equation}
are well-defined $G$-equivariant maps, and the
universal property of geometric factor implies
that they are morphisms.

\begin{lem}\label{lem.cartan} {\rm (a)}
%%In {\rm \eqref{cartan-morphisms}}, the
The morphisms $\gamma_{{}_C}$ and
$\gamma_{\lc}$ in {\rm
\eqref{cartan-morphisms}} are birational
$G$-isomorphisms.

{\rm (b)} Any rational $G$-maps
$G\times^N\hskip -.7mm C\dashrightarrow
G\times^N\hskip -.7mm \lc$ and
$G\times^N\hskip -.7mm \lc\dashrightarrow
G\times^N\hskip -.7mm C$
%%is a
are rational maps over $G/N$.
\end{lem}
\begin{proof} (a) Since the Cartan subgroups
of $G$ are all conjugate and every element of
a dense open set $U$ in $G$ belongs to a
unique Cartan subgroup, \cite[\S 12]{borel},
every fiber $\gamma_{ C}^{-1}(u)$, where $u\in
U$, is a single point. As ${\rm char}\,k=0$,
this means that $\gamma_{{}_C}$ is a
birational isomorphism. For $\gamma_{{\lc}}$
the arguments are analogous because $\lc$ is a
Cartan subalgebra in $\g$, Cartan subalgebras
in $\g$ are all ${\rm Ad}^{}_G G$-conjugate
and a general element of $\g$ is contained in
a unique Cartan subalgebra, \cite[Ch.
VII]{Bourbaki2}.

(b) Since a general element of $T$
(respectively, $\lt$) is regular, $C$
(respectively, $\lc$) is the unique Cartan
subgroup (respectively, subalgebra) containing
$T$ (respectively, $\lt$), \cite[\S
13]{borel}; see \cite[Ch.~VII]{Bourbaki2}.
This implies that $C$ and $\lc$ are the fixed
point sets of the actions of $T$ on
$G\times^N\hskip -.7mm C$ and $G\times^N\hskip
-.7mm \lc$, respectively. Since the maps under
consideration are $G$-equivariant, this
immediately implies the claim.\quad $\square$
\renewcommand{\qed}{}
\end{proof}

\begin{remark}  The group varieties of $C$
%%(\,respec\-ti\-ve\-ly,
and $\lc$
%%\,) is
are the ``standard relative sections'' of,
respectively, $G$
%%(respectively,
and $\g$
%%)
induced
%%respectively
by the rational $G$-map $\pi_{{}_{G, N,
C}}\circ\gamma_C^{-1}: G\dashrightarrow G/N$
%%(respectively,
and $\pi_{{}_{G, N,
\lc}}\circ\gamma_{\lc}^{-1}: \g\dashrightarrow
G/N$; in particular, this yields the following
isomorphisms of invariant fields:
\begin{equation}\label{restriction}
k(G)^G\stackrel{\simeq}{\longrightarrow}k(C)^N,
f\mapsto f|_{C},\qquad
k(\g)^G\stackrel{\simeq}{\longrightarrow}k(\lc)^N,
f\mapsto f|_{\lc};
\end{equation}
see \cite[Definition (1.7.6) and Theorem
(1.7.5)]{popov}.
\end{remark}

\begin{lem} \label{cor.caley}
%%The following equivalencies hold:
%%\begin{enumerate}
%%\item[\rm (a)]
{\rm (a)} $G$ is Cayley if and only if $C$ and
$\lc$ are birationally $N$-isomorphic.

 {\rm (b)} $G$ is stably Cayley
 if and only if $C$ and $\lc$ are
stably birationally $N$-isomorphic.
\end{lem}
\begin{proof} (a)
By Lemma \ref{properties}, the existence of a
birational $N$-isomorphism
%%$\varphi$
$\varphi  \colon C
\stackrel{\simeq}{\dasharrow} \lc$ implies the
existence of a birational $G$-isomorphism
$\psi: G\times^N\hskip -.7mm
C\stackrel{\simeq}{\dashrightarrow}
G\times^N\hskip -.7mm \lc$. Then Lemma
\ref{lem.cartan} shows that
$\gamma_{\lc}\circ\psi\circ\gamma_C^{-1}:
G\stackrel{\simeq}{\dashrightarrow}\g $ is a
Cayley map.

%%Vice versa
Conversely, let $\lambda:
G\stackrel{\simeq}{\dashrightarrow} \g$ be a
Cayley map. Then $\psi:=\gamma_{\lc}^{-1}\circ
\lambda\circ\gamma_{{}_C}:G\times^N\hskip
-.7mm C\stackrel{\simeq}{\dashrightarrow}
G\times^N\hskip -.7mm \lc$ is a birational
$G$-isomorphism. By Lemma \ref{lem.cartan},
$\psi$ is a rational map over $G/N$. Hence, by
Lemma \ref{properties}, $\psi|_{C}:
C\stackrel{\simeq}{\dashrightarrow} \lc$ is a
birational $N$-isomorphism.

(b) If $C$ and $\lc$ are stably birationally
$N$-isomorphic, it follows from the rationality of
the
%%group
underlying variety of any linear algebraic
torus that for some natural $d$ there exists a
birational $N$-isomorphism
\begin{equation}\label{isom}
%%\varphi:
C\times \T_m^d\stackrel{\simeq}{\dasharrow}\lc
\oplus k^d,
\end{equation}
%%for some~$d$;
where $k^d$ is the Lie algebra of $\T_m^d$ and
$N$ acts on $C\times \T_m^d$ and $\lc \oplus
k^d$ via
%%$\T_m^d$ and $k^d$
$C$ and $\lc$, respectively. Clearly $C\times
\T_m^d$ is the Cartan subgroup of $G\times
\T_m^d$ with normalizer $N\times \T_m^d$ and
Lie algebra $\lc\oplus k^d$, and the
birational isomorphism \eqref{isom} is
$N\times \T_m^d$-equivariant. Now~(a) implies
that $G\times \T_m^d$ is Cayley and hence $G$
is stably Cayley.

%%The other way round
Conversely, assume that $G\times \T_m^d$ is
Cayley for some $d$. Then the above arguments
and (a) show that there exists a birational
$N$-isomorphism \eqref{isom}. Since the group
varieties of $\T_m^d$ and $k^d$ are rational,
this means that $C$ and $\lc$ are stably
birationally $N$-isomorphic. \quad $\square$
\renewcommand{\qed}{}
\end{proof}

For reductive groups, Lemma~\ref{cor.caley}
translates into the statement resulting also
from \cite[p.\,13, Proposition]{lu1}:

\medskip

\noindent{\bf Corollary.} \label{cor.torus}
{\it Let $G$ be a connected reductive linear
algebraic group.
\begin{enumerate}
\item[\rm (a)] %%
$G$ is Cayley if and only if $T$ and $\lt$ are
birationally $W$-isomorphic.
 \item[\rm (b)] $G$
is stably Cayley if and only if $T$ and $\lt$
are stably birationally $W$-isomor\-phic.
\end{enumerate}
}
%%\end{cor}

\begin{proof} Since $G$ is reductive,
$C = T$ and $\lc=\lt$. As $T$ is commutative,
this implies that the actions of $N$ on $T$
and $\lt$ descend to the actions of $W$. The
claim now follows from Lemma~\ref{cor.caley}.
\quad $\square$
\renewcommand{\qed}{}
\end{proof}

\subsection { Generic tori.}
\label{sect.generic-tori} % Subsection 3.2
We now recall the definition of generic tori
in a form suitable for our purposes;
see~\cite[4.1]{voskresenskii}
or~\cite[p.\,772]{ck}. We maintain the
notation of Subsections \ref{fiber-spaces} and
\ref{sect.cartan}.

Assume that $G$ is a connected reductive
linear algebraic group; then $C=T$ and
$\lc=\lt$. According to the discussion in the
previous subsection, $G/N$ may be interpreted
in two ways: first, as the {\it variety of all
maximal tori in} $G$ and, second, as the {\it
variety of all maximal tori in} $\g$. The
maximal torus in $G$ (respectively, in $\g$)
assigned to a point $g(o)\in G/N$ is
$gTg^{-1}$ (respectively, $\ad_G\hskip -.25mm
g(\lt)$); it is naturally identified with the
fiber over $g(o)$ of the morphism $\pi_{{}_{G,
N, T}}: G\times^N\hskip -.7mm T \lra G/N$
(respectively, $\pi_{{}_{G, N, \lt}}:
 G\times^N\hskip
-.7mm \lt \lra G/N$).

\begin{defn} \label{def.generic-torus} The triples
$${\GT}_G:=\big(G\times^N\hskip -.7mm T,
\pi_{{}_{G, N, T}}, G/N\big)\quad \text{and}
\quad
{\boldsymbol{\lt}}_{\g}:=\big(G\times^N\hskip
-.7mm \lt, \pi_{{}_{G, N, \lt}}, G/N\big)$$
are called, respectively, the {\it generic torus
of} $G$ and the {\it generic torus of} $\g$.
\end{defn}

We identify the field $k(G/N)$ with its image
in $k(G\times^N\hskip -.7mm T)$ under the
embedding $\pi_{G, N, T}^*$.

\begin{defn} \label{rational-torus}
The generic torus $\GT_G$ is called {\it
rational} if $k(G\times^N\hskip -.7mm T)$ is a
purely trans\-cendental extension of $
k(G/N)$. If $\GT_{G\times {\bf G}_m^d}$ is
rational for some~$d$, then $\GT_G$ is called
{\it stably rational}.
\end{defn}

Equivalently, $\GT_G$ is called rational if
there exists a birational isomorphism
\begin{equation}\label{rational-T}
G\times^N\hskip -.7mm
T\stackrel{\simeq}{\dashrightarrow}G/N\times
\bbA^r
\end{equation} over $G/N$ (then $r=\dim T$). The
arguments used in the proof of Lemma
\ref{cor.caley}(b) show that stable
rationality of $\GT_G$ is equivalent to the
property that there exists a purely
transcendental field extension $E$ of
$k(G\times^N\hskip -.7mm T)$ such that $E$ is
a purely transcendental extension of $k(G/N)$.
There are groups $G$ such that the generic
torus ${\GT}_G$ is not stably rational (and
hence not rational), \cite{voskresenskii},
\cite{ck}.

Of course, for the generic torus
${\boldsymbol{\lt}}_{\g}$ in $\g$, one could
also introduce the notions analogous to that
in Definition \ref{rational-torus}. However in
the Lie algebra context the rationality
problem of generic tori is quite easy: since
$\pi_{{}_{G, N, \lt}}: G\times^{N}\hskip -.6mm
\lt\lra G/N$ is a vector bundle, it is locally
trivial in the Zariski topology, and hence
${\boldsymbol{\lt}}_{\g}$ is always rational;
i.e., there exists a birational isomorphism
\begin{equation}\label{rational-t}
G\times^N\hskip -.7mm
\lt\stackrel{\simeq}{\dashrightarrow}G/N\times
\bbA^r
\end{equation} over $G/N$.

\subsection { Proof of Theorem~\ref{thm1}.}
\label{sect.pr-of-thm1} % Subsection 3.3

{\it Implication} (a): By the Corollary of
Lemma~\ref{cor.caley}, it is enough to
construct a $W$-equiva\-ri\-ant birational
isomorphism $\varphi \colon T
\!\stackrel{\simeq}{\dasharrow} \!\lt$.

Using the sign-permutation basis of
$\widehat{T}$, we can $W$-equivariantly
identify the maximal torus $T$ with $\T_m^r$,
where $r$ is the rank of $G$ and every $w \in
W$ acts on $\T_m^r$ by
\begin{equation} \label{sign}
(t_1, \dots, t_r) \mapsto
(t_{\sigma(1)}^{\varepsilon_1}, \dots,
t_{\sigma(r)}^{\varepsilon_r} ),
%%\]
\end{equation}
for some $\sigma \in \Sym_r$ and some
$\varepsilon_1, \dots, \varepsilon_r \in \{
\pm 1 \}$ (depending on $w$). The Lie algebra
$\lt$ is the tangent space to $\T_m^r$ at
$e=(1, \dots, 1)$; it follows from
\eqref{sign} that we can identify it
with~$k^r$ where $w$ acts by
\begin{equation} \label{sign1}
(x_1, \dots, x_r) \mapsto (\varepsilon_1
x_{\sigma(1)}, \dots, \varepsilon_r
x_{\sigma(r)} ).
\end{equation}
 From \eqref{sign} and \eqref{sign1} we
easily deduce that the formula
\[ (t_1, \dots, t_r) \mapsto \bigl( (1-t_1) (1+t_1)^{-1},
\dots, (1-t_r) (1+t_r)^{-1} \bigr)
\]
defines a desired birational $W$-isomorphism
$\varphi \colon T
\!\stackrel{\simeq}{\dasharrow} \! \lt$. This
completes the proof of implication (a).

To see that implication (a) cannot be
reversed, consider the group $G:=\SL_3$. First
note that this group is Cayley; see
Proposition~\ref{prop.sl_3}. On the other
hand, $W\simeq\Sym_3$ and since the character
lattice ${\mathcal X}_G$ has rank $2$, it cannot 
be sign-permutation. Indeed, if it were,
then $\Sym_3$ would embed into $(\bbZ/2
\bbZ)^2 \rtimes \Sym_2$, which is impossible.

\smallskip

{\it Implication} (b): By the Corollary of
Lemma~\ref{cor.caley}, there is a birational
$N$-isomor\-phism $T
\!\stackrel{\simeq}{\dasharrow} \! \lt$. By
Lem\-ma~\ref{properties}, this implies that
there is a birational $G$-isomorphism
$G\times^N\hskip -.6mm
T\stackrel{\simeq}{\dasharrow}G\times^N\hskip
-.6mm \lt$ over $G/N$. Its composition with
the birational isomorphism \eqref{rational-t}
is a birational isomorphism \eqref{rational-T}
over $G/N$. Hence ${\GT}_G$ is rational.

To see that implication (b) cannot be
reversed, consider the exceptional group ${\bf
G}_2$.  The generic torus of ${\bf G}_2$ is
rational; see~\cite[4.9]{voskresenskii}. On
the other hand, ${\bf G}_2$ is not a Cayley
group; see~Proposition~\ref{prop.iskovskih}.

\smallskip

{\it Implication} (c):  This is obvious from
the definition.

\smallskip

{\it Equivalence} (d): This is well known;
see,~e.g.,~\cite[Theorem
4.7.2]{voskresenskii}.

\smallskip

{\it Equivalence} (e): Let $V$ be any finite-dimensional 
faithful permutation $W$-module
over $k$ (for instance, the one determined by
the regular representation of $W$). Then
clearly $k(V)=k(P)$ for some permutation
$W$-lattice $P$. Since the action of $W$ on
$\lt$ is faithful, \cite{borel}, we deduce
from Lemma~\ref{lem.no-name}(c) that $k(\lt)$
and $k(P)$ are stably $W$-isomorphic over $k$.
Therefore, since $k(T) = k(\widehat{T})$,
applying the Corollary of
Lemma~\ref{cor.caley} implies that $G$ is
stably Cayley if and only if $k(\widehat{T})$
and $k(P)$ are stably $W$-isomorphic over $k$.
On the other hand, the latter property holds
if and only if the $W$-lattice $\widehat{T}$
is quasi-permutation; see the Corollary of
Lemma \ref{stabiso}, whence the claim. \quad
$\square$

\begin{example}
\label{ex.pgln-a} The character lattice
$\bbZ\A_{n-1}$ of $\PGL_n$ is defined by the
exact sequence
\[ 0 \lra \bbZ\A_{n-1} \lra
\bbZ[\Sym_n/\Sym_{n-1}] \stackrel{ \epsilon}
{\lra} \bbZ \lra 0,
\] where $\epsilon$ is the augmentation
map and the Weyl group $W = \Sym_n$ acts
trivially on $\bbZ$ and naturally on
$\bbZ[\Sym_n/\Sym_{n-1}]$; see Subsection
\ref{epsilon}. Thus $\bbZ\A_{n-1}$ is
quasi-permutation. By Theorem~\ref{thm1}, we
conclude that $\PGL_n$ is stably Cayley. We
know that in fact $\PGL_n$ is even Cayley; see
Example~\ref{ex.pgl_n}. Note though that
$\bbZ\A_{n-1}$ is not sign-permutation if
$n>2$.  Indeed, from the sequence above, we
can show that ${\rm
H}^1(\Sym_n,\bbZ\A_{n-1})\simeq \bbZ/n\bbZ$,
whereas by~\cite[Lemma 4.4]{ll}, a
sign-permutation $\Gamma$-lattice $L$ would
have ${\rm H}^1(\Gamma,L)\cong (\bbZ/2\bbZ)^d$
for some $d\geqslant 0$. \quad $\square$
\end{example}

\section{\bf Reduction theorems}
\label{sect.levi}

The purpose of this section is to show that to
a certain extent classifying arbitrary Cayley
groups is reduced to classifying simple ones.

As before, let $G$ be a connected linear
algebraic group. Denote by $R$ and $R_u$,
respectively, the radical and the unipotent
radical of $G$. Recall that a {\it Levi
subgroup} of $G$ is a connected subgroup $L$,
necessarily reductive, such that $G=L\ltimes
R_u$; since ${\rm char}\,k=0$, Levi subgroups
exist and are conjugate, \cite[11.22]{borel}.

In this section we will address the following
questions:
\begin{enumerate}
\item[(a)] If a Levi subgroup of $G$ is
(stably) Cayley, is $G$ (stably) Cayley?
\item[(b)] Let $G$ be reductive. If $G/R$ is
(stably) Cayley, is $G$ (stably) Cayley? %%
\item[(c)] Let $G$ be reductive and let
$H_1,\ldots,H_n$ be a complete list of its
connected normal simple subgroups. What is the
relation between the (stably) Cayley property of
$G$ and that of $H_1,\ldots, H_n$?
\end{enumerate}

\subsection { Unipotent normal subgroups.}
\label{subsection4.1} We will need a
generalization of Example \ref{uni}. Let $U$
be a normal unipotent subgroup of $G$. Denote
by $\lu$ the Lie algebra of $U$. The group $G$
acts on $U$ by conjugation and on $\lu$ by
$\operatorname{Ad}_G\!|_{\lu}$.

\begin{lem} \label{lem.unipotent}
There exists a $G$-isomorphism of
$G$-varieties $U\to \lu$.
\end{lem}

\begin{proof}
We may assume without loss of generality that
$G \subset \GL_n$. Since
$\operatorname{Ad}^{}_G$ is given by
\eqref{adj}, it follows from \eqref{e.ln} that
$\operatorname{ln}\!: U\to \lu$ is a
$G$-morphism. By Example \ref{uni}, it is an
isomorphism, whence the claim. \quad $\square$
\renewcommand{\qed}{} \end{proof}

\subsection { The Levi decomposition.}
\label{subsection.levi} % Subsection 4.2

\begin{prop} \label{prop.levi}
Let $L$ be a Levi subgroup of $G$. %%
\begin{enumerate}
\item[\rm (a)] If $L$ is Cayley, then so is
$G$. %%
\item[\rm (b)] $G$ is stably Cayley if and
only if $L$ is stably Cayley.
\end{enumerate}
\end{prop}

\begin{proof} Let $T$ be
a maximal torus of $L$. It is a maximal torus
of $G$ as well, \cite[11.20]{borel}. Using the
notation of \eqref{torus-notations} and
Subsection~\ref{sect.cartan}, we have
$C=T\times U$ where $U$ is a unipotent group,
\cite[12.1]{borel}. Let  $\lu$ be the Lie
algebra of $U$ and let $d=\dim U$. As $T$ and
$U$ are, respectively, the semisimple and
unipotent parts of the nilpotent group $C$,
they are stable under the conjugating action
of $N$, and $C$, as an $N$-variety, is the
product of the $N$-varieties $T$ and $U$.
Consequently, $\lt$ and $\lu$ are stable under
the adjoint action of $N$, and $\lc$, as an
$N$-variety, is the product of $N$-varieties
$\lt$ and $\lu$. By Lemma~\ref{lem.unipotent},
there exists an isomorphism of $N$-varieties
\begin{equation}\label{tau}
\tau \colon U \lra\lu.
\end{equation}

(a) Assume that $L$ is Cayley. Then by
the Corollary of Lemma~\ref{cor.caley},
%%\ref{cor.torus},
there is a birational $W_{L,T}$-iso\-morphism
$\varphi \colon T \overset{\simeq}
{\dashrightarrow}\lt$. Since the action of
$W_{L,T}$ on $T$ (respectively, $\lt$) is
faithful, $W_{L,T}$ can be considered as a
transformation group of $T$ (respectively,
$\lt$). By \cite[11.20]{borel}, it coincides
with the transformation group $\{T\to T,\
t\mapsto ntn^{-1}\mid n\in N\}$ (respectively,
$\{\lt\to \lt,\ x\mapsto {\rm Ad}^{}_G
n(x)\mid n\in N\}$). Therefore the map
$\varphi$ is $N$-equivariant. Hence
\[ \varphi \times \tau \colon C = T \times U
{\dashrightarrow} \lt \oplus \lu = \lc
\] is a birational $N$-isomorphism.
Lemma~\ref{cor.caley} now implies that $G$ is
Cayley.

\smallskip

(b) Since $L\times \T_m^d$ is the Levi
subgroup of $G\times \T_m^d$, it follows from
(a) that if $L$ is stably Cayley, then $G$ is
stably Cayley.

To prove the converse, it suffices to show
that if $G$ is Cayley, then $L$ is stably
Cayley. In turn, Lemma \ref{cor.caley} and its
Corollary
%%\ref{cor.torus}
reduce this to proving that if there exists a
birational $N$-isomorphism
\begin{equation*}
\alpha \colon C = T \times U \overset{\simeq}
{\dashrightarrow}
%%\dasharrow
\lt \times \lu = \lc,
\end{equation*}
then $T$ and $\lt$ are stably birationally
$W_{L, T}$-isomorphic. We shall prove this
%%latter
last statement.

 Since $T$ is the identity component
 of $N_{L, T}=N\cap L$ and $T$ acts trivially on $C$
and $\lc$, the actions of $N_{L, T}$ on $C$,
$\lc$, $T$, $\lt$, $U$, and $\lu$ descend to
actions of $W_{L, T}=N_{L, T}/T$. Moreover,
$C$ (respectively, $\lc$), as a $W_{L,
T}$-variety, is the product of $W_{L,
T}$-varieties $T$ and $U$ (respectively, $\lt$
and $\lu$), and $\alpha$ is a birational
$W_{L, T}$-isomorphism.

Since $W_{L, T}$ acts linearly on $\lu$,
Lemma~\ref{lem.no-name}(b) implies that there
are birational $W_{L, T}$-isomorphisms
$$\beta\!: T\times \bbA^d
\overset{\simeq} {\dashrightarrow}
%%\dasharrow
T\times \lu \quad \text{and}\quad \gamma\!:
\lt\times \lu
%%\dasharrow
\overset{\simeq} {\dashrightarrow}\lt\times
\bbA^d \, ,$$ where $W_{L, T}$ acts on
$T\times \bbA^d$ and $\lt\times \bbA^d$ via
the first factors. Considering the composition
of the following birational $W_{L,
T}$-isomorphisms
$$
\xymatrix{ T\times \bbA^d\ar@{-->}[r]^\beta&
T\times \lu\ar[r]^{{\rm id}\times \tau^{-1}}&
T\times U\ar@{-->}[r]^\alpha& \lt\times
\lu\ar@{-->}[r]^\gamma& \lt\times \bbA^d},
$$
we now see that $T$ and $\lt$ are indeed
stably birationally $W_{L, T}$-isomorphic.
\quad $\square$
\renewcommand{\qed}{} \end{proof}

\begin{remark} The converse to Proposition
\ref{prop.levi}(a) fails for $G:={\bf
G}_2\times {\bf G}_a^2$. Indeed, the first
factor is the Levi subgroup of $G$. By
Proposition \ref{prop.iskovskih}, it is not
Cayley. Consider the group $H:={\bf G}_2\times
{\bf G}_m^2$. Both $G$ and $H$ have the same
Lie algebra $\g$. By Proposition
\ref{prop.g_2}, $H$ is Cayley; let $\lambda:
H\stackrel{\simeq}{\dasharrow} \g$ be a Cayley
map. Fix a birational isomorphism of group
varieties $\delta: {\bf
G}_a^2\stackrel{\simeq}{\dasharrow} {\bf
G}_m^2$. Since the second factors of  $G$ and
$H$ lie in the kernels of conjugating and
adjoint actions, $\lambda\circ({\rm id}\times
\delta)\!:G\stackrel{\simeq}{\dasharrow}\g$ is
a Cayley map. Thus $G$ is Cayley.
\end{remark}

\noindent{\bf Corollary.}\label{cor.solvable}
{\it Every connected solvable linear algebraic
group $G$ is Cayley.}

\begin{proof}
%%Let $G$ be a connected solvable group.
A Levi subgroup $L$ of $G$ is a
torus,\,\cite[10.6]{borel}. By
Example~~\ref{ex.abelian}, $L$ is Cayley.
Hence by Proposition~\ref{prop.levi}(a), $G$
is Cayley as well. \quad $\square$
\renewcommand{\qed}{} \end{proof}

\subsection { From reductive to
semisimple.} \label{red-to-semisimple} % Subsection 4.3
\begin{prop} \label{prop.ss}
Let $G$ be a connected reductive group and let
$Z$ be a connected closed central subgroup of
$G$.
\begin{enumerate}
\item[\rm(a)] If $G/Z$ is Cayley, then so is
$G$. \item[\rm(b)] $G$ is stably Cayley if and
only if $G/Z$ is stably Cayley.
\end{enumerate}
\end{prop}

\begin{proof} Since $G$ is reductive, $R$ is
a torus and the identity component of the
center of $G$; see \cite[11.21]{borel}. Thus
$Z$ is a subtorus of $R$. Let $T$ be a maximal
torus of $G$. We have $R \subset T$ (see
\cite[11.11]{borel}), $T/Z$ is a maximal torus
of $G/Z$ and the natural epimorphism $G\to
G/Z$ identifies $W$ with $W_{G/Z, T/Z}$ (we
use the notation of \eqref{torus-notations}
and Subsection \ref{sect.cartan}); see
\cite[11.20]{borel}. Since $Z$ is central, it
is pointwise fixed with respect to the action
of $W$. Thus we have the following exact
sequence of $W$-homomorphisms of tori:
\begin{equation*}
e\lra Z\lra T\lra T/Z\lra e
\end{equation*}
which in turn yields the exact sequence of
$W$-lattices of character groups
\[ 0 \lra \widehat{
T/Z} \lra \widehat{T} \lra \widehat{ Z} \lra
0.
\]
Note that $W$ acts trivially on $\widehat{Z}$.
In particular, $\widehat{Z}$ is a permutation
$W$-lattice, and the last exact sequence tells
us that the character lattices $\widehat{T}$
and $ \widehat{T/Z }$ are equivalent; see
Definition~\ref{equivalent}. Thus
% by Lemma \ref{stabiso} and its Corollary,
if one of them is quasi-permutation, then so
is the other. Part~(b) now follows from
Theorem~\ref{thm1}.

Since the $W$-fields $k(T)$ and $k(T/Z)$ are
$W$-isomorphic to $k(\widehat{T})$ and
$k(\widehat{T/Z})$, respectively, we deduce
from Lemma \ref{lem.no-name}(d) that $T$ is
birationally $W$-isomorphic to $T/Z\times
\bbA^m$, where $W$ acts on $T/Z\times \bbA^m$
via the first factor and $m=\dim Z$.

On the other hand, let $\mathfrak f$ and
$\mathfrak z$ be the Lie algebras of $T/Z$ and
$Z$, respectively. Then, since the Lie algebras
$\lt$ and $\mathfrak f\oplus \mathfrak z$ are
$W$-equivariantly isomorphic and $W$ acts on
$\mathfrak z$ trivially,
%%Since $W$
%%acts trivially on $\lr$,
we see that $\lt$, as a $W$-variety, is
isomorphic to $\mathfrak f\times \bbA^m$,
where $W$ acts on $\mathfrak f\times \bbA^m$
via the first factor.

Now to prove part (a), assume that $G/Z$ is
Cayley. Then by the Corollary of
Lemma~\ref{cor.caley},
%%\ref{cor.torus},
there is a birational $W$-isomorphism
$\varphi:
T/Z\stackrel{\simeq}{\dasharrow}\mathfrak f$.
This gives a birational $W$-isomorphism
$\xymatrix{T/Z \times \bbA^m \ar@{-->}[r]^{\
\hskip 1mm\varphi\times{\rm id}}& \mathfrak
f\times \bbA^m}$. Applying the Corollary of
Lemma~\ref{cor.caley} once again, we conclude
that $G$ is Cayley. This completes the proof
of part (a). \quad $\square$
\renewcommand{\qed}{} \end{proof}

%%\smallskip

%%\begin{cor}

Setting $Z=R$, we obtain

\medskip

\noindent{\bf Corollary.} \label{ss} {\it Let
$G$ be a connected reductive group and
$G\hskip -.3mm_{\rm s
%%\hskip -.3mm
s}:=G/R$.
\begin{enumerate}
\item[\rm(a)] If $G\hskip -.3mm_{\rm s
%%\hskip-.3mm
s}$ is Cayley, then so is $G$.
\item[\rm(b)] $G$ is stably Cayley if and only
if $G\hskip -.3mm_{\rm s
%%\hskip -.3mm
s}$ is
stably Cayley. \quad $\square$
\end{enumerate}
}

\begin{remark} The converse to part (a)
of the Corollary fails for $G={\bf G}_2\times
\T_m^2$. Indeed, $G$ is Cayley by Proposition
\ref{prop.g_2} and $G/R\simeq{\bf G}_2$ is not
Cayley by Proposition \ref{prop.iskovskih}.
\end{remark}

\subsection { From semisimple to simple.}
\label{semisimple-to-simple} % Subsection 4.4
Let $G_1, \ldots, G_n$ be connected linear
algebraic groups and let $\g_i$ be the Lie
algebra of $G_i$. If each $G_i$ is Cayley,
then so is $G_1\times\ldots\times G_n$; see
Example~\ref{product}. The converse fails for
$n=2$, $G_1={\bf G_2}$, $G_2={\bf G}_m^2$; see
Propositions~\ref{prop.iskovskih}
and~\ref{prop.g_2}.

\begin{lem} \label{lem.direct-product}
$G_1\times\ldots\times G_n$ is stably Cayley
if and only if each $G_i$ is stably Cayley.
\end{lem}

\begin{proof}
The ``if" direction follows from Definition
\ref{stably-cayley} and Example \ref{product}.
To prove the converse, we use the fact that
the underlying variety of each $G_i$ is
rational over $k$; see \cite{c1}. This implies
that the underlying variety of
$G_1\times\ldots \times G_n$, as a
$G_i$-variety, is birationally isomorphic to
$G_i \times \T_m^{d_i}$ with the  conjugating
action via the first factor and $d_i=\sum_{j
\neq i} \dim G_j$. The ``only if'' direction
now follows from Definition
\ref{stably-cayley} and the fact that the
underlying variety of the Lie algebra of
$G_1\times\ldots \times G_n$, as a
$G_i$-variety, is isomorphic to $\g_i\oplus
k^{d_i}$ with the adjoint action via the first
summand. \quad $\square$
\renewcommand{\qed}{} \end{proof}

As usual, given subgroups $X$ and $Y$ of $G$,
we denote by $(X, Y)$ the subgroup generated
by the commutators $xyx^{-1}y^{-1}$ with $x\in
X$, $y\in Y$.

\begin{prop}\label{decomposition}
Assume $G$ is a connected reductive group and
let $H_1,\ldots, H_m$ be the connected closed
normal
%%{\rm(}automatically
%%normal\,{\rm)}
subgroups of $G$ such that
\begin{enumerate}
\item[\rm(i)] $(H_i, H_j)=e$ for all $i\neq
j$, \item[\rm(ii)] $G=H_1\ldots H_m$.
\end{enumerate}
\noindent Let $\widetilde H_i$ be
%%{\rm``}\!complement\,{\rm''} to $H_i$ in
%%$G$, i.e.,
the subgroup of $G$ generated by all $H_j$'s
with $j\neq i$. If $G$ is stably Cayley, then
each $G/\widetilde H_i\simeq H_i/(H_i\cap
\widetilde H_i)$ is stably~Cayley.
\end{prop}
\begin{proof} Since
$H_1,\ldots, H_m$ are connected, each
$\widetilde H_i$ is connected; see
\cite[2.2]{borel}. Since $G$ is reductive, all
$H_i$ and $\widetilde H_i$ are reductive.

It follows from (i) and (ii) that
\begin{equation*}\label{epi}
H_1\times\ldots\times H_m\to G,\quad
(h_1,\ldots, h_m)\mapsto h_1\ldots h_m,
\end{equation*}
is an epimorphism of algebraic groups. Let
$T_i$ be a maximal torus of $H_i$. Then
$T_1\times\ldots\times T_m$ is a maximal torus
of $H_1\times\ldots\times H_m$. Therefore its
image $T:=T_1\ldots T_m$
 under the above
epimorphism is a maximal torus of $G$; see
\cite[11.14]{borel}. The same argument shows
that the group $S_i$ of $T$ generated by all
$T_j$'s with $j\neq i$ is a maximal torus of
$\widetilde H_i$.

It follows from (i) that $S_i$ is pointwise
fixed under the conjugating action of
$N_i:=N_{H_i, T_i}$ on $T$. This action
clearly descends to an action of $W_i:=W_{H_i,
T_i}=N_i/T_i$. Since $H_i$ is connected
reductive, any maximal torus of $H_i$
coincides with its centralizer in $H_i$; see
\cite[13.17]{borel}. Consequently, $T\cap
H_i=T_i$ and $W_i$, considered as a
transformation group of $T$, is the image of
$N_i$ under the natural projection $N\to
N/T=W$. The natural epimorphism $\pi_i: H_i\to
H_i/(H_i\cap \widetilde H_i)$ identifies $W_i$
with $W_{H_i/(H_i\cap \widetilde H_i),
\pi_i(T_i)}$, so that the isomorphism
$T_i/(T_i\cap \widetilde H_i)\to \pi_i(T_i)$
induced by $\pi_i$ is $W_i$-equivariant; cf.,
e.g., \cite[11.20, 11.11]{borel}.

The same argument applied to $\widetilde H_i$
and $S_i$ instead of $H_i$ and $T_i$ shows
that $T\cap \widetilde H_i=S_i$,
\begin{equation*}\label{ht} T_i\cap
\widetilde H_i=T_i\cap S_i,
\end{equation*}
and that a maximal torus of $H_i/(H_i\cap
\widetilde H_i)$ is $W_i$-isomorphic to
$T_i/(T_i\cap S_i)$. Now observe that
$T_i/(T_i\cap S_i)$ is $W_i$-isomorphic to
$T/S_i$ because $T=T_iS_i$. Therefore there is
an exact sequence of $W_i$-homomorphisms of
tori:
\begin{equation*}\label{exact-tori}
e \lra S_i\lra T\lra T_i/(T_i\cap S_i) \lra e
\, .
\end{equation*}
Passing to the character groups, we deduce
from it the following exact sequence of
$W_i$-lattices:
\begin{equation*}
0\lra\widehat{T_i/(T_i\cap S_i)}\lra
\widehat{T}\lra\widehat{S_i}\lra 0.
\end{equation*}
As the action of $W_i$ on $S_i$ is trivial,
$\widehat{S_i}$ is a trivial and, in
particular, a permutation $W_i$-lattice. Hence
the above exact sequence shows that
$\widehat{T_i/(T_i\cap S_i)}$ and
$\widehat{T}$ are equivalent $W_i$-lattices.

Assume now that $G$ is stably Cayley. Then
Theorem \ref{thm1} implies that $\widehat T$
is quasi-permuta\-ti\-on as a $W$-lattice, and
hence
%%$\widehat T$ is quasi-permutation
as a $W_i$-lattice because $W_i$ is a subgroup
of $W$. Therefore the equivalent $W_i$-lattice
$\widehat{T_i/(S_i\cap T_i)}$ is
quasi-permutation as well. Since the latter is
the character lattice of $H_i/(H_i\cap
\widetilde H_i)$, Theorem \ref{thm1} implies
that $H_i/(H_i\cap \widetilde H_i)$ is stably
Cayley.\quad $\square$
\renewcommand{\qed}{}
\end{proof}

\noindent{\bf Corollary.} {\it Let $G$ be a
connected semisimple group. Let $H_1,\ldots,
H_m$ be the minimal elements among its
connected closed normal subgroups. Define
$\widetilde H_i$ as in Proposition {\rm
\ref{decomposition}}. If $G$ is stably Cayley,
then each $H_i/(H_i\cap \widetilde H_i)$ is
stably Cayley.}
%%\end{cor}
\begin{proof}
By \cite[14.10]{borel}, the assumptions of
Proposition \ref{decomposition} hold.
 \quad $\square$
\renewcommand{\qed}{}
\end{proof}

\begin{remark}
In Proposition \ref{decomposition}, if $G$ is
stably Cayley, $H_i$ is not necessarily stably
Cayley. For example, take $G={\bf GL}_n$,
$m=2$, $H_1={\bf G}_m$ diagonally embedded in
${\bf GL}_n$ and $H_2={\bf SL}_n$. Then $G$ is
Cayley by Example \ref{ex.gl_n}, and $H_2$ is
not stably Cayley for $n>3$ by
Theorem~\ref{thm2}.
\end{remark}

\section{\bf Proof of Theorem~\ref{thm2}: An overview}
\label{stablycayley} % Section 5

In this section we outline a strategy for
proving Theorem~\ref{thm2}; the technical
parts of the proof will be carried out in
Sections~\ref{sect.intA_n1}--\ref{sect.intD_n}.

By Theorem~\ref{thm1}, it will suffice to
determine which connected simple groups have a
stably rational generic torus (or,
equivalently, a quasi-permutation character
lattice). {\sc Cortella} and {\sc
Ku\-nyavski\v\i} in \cite[Theorem 0.1]{ck}
have classified all  simply connected and all
adjoint connected simple groups that have a
quasi-permutation character lattice. These are
precisely $\SO_{2n+1}$, $\Symp_{2n}$,
$\PGL_n$, $\SL_3$, and ${\bf G}_2$. Therefore
in order to complete the proof of
Theorem~\ref{thm2}, we need to determine which
intermediate (i.e., neither simply connected
nor adjoint) connected simple groups have a
quasi-permutation character lattice.

Recall that intermediate connected simple
groups exist only for types $\A_n$ and $\D_n$.
Connected simple groups of type $\A_{n-1}$ are
precisely the groups $\SL_n/\m_d$, where $d$
is a divisor of~$n$. Among them, intermediate
groups are those with $1 < d < n$. In
Section~\ref{sect.intA_n2} we will prove the
following.

\begin{prop} \label{prop.A_n}
Let $d$ be a divisor of $n$, where $1 < d < n$
and $(n,d)\ne (4,2)$. Then the character
lattice of the group $\SL_n/\m_d$ is not
quasi-permutation.
\end{prop}

As we saw in Example~\ref{so}, the group
$\SL_4/ \m_2 $ is Cayley; in particular, by
Theorem~\ref{thm1}, its character lattice is
quasi-permutation.

The intermediate connected simple groups of
type $\D_n$ are $\SO_{2n}$ for any $n
\geqslant 3$ and the half-spinor groups
$\Spin_{2n}^{1/2}$ for even $n\geqslant 4$.
The latter are defined as follows. Consider
the spinor group $\Spin_{2n}$ for even
$n\geqslant 4$. Its
%%the
center is isomorphic to $\m_2\times \m_2$,
see\;\cite{c2}, \cite[\S25]{knus_et_al}, and
consequently contains precisely three
subgroups of order~$2$. One of them is the
kernel of the vector representation, so the
quotient of $\Spin_{2n}$ modulo it is
$\SO_{2n}$. Two others are the kernels of the
half-spinor representations of $\Spin_{2n}$.
They are mapped to each other by an outer
automorphism of $\Spin_{2n}$, so the images of
the half-spin representations are isomorphic
to the same group that is $\Spin_{2n}^{1/2}$.

By Example \ref{so}, the groups $\SO_{2n}$ are
Cayley. If $n=4$, the group of outer
automorphisms of $\Spin_{2n}$ is isomorphic to
${\rm S}_3$ (for $n>4$, it is isomorphic to
${\rm S}_2$) and acts transitively on the set
of all subgroups of order $2$ of the center of
$\Spin_{2n}$. Therefore
$\Spin_{8}^{1/2}\simeq\SO_8$, whence it is
Cayley. Thus we only need to consider the
half-spin groups $\Spin_{2n}^{1/2}$ for even
$n>4$. In Section~\ref{sect.intD_n} we will
prove the following.

\begin{prop} \label{prop.D_n}
The character lattice of the half-spinor group
$\Spin^{1/2}_{2n}$ for even $n>4$ is not
quasi-permutation.
\end{prop}

Thus in order to complete the proof of
Theorem~\ref{thm2}, we need to prove
Propositions~\ref{prop.A_n}
and~\ref{prop.D_n}. This will be done in the
next three sections.

\section{\bf The groups \boldmath$\SL_n/\m_d$
and their character lattices}
\label{sect.intA_n1} % Section 6

\subsection { Lattices \boldmath
$Q_n(d)$.} \label{epsilon} % Subsection 6.1
For any divisor $d$ of $n$, the Weyl group $W$
of the group $G= \SL_n/\m_d$ is isomorphic
to~the permutation group $\Sym_n$ of the set
of integers $\{1,\ldots, n\}$. The character
lattice $\mathcal X_G$ is described as
follows.

Let $\ve_1,\dots,\ve_n$ be the standard basis
for the permutation $\Sym_n$-lattice $\bbZ
[\Sym_n/\Sym_{n-1}]$ on which $\sigma\in
\Sym_n$ acts via
\begin{equation}\label{action} \text{ $\sigma
(\ve_i)=\ve_{\sigma (i)}$\quad for all
$i=1,\dots, n$.}
\end{equation}
%%n-1$.
We naturally embed $\bbZ [\Sym_n/\Sym_{n-1}]$
into the $\mathbb Q$-vector space $\bbZ
[\Sym_n/\Sym_{n-1}]\otimes _{\bbZ}\mathbb Q$
endowed with the Euclidean structure such that
$\ve_1,\dots,\ve_n$ is the orthonormal basis
and we naturally extend the action of $\Sym_n$
to this space.

  The root system of type $\A_{n-1}$
is the subset
$$
\A_{n-1}:=\{\ve_i-\ve_j\mid 1\leqslant i\ne
j\leqslant n\}
$$
of $\bbZ [\Sym_n/\Sym_{n-1}]\otimes
_{\bbZ}\mathbb Q$. The Weyl
 group
 $W(\A_{n-1})$ of $\A_{n-1}$
 is $\Sym_n$ acting  by
 \eqref{action}, and
the standard base of $\A_{n-1}$ is $\alpha_1,
\ldots, \alpha_{n-1}$, where
\begin{equation}\label{Aalpha} \alpha_i=\ve_i-\ve_{i+1}, \quad
i=1,\dots,n-1;
\end{equation}
see \cite{Bourbaki-Cox}. The kernel of
augmentation map
$$\textstyle
\bbZ [\Sym_n/\Sym_{n-1}]
\stackrel{\epsilon}{\lra}\bbZ,\quad
\sum_{i=1}^na_i\ve_i\mapsto \sum_{i=1}^na_i,$$
is the root $\Sym_n$-lattice $\mathbb
Z\A_{n-1}$ of $\A_{n-1}$,
\begin{equation}
\label{ZA} \textstyle \bbZ\A_{n-1}:=\bbZ
\alpha_1\oplus\ldots\oplus \bbZ\alpha_{n-1}
=\bigl\{\sum_{i=1}^na_i\ve_i \mid
\sum_{i=1}^na_i=0\bigr\}.
\end{equation}
%%$$
%%is the root lattice of $\A_{n-1}$.
The character lattice of $\SL_n/\m_d$ is
isomorphic
%%(see Subsection
%%\ref{sect.lattices})
to the following $\Sym_n$-lattice:
%%$Q_n(d)$ in $\bbZ
%%[\Sym_n/\Sym_{n-1}]\otimes _{\bbZ}\mathbb Q$:
\begin{equation} \label{e.q_n(d)}
\textstyle \text{$Q_n(d):=
%%Q_n
\mathbb Z\A_{n-1}+\bbZ d \vp_1$, where
$\vp_1=\ve_1-\frac{1}{n} \sum_{i=1}^n\ve_i$}.
\end{equation}
The vector $\vp_1$ is the first fundamental
dominant weight of the root system $\A_{n-1}$
with respect to the
%%standard
base $\alpha_1, \ldots,\alpha_{n-1}$.

Observe that the character lattice of
$\SL_n/\m_n=\PGL_n$ is the root
$\Sym_n$-lattice $Q_n(n)=\mathbb Z\A_{n-1}$,
the character lattice of $\SL_n/\m_1=\SL_n$ is
the weight $\Sym_n$-lattice $\Lambda_n$ of
type $\A_{n-1}$, and that the following
sequences of homomorphisms of
$\Sym_n$-lattices are exact:
\begin{eqnarray}\textstyle
\label{exact} &0\lra \mathbb Z\A_{n-1}\lra
Q_n\left(n/d
%%\frac{n}{d}
\right)
\lra \bbZ/d\bbZ\lra 0\label{Qnd}, \\
&0\lra Q_n(d)\lra \Lambda_n\lra \bbZ/d\bbZ\lra
0\label{Qd}.
\end{eqnarray}
%%{eqnarray}
Here $\bbZ/d \bbZ$ denotes the cyclic group of
order $d$ with trivial $\Sym_d$-action. Note
that
\begin{equation}\label{qddual}\textstyle
Q_n(d)^*\simeq Q_n(n/d).
%%\frac{n}{d}
\end{equation}
In this section we will prove a number of
preliminary results about the lattices
$Q_n(d)$. In the next section we will use
these results to prove
Proposition~\ref{prop.A_n}.

\subsection { Properties of
\boldmath $Q_n(d)$.} \label{subsection6.2} We
begin by recalling a simple lemma which
computes the cohomology ${\rm
H}^1(\Gamma,\mathbb Z\A_{n-1})$ for all
subgroups $\Gamma$ of $\Sym_n$. The first part
is extracted from \cite[Lemma 4.3]{ll}.

\begin{lem} \label{Ancoh}
For any subgroup $\Gamma$ of $\Sym_n$, we have
$$\textstyle {\rm H}^1(\Gamma,\mathbb Z
\A_{n-1})\simeq \bbZ/\sum_{\O}|{\O}|\bbZ,$$
where $\O$ runs over the orbits of $\Gamma$ in
$\{1,\ldots,n\}$. More explicitly, the
connecting homomorphism of the cohomology
sequence induced by the augmentation sequence
\begin{equation}\label{aug}
0\lra \mathbb Z\A_{n-1}\lra
\bbZ[\Sym_n/\Sym_{n-1}] \stackrel{\epsilon
}{\lra}\bbZ\lra 0
\end{equation}
is given by
\begin{eqnarray*}
\bbZ=\bbZ[\Sym_n/\Sym_{n-1}]/\mathbb
Z\A_{n-1}\stackrel{\partial}{\lra} {\rm
H}^1(\Gamma,\mathbb Z\A_{n-1}), \quad m \ve_1+
\mathbb Z\A_{n-1}\mapsto [\sigma \mapsto
m(\ve_{\sigma(1)} - \ve_1)],
\end{eqnarray*}
where the image is the class of the given
$1$-cocycle from $\Gamma$ to $\mathbb
Z\A_{n-1}$.
\end{lem}

\begin{proof}
From the cohomology sequence that is
associated with \eqref{aug}, one obtains the
exact sequence
$\bbZ[\Sym_n/\Sym_{n-1}]^{\Gamma}\stackrel{\epsilon}{\to}\bbZ
\stackrel{\partial}{\to} {\rm
H}^1(\Gamma,\mathbb Z\A_{n-1})\to 0$ which
implies the asserted description of ${\rm
H}^1(\Gamma,\mathbb Z\A_{n-1})$. The
calculation of the connecting homomorphism
$\partial$ follows directly from the
identification of $\bbZ$ with
$\bbZ[\Sym_n/\Sym_{n-1}]/\mathbb Z\A_{n-1}$
and an application of the Snake Lemma. \quad
$\square$
\renewcommand{\qed}{}
\end{proof}

\begin{lem} \label{connecthom}
For any subgroup $\Gamma$ of $\Sym_n$, the
exact sequence \eqref{Qnd} induces the
following connecting homomorphism in
cohomology:
\begin{equation*}\textstyle
%%\partial:
\bbZ/d\bbZ=Q_n(n/d
%%\frac{n}{d}
)/\mathbb Z\A_{n-1} \stackrel{\partial}{\lra}
{\rm H}^1(\Gamma,\mathbb Z\A_{n-1}), \quad
%%\\
m+d\bbZ\mapsto \frac{mn}{d}
%%mn/d
+\sum_{\O}|\O|\bbZ,
\end{equation*}
where the sum on the right runs over the
orbits $\O$ of $\Gamma$ in $\{1,\ldots,n\}$.
In particular, if $|{\rm H}^1(\Gamma,\mathbb
Z\A_{n-1})|$ divides
%%$\frac{n}{d}$,
$n/d$, then $\partial$ is the zero map.
\end{lem}

\begin{proof}
Since $Q_n(n/d
%%\frac{n}{d}
)$ has $\bbZ$-basis $\frac{n}{d}\vp_1,
\ve_1-\ve_2,\dots, \ve_{n-2}-\ve_{n-1}$ where
$\vp_1$ is given by \eqref{e.q_n(d)}, we
conclude that
%%then
$Q_n(n/d
%%\frac{n}{d}
)/\mathbb Z\A_{n-1}$ is generated by
$\frac{n}{d}\vp_1+\mathbb Z\A_{n-1}$. Using
the Snake Lemma, one sees that the connecting
homomorphism $$
%%\partial:
\bbZ/d\bbZ=Q_n(n/d
%%\frac{n}{d}
)/\mathbb Z\A_{n-1} \stackrel{\partial}
%%{\lra}
{\to} {\rm H}^1(\Gamma,\mathbb Z\A_{n-1})$$
sends $\frac{n}{d}\vp_1+\mathbb Z\A_{n-1}$ to
the class of the $1$-cocycle $[\sigma\mapsto
\frac{n}{d}(\ve_{\sigma (1)}-\ve_1)]$ in ${\rm
H}^1(\Gamma,\mathbb Z\A_{n-1})$. An
application of Lemma~\ref{Ancoh} and the
identification $\bbZ/d\bbZ = Q_n(n/d
%%\frac{n}{d}
)/\mathbb Z\A_{n-1}$ completes the proof of
the first statement. The second statement
follows directly from the first. \quad
$\square$
%%\qed
\renewcommand{\qed}{}\end{proof}

\begin{lem}\label{H1triv} Let $\Gamma$
be a subgroup of ${}\,\Sym_n$ which fixes at
least one integer $i\in \{1,\dots,n\}$. Then
${\rm H}^1(\Gamma,Q_n(d))=0$.
\end{lem}

\begin{proof}
Note that in this case, $\{\ve_t-\ve_i\mid
t\ne i\}$ is a permutation basis for $\mathbb
Z\A_{n-1}$ so that both $\mathbb Z\A_{n-1}$
and $\Lambda_n=(\mathbb Z\A_{n-1})^*$ are
permutation $\Gamma$-lattices. This implies
that ${\rm H}^1(\Gamma,\mathbb
Z\A_{n-1})=0={\rm H}^1(\Gamma,\Lambda_n)$.
Observe that
$\nu_i=\ve_i-\frac{1}{n}\sum_{t=1}^n\ve_t\in
\Lambda_n^\Gamma$ and that
$\nu_i+Q_n(d)=\vp_1+Q_n(d)$ since
$\nu_i-\vp_1=\ve_i-\ve_1\in \mathbb Z\A_{n-1}
\subseteq Q_n(d)$. Then applying cohomology to
the exact sequence \eqref{Qd}, we obtain
$$\Lambda_n^{\Gamma} \lra
\bbZ/d\bbZ\lra {\rm H}^1(\Gamma,Q_n(d))\lra
{\rm H}^1(\Gamma,\Lambda_n)=0. $$ Since
$\Lambda_n/Q_n(d)=\bbZ/d\bbZ$ is generated by
$\vp_1+Q_n(d)$, the above argument shows that
the map $\Lambda_n^{\Gamma} \to \bbZ/d\bbZ$ is
surjective so that ${\rm H}^1(\Gamma,
Q_n(d))=0$, as required. \quad $\square$
\renewcommand{\qed}{} \end{proof}

For a sequence of integers $1\leqslant
i_1<\ldots <i_r\leqslant n$, set
%%put
$$
\Sym_{\{i_1, \ldots, i_r\}}:=
\bigl\{\sigma\in\Sym_n\mid \sigma(j)=j\hskip
2mm \text {for every} \hskip 2mm j\notin
\{i_1, \ldots, i_r\}\bigr\}.
$$
This is a subgroup of $\Sym_n$; in particular,
$\Sym_{\{1,\ldots,n\}}=\Sym_n$. The map
$$\iota_{\{i_1,\ldots, i_r\}}:\
\Sym_r\lra \Sym_{\{i_1, \ldots, i_r\}},
%%\subset \Sym_n,
\quad \iota_{\{i_1,\ldots, i_r\}}(\sigma
)(i_s)=i_{\sigma(s)}
%%\sigma\mapsto \big(i_1\rightarrow i_{\sigma(1)},\ldots,
%%i_d\rightarrow i_{\sigma(d)}\big)
\quad \text{for all $\sigma$
%%\in \Sym_d,\
and $s$,}
%%=1,\ldots, d,
$$ is an isomorphism. In what follows, the
subgroup $\Sym_{\{1,\ldots, m\}}\times
\Sym_{\{m+1,\ldots, 2m\}}$ of $\Sym_{2m}$ is
denoted simply by $\Sym_m\times \Sym_m$.
%%a group embedding whose image is
%%$\Sym_{\{i_1, \ldots, i_d\}}$.
For a sequence of integers
$$1\leqslant i_1<\ldots <i_r<j_1
<\ldots <j_r<\ldots <l_1<\ldots <l_r\leqslant
n,$$ the image of the
%%homomorphism
embedding
$$
\Sym_r\lra \Sym_n,
%%{\{i_1,\ldots, i_d, j_1,\ldots, j_d,\ldots,
%%l_1,\ldots, l_d\}},
\quad \sigma\mapsto \iota_{\{i_1,\ldots,
i_r\}}(\sigma) \iota_{\{j_1,\ldots,
j_r\}}(\sigma)\ldots \iota_{\{l_1,\ldots,
l_r\}}(\sigma),
$$
is called the {\it copy of $\Sym_r$ diagonally
embedded  in} $\Sym_{\{i_1,\ldots, i_r,
j_1,\ldots, j_r,\ldots, l_1,\ldots, l_r\}}$.

\begin{lem}\label{anrest} Let $n=td$. Then the following
properties hold:
%%\begin{enumerate}
%%\item[(a)]

{\rm (a)} Let $
%%H
X_d$
%%\simeq \Sym_d$
be the copy of $\Sym_d$ diagonally embedded in
$\Sym_n$. Then
\begin{equation*}
\label{Hd} \mathbb Z\A_{n-1}\vert_{X_d}\simeq
\mathbb Z\A_{d-1} \oplus
\bbZ[\Sym_d/\Sym_{d-1}]^{t-1}.\end{equation*}

{\rm (b)} Let $Y_d:= \Sym_{\{1,\dots,d\}}
\times \widetilde X_d$ where $\widetilde X_d$
is the copy of $\,\Sym_d$ diagonally embedded
in $\Sym_{\{d+1,\dots,n\}}$. Then
\begin{equation*}\label{H2d}
\mathbb Z\A_{n-1}\vert_{Y_d} \simeq\mathbb
Z\A_{2d-1}\vert_{\Sym_d\times \Sym_d} \oplus
\bbZ[(\Sym_d\times \Sym_d)/(\Sym_d\times
\Sym_{d-1})]^{t-2}.
\end{equation*}
\end{lem}

%%\smallskip

\begin{proof}
For the first statement, note that
$${\mathcal B}_1=
\{\ve_1-\ve_2,\dots,\ve_{d-1}-\ve_d\}\cup
\{\ve_i-\ve_{d+i} \mid i=1,\dots,(t-1)d\}$$ is
a basis for $\mathbb Z\A_{n-1}$, since
${\mathcal B}_0=\{\alpha_i=\ve_i-\ve_{i+1}\mid
i=1,\dots,n-1\}$ is a basis for $\mathbb
Z\A_{n-1}$, and the equations
$$\textstyle
\ve_i-\ve_{d+i}=\sum_{t=i}^{d+i-1} \alpha_k$$
for $i=1,\dots,(t-1)d$ show that the change of
coordinates matrix relating ${\mathcal B}_1$
to ${\mathcal B}_0$ is
%%an
upper triangular
%%matrix
with coefficients in $\mathbb Z$ and
diagonal entries $1$.
%%element of
%%$\U_n(\bbZ)$.
But then
\begin{eqnarray*}\textstyle
{\mathbb
Z\A_{n-1}}\vert_{X_d}&=&\bigoplus_{i=1}^{d-1}\bbZ(\ve_i-\ve_{i+1})\oplus
\bigoplus_{r=1}^{t-1}
\bigl(\bigoplus_{i=(r-1)d+1}^{rd}
\bbZ(\ve_i-\ve_{d+i})\bigr)\\
 &\simeq& \mathbb
Z\A_{d-1}\oplus \mathbb
Z[\Sym_d/\Sym_{d-1}]^{t-1} .
\end{eqnarray*}

For the second statement, similarly note that
$$\{\ve_1-\ve_2,\dots,\ve_{2d-1}-\ve_{2d}\}\cup
\{\ve_{i}-\ve_{d+i} \mid
i=d+1,\dots,(t-1)d\}$$ is a basis for $\mathbb
Z\A_{n-1}$ so that
%%$$
\begin{gather*}\textstyle
{\mathbb
Z\A_{n-1}}\vert_{Y_d}=\bigoplus_{i=1}^{2d-1}
\bbZ(\ve_i-\ve_{i+1})\oplus
\bigoplus_{r=2}^{t-1}
\bigl(\bigoplus_{i=(r-1)d+1}^{rd}
\bbZ(\ve_i-\ve_{d+i})\bigr)
\\
%%$$
%%$$
\simeq {\mathbb Z\A_{2d-1}}
\vert_{\Sym_d\times \Sym_d}\oplus
\bbZ[(\Sym_d\times \Sym_d)/(\Sym_d\times
\Sym_{d-1})]^{t-2}. \quad \qed
\end{gather*}
\renewcommand{\qed}{}\end{proof}

\section{ \bf Stably Cayley groups of
type \boldmath$\A_n$ } \label{sect.intA_n2}

\subsection {\bf Restricting \boldmath
$Q_n(d)$ to some subgroups.}
\label{subsection7.1} In this section we will
prove Proposition~\ref{prop.A_n}. We will
first show that $Q_n(d)$ restricted to certain
appropriate subgroups of $\,\Sym_n$ is
equivalent in each case to a smaller more
ma\-nageable sub\-lattice. We will then show
that the smaller lattices are not
quasi-permutation.

%Strategy moved to after end of document as it has changed.

\begin{prop} \label{qndrest} Suppose $d|n$ and let
$p$ be a prime divisor of
%%$\frac{n}{d}$,
$n/d$.
%%, $n=lp$.
Let $X_p$ be the copy of $\,\Sym_p$ diagonally
embedded
%%, diagonally embedded
in $\,\Sym_n$, and let
$Y_{p}=\Sym_{\{1,\dots,p\}}
%%_{1,\dots,p}
\times \widetilde X_p $,
%%\simeq \Sym_p\times \Sym_p$
where $\widetilde X_p$ is the copy of $\Sym_p$
dia\-gonally embedded
%%diagonally embedded
in $\Sym_{\{p+1,\dots,n\}}$.
%%_{p+1,\dots,n}$.
Then the following equivalencies hold:
\begin{enumerate}
\item[(a)]$Q_n(d)\vert_{X_p} \sim \Lambda_p$.
\item[(b)]$Q_n(d)\vert_{Y_{p}}\sim
{\Lambda_{2p}}\vert_{\Sym_p\times \Sym_p}$.
\end{enumerate}
\end{prop}

\begin{proof}
Recall that we have the exact sequence
\eqref{exact}. The definition of $p$ implies
that $n=lp$ for a positive integer $l$.  By
Lemma~\ref{anrest},
\begin{gather*}
{\mathbb Z\A_{n-1}}\vert_{X_p} \simeq\mathbb
Z\A_{p-1}\oplus
\bbZ[\Sym_p/\Sym_{p-1}]^{l-1},\\
{\mathbb Z\A_{n-1}}\vert_{Y_{p}}\simeq
{\mathbb Z\A_{2p-1}}\vert_{\Sym_p\times
\Sym_p}\oplus \bbZ[(\Sym_p\times
\Sym_p)/(\Sym_p\times \Sym_{p-1})]^{l-2}.
\end{gather*}

Using this and Lemma~\ref{Ancoh}, we see that
${\rm H}^1(\Gamma,\mathbb Z\A_{n-1})={\rm
H}^1(\Gamma,\mathbb Z\A_{p-1})=0$ or
$\bbZ/p\bbZ$ for all subgroups $\Gamma$ of
$X_p$ and that ${\rm H}^1(\Gamma,\mathbb
Z\A_{n-1})={\rm H}^1(\Gamma,\mathbb
Z\A_{2p-1})=0$ or $\bbZ/p\bbZ$ for all
subgroups $\Gamma$ of $Y_{p}$. Then 
Lemma~\ref{connecthom} and the fact that $p$
divides $n/d$ show that the connecting
homomorphism $\bbZ/d\bbZ \to {\rm
H}^1(\Gamma,\mathbb Z\A_{n-1})$ is zero for
all subgroups $\Gamma$ of $X_p$ or of $Y_{p}$.
But then the sequence \eqref{exact} restricted
to $X_p$ or $Y_{p}$ satisfies the conditions
of Proposition~\ref{equivlat}(b). This means
that
\begin{gather*}
{Q_{n}(d)}\vert_{X_{p}}={Q_{n}(n/d
%%\frac{n}{d}
)^*}\vert_{X_p} \sim {(\mathbb
Z\A_{n-1})^*}\vert_{X_p}\sim (\mathbb
Z\A_{p-1})^* = \Lambda_p,
\\
{Q_{n}(d)}\vert_{Y_{p}}={Q_{n}(n/d
)^*}\vert_{Y_{p}} \sim {(\mathbb
Z\A_{n-1})^*}\vert_{Y_{p}}\sim {(\mathbb
Z\A_{2p-1})^*}\vert_{\Sym_p\times \Sym_p} =
{\Lambda_{2p}}\vert_{\Sym_p\times \Sym_p}.
\hfill\qed
\end{gather*}
\renewcommand{\qed}{}
\end{proof}

\medskip

\subsection { Lattices
\boldmath$\Lambda_p$ and
\boldmath$\Lambda_{2p}$.}
\label{subsection7.2} The following lemma is
essentially a rephrasing of a result proved by
{\sc Bessenrodt} and {\sc Le Bruyn}
in~\cite{bl}:

\medskip

\begin{lem}~\label{lamp} Let $p>3$ be prime. Then
$\Lambda_p$ is not a quasi-permutation
$\Sym_p$-lattice.
\end{lem}

\begin{proof}
Tensoring the augmentation sequence for
%%$\Sym_n$
$\bbZ[\Sym_n/\Sym_{n-1}]$ with $\mathbb
Z\A_{n-1}$, we obtain the exact sequence
\begin{equation}\label{an2}
0\lra (\mathbb Z\A_{n-1})^{\otimes 2} \lra
\mathbb Z\A_{n-1}\otimes
\bbZ[\Sym_n/\Sym_{n-1}] \stackrel{
%%\chi
\tau}{\lra} \mathbb Z\A_{n-1}\lra 0.
\end{equation}
%%Here
We have
$$\mathbb Z\A_{n-1}\otimes \bbZ[\Sym_n/\Sym_{n-1}]
\simeq \bbZ[\Sym_n/\Sym_{n-2}].$$ One can show
that $\{
 %%y_{ij}\equiv
 (\ve_i-\ve_j)\otimes \ve_j \mid
 i\ne j\}$ is the set of elements of a permutation
basis for $\mathbb Z\A_{n-1}\otimes
\bbZ[\Sym_n/\Sym_{n-1}]$.
%% is .
 The map $
 %%\chi
 \tau$ then sends $(\ve_i-\ve_j)\otimes
 \ve_j$
 %%$y_{ij}$
 %%\equiv (\ve_i-\ve_j)\otimes \ve_j$
 to $\ve_i-\ve_j$.

For $p$ prime, {\sc Bessenrodt} and {\sc Le
Bruyn} in \cite{bl} show that
$$0\lra (\mathbb Z\A_{p-1})^{\otimes 2}\lra
\bbZ[\Sym_p/\Sym_{p-2}]\lra \mathbb
Z\A_{p-1}\lra 0$$ is a coflasque resolution of
$\mathbb Z\A_{p-1}$ as an $\Sym_p$-lattice.
They also show that $(\mathbb
Z\A_{p-1})^{\otimes 2}$ is permutation
projective as an $\Sym_p$-lattice but is only
$\Sym_p$-stably permutation if $p=2,3$. By
duality, the stable equivalence class of
%%$\rho(\Lambda_p)=
$\bigl((\mathbb Z\A_{p-1})^{\otimes
2}\bigr)^*$ is
%%a representative of
$\rho(\Lambda_p)$; see Subsection
\ref{stableeq}). The statements above then
imply that $\Lambda_p$ is not a
quasi-permutation $\Sym_p$-lattice for any
$p>3$. \quad $\square$
\renewcommand{\qed}{} \end{proof}

\medskip

\begin{prop}~\label{lambda2p}
Let $p$ be a prime and let
$$\Gamma := \langle (1,\dots,p),
(p+1,\dots,2p)\rangle\leqslant \Sym_p\times
\Sym_p\leqslant \Sym_{2p}.$$ Then the following hold.
\begin{enumerate}
\item[(a)] $\Sha^2(\Gamma, \Lambda_{2p})=0$.
In particular, a lattice in the stable
equivalence class $\rho(\Lambda_{2p})$ is
coflasque as an $\Gamma$-lattice. \item[(b)]
If $p$ is odd, $\Lambda_{2p}$ is not
quasi-permutation as an $\Gamma$-lattice and
hence is not quasi-permutation as an
$\Sym_p\times \Sym_p$-lattice.
%%On the other hand,
%%$\Lambda_{4}$ is quasi-permutation as an
%%$\Sym_2\times \Sym_2$-lattice.
\end{enumerate}
\end{prop}

\begin{proof}
(a) The second statement follows from the
first. Note that any proper subgroup of
$\Gamma$ is cyclic, so that by the claim
$\Sha^2(S,\Lambda_{2p})=0$ for all subgroups
$S$ of $\Gamma$. Then if
$$0\lra \Lambda_{2p}\lra Q\lra M\lra 0$$
is a flasque resolution of $\Lambda_{2p}$
considered as an $S$-lattice,
 then
${\rm H}^1(S, M) =\Sha^2(S,\Lambda_{2p})\break =0$ by
Lemma~\ref{shalem}.

To prove the first statement, we need to first
compute ${\rm H}^1(\Gamma, \Lambda_{2p})$ and
${\rm H}^2(\Gamma,\Lambda_{2p})$.

We have ${\rm H}^1(\Gamma,\Lambda_{2p})={\rm
H}^{-1}(\Gamma,\mathbb Z\A_{2p-1})$ by
duality. Then
$$
{\rm H}^{-1}(\Gamma,\mathbb Z\A_{2p-1})=
\Ker_{\mathbb Z{\mbox{\tt
\fontsize{9pt}{0mm}\selectfont
A}}_{2p-1}}(N_{\Gamma})/I_{\Gamma}\mathbb
Z\A_{2p-1},
$$
where $N_{\Gamma}$ is the endomorphism
$l\mapsto\sum_{a\in \Gamma} al$, and
$I_{\Gamma}$ is the augmentation ideal of
$\bbZ [\Gamma]$ (\cite{br}). We need to compute
$N_{\Gamma}$ on a basis for $\mathbb
Z\A_{2p-1}$: we have
$N_{\Gamma}(\ve_i-\ve_{i+1})=0$ for
$i=1,\ldots, p-1,
%%\linebreak
p+1,\dots,2p-1$, and
$N_{\Gamma}(\ve_{p}-\ve_{p+1})=p(\ve_1+\dots +
\ve_p-\ve_{p+1}-\dots-\ve_{2p})$. Then
$$\Ker N_{\Gamma} = \Span\{\ve_1-\ve_2,\dots,
\ve_{p-1}-\ve_p,\ve_{p+1}-\ve_{p+2},\dots,\ve_{2p-1}-
\ve_{2p}\}.$$ But $I_{\Gamma} \mathbb
Z\A_{2p-1}=\Ker N_{\Gamma}$ as $((1,\ldots,
p)-{\rm
id})(\ve_{p+1}-\ve_i)=\ve_i-\ve_{i+1}$,
$i=1,\dots,p-1$, $((p+1,\ldots, 2p)-{\rm
id})(\ve_1-\ve_i)=\ve_i-\ve_{i+1}$,
$i=p+1,\dots,2p-1$. This shows that ${\rm
H}^1(\Gamma,\Lambda_{2p})={\rm
H}^{-1}(\Gamma,\mathbb Z\A_{2p-1})=0$.

To determine ${\rm H}^2(\Gamma,\Lambda_{2p})$,
we use the restriction of the sequence
$$
0\lra \bbZ\lra \bbZ[\Sym_{2p}/\Sym_{2p-1}]\lra
\Lambda_{2p}\lra 0
$$
to $\Gamma$.  Let
\begin{equation}\label{hh}
%%H_1
C_1=\langle (1,\dots,p)\rangle, \hskip 2mm
%%H_2
C_2=\langle (p+1,\dots,2p)\rangle \quad
\text{and} \quad P_1=\bbZ[\Gamma/C_2], \hskip
2mm P_2=\bbZ[\Gamma/C_1].
\end{equation}
Then we have the following exact sequence of
$\Gamma$-lattices:
$$
0\lra \bbZ\lra P_1\oplus P_2\lra
\Lambda_{2p}\lra 0.
$$
Taking cohomology of this sequence, we get
\begin{multline*}
\ \hskip 10mm 0={\rm
H}^1(\Gamma,\Lambda_{2p})\lra {\rm
H}^2(\Gamma,\bbZ)\lra {\rm H}^2(\Gamma,P_1)
\oplus {\rm H}^2(\Gamma,P_2)\\
\lra {\rm H}^2(\Gamma,\Lambda_{2p})\lra {\rm
H}^3(\Gamma, \bbZ)\lra {\rm H}^3(\Gamma,
P_1\oplus P_2).\hskip 10mm \
\end{multline*}
\noindent But by Shapiro's Lemma, we have
${\rm H}^2(\Gamma, P_i)={\rm H}^2(\bbZ/p\bbZ,
\bbZ) = \bbZ/p\bbZ $ and ${\rm H}^3(\Gamma,
P_i)\break ={\rm H}^3(\bbZ/p\bbZ,\bbZ) =0$ for $i=1,
2$. Also, by the K\"{u}nneth
formula,~\cite[p.\,166]{weibel},
\begin{multline*}
\textstyle \ \hskip 10mm {\rm
H}^n(\Gamma,\bbZ)\!=\!\bigoplus_{i+j=n} {\rm
H}^i(\bbZ/p\bbZ,\bbZ)\otimes {\rm
H}^j(\bbZ/p\bbZ,\bbZ)\\
\textstyle \oplus \bigoplus_{i+j=n+1}
\!\Tor^1_{\bbZ}\bigl({\rm
H}^i(\bbZ/p\bbZ,\bbZ), {\rm
H}^j(\bbZ/p\bbZ,\bbZ)\bigr),\ \hskip 10mm
\end{multline*}
so that, in particular,
 ${\rm H}^3(\Gamma,\bbZ)=\bbZ/p\bbZ$ and
 ${\rm H}^2(\Gamma,\bbZ)=(\bbZ/p\bbZ)^2$.
This all yields an
%%shows that the
exact sequence
$$0\lra (\bbZ/p\bbZ)^2\lra
(\bbZ/p\bbZ)^2\lra {\rm
H}^2(\Gamma,\Lambda_{2p})\lra \bbZ/p\bbZ\lra
0,$$
%%is exact
and so ${\rm
H}^2(\Gamma,\Lambda_{2p})=\bbZ/p\bbZ$.

To show that $\Sha^2(\Gamma,\Lambda_{2p})=0$,
it would suffice to find a cyclic subgroup $C$
of $\Gamma$ for which $\Res^{\Gamma}_C: {\rm
H}^2(\Gamma,\Lambda_{2p}) \to {\rm
H}^2(C,\Lambda_{2p})$ is injective.

Take $C=
%%H_1
C_1$.  Since ${\rm
H}^1(\Gamma,\Lambda_{2p})=0$, we have that the
sequence
$$
0\lra {\rm
H}^2(\Gamma/C,\Lambda_{2p}^C)\stackrel{\Inf}{\lra}
{\rm H}^2(\Gamma,\Lambda_{2p})
\stackrel{\Res}{\lra} {\rm
H}^2(C,\Lambda_{2p})
$$
is exact. So it suffices to show that ${\rm
H}^2(\Gamma/C,\Lambda_{2p}^C)=0$.

The fundamental dominant weights for
$\Lambda_{2p}$ are
$$\textstyle
\vp_t=\sum_{i=1}^t\ve_i -
\frac{t}{2p}\sum_{i=1}^{2p}\ve_i, \quad
t=1,\ldots,2p-1.
$$
Let
$\nu_i=\ve_i-\frac{1}{2p}\sum_{i=1}^{2p}\ve_i$,
$i=1,\dots,2p$. Note that
$$
\nu_1=\vp_1,\hskip 2.5mm
\nu_t=\vp_t-\vp_{t-1}, t=2,\dots,2p-1,\hskip
2.5mm \nu_{2p}=-\vp_{2p-1}.
$$
This shows that
$\nu_1,\dots,\nu_p,\vp_{p+1},\dots,\vp_{2p-1}$
is another basis for $\Lambda_{2p}$ and that
$$\textstyle \Lambda_{2p}\vert_C=\bigoplus_{i=1}^p\bbZ \nu_i
\oplus \bigoplus_{i=p+1}^{2p-1}\bbZ
\vp_i\simeq \bbZ [C]\oplus \bbZ^{p-1}.$$

This shows that
$$\textstyle
\Lambda_{2p}^C=\bbZ (\sum_{i=1}^p\nu_i)\oplus
\bigoplus_{i=p+1}^{2p-1}\bbZ \vp_i
=\bigoplus_{i=p}^{2p-1}\bbZ\vp_i=
\bigoplus_{i=p+1}^{2p}\bbZ \nu_i.$$ But
$\Gamma/C$ permutes $\nu_{p+1},\ldots,
\nu_{2p}$ cyclically so that
$\Lambda_{2p}^C\simeq \bbZ [\Gamma/C]$. This
implies that ${\rm
H}^2(\Gamma/C,\Lambda_{2p}^C)=0$ as required.

\smallskip

(b) To prove that $\Lambda_{2p}$ is not
$\Gamma$-quasi-permutation, we will construct
a coflasque $\Gamma$-resolution of $\mathbb
Z\A_{2p-1}$. By duality, this will give us  a
flasque resolution of $\Lambda_{2p}$.  We will
then show that the lattice in the stable
equivalence class  $\rho(\Lambda_{2p})$ is not
permutation projective as a $\Gamma$-lattice.

As $\alpha_1,\ldots,\alpha_{p-1}$ and
$\alpha_{p+1},\ldots, \alpha_{2p-1}$ are the
standard bases of the root subsystems of type
$\A_{p-1}$, we denote the $\Gamma$-sublattice
of $\bbZ\A_{2p-1}$ generated by them simply by
$\bbZ\A_{p-1}\oplus \bbZ\A_{p-1}$. Let $\iota$
be its natural embedding into $\bbZ\A_{2p-1}$.
%% be the natural embedding.
It is easily seen that
$\alpha_p+\bbZ\A_{p-1}\oplus \bbZ\A_{p-1}$ is
$\Gamma$-stable. This implies that there is an
exact sequence of $\Gamma$-lattices
$$0\lra \mathbb Z\A_{p-1}
\oplus \mathbb Z\A_{p-1}\stackrel{
\iota}{\lra} \mathbb Z\A_{2p-1}\lra \bbZ\lra
0.$$ A coflasque resolution of the
$\Gamma$-lattice $\mathbb Z\A_{p-1}\oplus
\mathbb Z\A_{p-1}$ is given by
$$0\lra \bbZ^2\lra P_1\oplus P_2 \lra \mathbb Z\A_{p-1}
\oplus \mathbb Z\A_{p-1}\lra 0$$ where $P_1$
and $P_2$ are defined by \eqref{hh} and the
generator of the $\Gamma$-lattice $P_1$
(respectively, $P_2$) is sent to $\alpha_1$
(respectively, $\alpha_{p+1}$).

We now extend $\iota$ to a coflasque
resolution of the $\Gamma$-lattice $\mathbb
Z\A_{2p-1}$. Let
$$
P_1\oplus P_2\oplus \bbZ [\Gamma] \oplus
\bbZ\stackrel{ \pi}{\lra} \mathbb Z\A_{2p-1}$$
be a map of $\Gamma$-lattices where $
\pi{}_{P_1\oplus P_2}= \iota$, $\pi $ sends
$1\in \bbZ[\Gamma]$ to $\alpha_p$, and $ \pi$
sends the $1\in \bbZ$ to
$\sum_{i=1}^p\ve_i-\sum_{i=p+1}^{2p}\ve_i=2\vp_p$.
It is easily verified that $ \pi $ is
surjective (in fact $\pi \vert_{\bbZ[\Gamma]}$
is surjective).

Let $L=\Ker \pi $. To check that $L$ is
coflasque and hence that $$0\lra L\lra
P_1\oplus P_2 \oplus \bbZ [\Gamma]\oplus
\bbZ\stackrel{\pi }{\lra} \mathbb
Z\A_{2p-1}\lra 0$$ is a coflasque resolution
of $\mathbb Z\A_{2p-1}$, we need only verify
that for $R:=P_1\oplus P_2\oplus \bbZ
[\Gamma]\oplus \bbZ$, we have $
\pi(R^K)=(\mathbb Z\A_{2p-1})^K$ for all
subgroups $K$ of $\Gamma$.

For $K= \Gamma$ or a cyclic subgroup generated
by a disjoint product of two $p$-cycles,
$(\mathbb Z\A_{2p-1})^K=\bbZ 2\vp_p$ so that $
\pi(\bbZ^K)= \pi(\bbZ)= (\mathbb
Z\A_{2p-1})^K$ and so $ \pi(R^K)=(\mathbb
Z\A_{2p-1})^K$.

The only other subgroups are $C_1$ and $ C_2$.
As the arguments are similar, we just consider
$C_1$: the lattice
$(\mathbb{Z}\A_{2p-1})^{C_1}$ has basis
$2\vp_p,\alpha_{p+1},\dots,\alpha_{2p-1}$, and
we have $ \pi(\bbZ)=\bbZ 2\vp_2$ and
$\pi(P_2^{C_1})= \pi(P_2)
=\bigoplus_{i=p+1}^{2p-1}\bbZ\alpha_{i}$. This
shows that
$$
0\lra L\lra P_1\oplus P_2\oplus \bbZ
[\Gamma]\oplus \bbZ\stackrel{
\pi}{\lra}\mathbb Z\A_{2p-1}\lra 0
$$
is a coflasque resolution. Dualizing, we
obtain a flasque resolution for
$\Lambda_{2p}$:
$$
0\lra \Lambda_{2p}\lra P_1\oplus P_2\oplus
\bbZ [\Gamma] \oplus \bbZ\lra L^*\lra 0.
$$
We have ${\rm H}^1(\Gamma,L^*) =
\Sha^2(\Gamma,\Lambda_{2p}) = 0$. This shows
that $L$ is flasque and coflasque as a
$\Gamma$-lattice.

We have the following commutative diagram with
exact rows and columns:
%%$$
\begin{gather}\label{diagra}
\begin{gathered}
\xymatrix@C=5mm@R=4mm{
&0\ar[d]&0\ar[d]&0\ar[d]&\\
0\ar[r]&\bbZ_{}^2 \ar[r]\ar[d] &P_1\oplus
P_2\ar[r]^{\hskip -6.5mm \iota}\ar[d] &\mathbb
Z\A_{p-1}\oplus \mathbb Z\A_{p-1}
\ar[r]\ar[d]& 0\\
0\ar[r] &L\ar[r]\ar[d] &P_1\oplus P_2\oplus
\bbZ [\Gamma]\oplus \bbZ \ar[d]\ar[r]^{\
\hskip 9mm \pi} &
\mathbb Z\A_{2p-1} \ar[r]\ar[d] &0\\
0\ar[r]& U(p)\ar[r]\ar[d] &\bbZ [\Gamma]\oplus
\bbZ \ar[r]^{ \hskip 3mm
\theta} \ar[d]&\bbZ \ar[r]\ar[d] &0\\
&0&0&0&}
\end{gathered}
\end{gather}
where $U(p)$ is the kernel of the induced map
$\theta$. Now $2\vp_p=\sum_{i=1}^{p-1}
i(\alpha_i+\alpha_{2p-i})+p\alpha_p$. So
$\theta$ sends  $1\in \bbZ[\Gamma]$ to
$\overline{\alpha_p}$ and sends $1\in \bbZ$ to
$p\overline{\alpha_p}$. This shows that
$$\{(h-1,0)\mid h\in \Gamma\}\cup \{(-p,1)\}$$
is  a $\bbZ$-basis for
$U(p)$. Note that $U(p)$ also satisfies
$$0\lra U(p)\lra \bbZ[\Gamma]\lra \bbZ/p\bbZ,$$
so that $\bbQ U(p)\simeq \bbQ [\Gamma]$.

  From the above diagram, we then see that $\bbQ L \simeq
\bbQ [\Gamma]\oplus \bbQ^2$. By~\cite[Lemmas
2 and  3]{cw}, to determine whether or not $L$ is
permutation projective is equivalent to
checking
 whether $\bbF_p L$
 is a permutation module
 for $\bbF_p [\Gamma]$.

Tensoring the diagram \eqref{diagra} with
$\bbF_p$ leaves it exact so we have the
following commutative diagram with exact rows
and columns:
\begin{gather*}\label{dia}
\begin{gathered}
\xymatrix@R=4mm@C=5mm{&0\ar[d]&0\ar[d]&0\ar[d]&\\
0\ar[r]&\bbF_p^2\ar[r]\ar[d] &\bbF_pP_1\oplus
\bbF_pP_2\ar[r]^{\ \hskip -8mm {\rm
id}\otimes\iota}\ar[d] &\bbF_p\A_{p-1} \oplus
\bbF_p\A_{p-1}
\ar[r]\ar[d]& 0\\
0\ar[r] &\bbF_p L\ar[r]\ar[d] &\bbF_pP_1\oplus
\bbF_pP_2\oplus \bbF_p [\Gamma]\oplus \bbF_p
\ar[d]\ar[r]^{\ \hskip 14mm {\rm
id}\otimes\pi} & \bbF_p\A_{2p-1}
\ar[r]\ar[d] &0\\
0\ar[r]& \bbF_p U(p)\ar[r]\ar[d] &\bbF_p
[\Gamma]\oplus \bbF_p \ar[r]^{\ \hskip 4mm
{\rm id}\otimes\theta}\ar[d] &\bbF_p
\ar[r]\ar[d]
&0\\
&0&0&0&}
\end{gathered}\hskip 2mm
\end{gather*}

Suppose that $\bbF_p L$ is permutation. Then
since $L$ is coflasque, the sequence
$$0\lra L^{\Gamma} \stackrel{p}{\lra} L^{\Gamma}
\lra (L/pL)^{\Gamma} \lra 0$$ is exact so that
$(\bbF_pL)^{\Gamma} = L^{\Gamma}
/pL^{\Gamma}$. Since $\bbQ [L]\simeq \bbQ
[\Gamma]\oplus \bbQ^2$, $\rank L^{\Gamma} =3$.
But then $\dim^{}_{\bbF_p} (\bbF_pL)^{\Gamma}
= 3$. This means that $\bbF_p L$ must then
have three transitive components. Since $\rank
L=p^2+2$ and $p>2$, this means that $\bbF_p
L=\bbF_p[\Gamma]\oplus \bbF_p^2$.

Looking at the $\bbZ$-basis for $U(p)$ given
above, it is clear that $\bbF_p U(p)\simeq
\bbF_p\oplus \bbF_p I_{\Gamma}$ where $\bbF_p
I_{\Gamma}$ is the augmentation ideal of
$\bbF_p[\Gamma]$. Then the left column of the
last commutative diagram implies that we have
a surjective map $\bbF_p [\Gamma]\oplus
\bbF_p^2\to \bbF_p\oplus \bbF_p I_{\Gamma}$.
Since $(\bbF_p I_{\Gamma})^{\Gamma}=0$, this
would
imply that we have a surjective map %%from
$\bbF_p[\Gamma]\to \bbF_pI_{\Gamma}$ or
equivalently that $\bbF_pI_{\Gamma}$ is a
cyclic $\bbF_p[\Gamma]$-module. But since
$\Gamma$ is a finite $p$-group,
$\bbF_p[\Gamma]$ is a local ring with unique
maximal ideal $\bbF_pI_{\Gamma}$
by~\cite[Corollary 1.4]{car}. Then Nakayama's
Lemma implies that $\bbF_pI_{\Gamma}$ is a
cyclic $\bbF_p[\Gamma]$-module if and only if
$\bbF_pI_{\Gamma}/(\bbF_pI_{\Gamma})^2$ is
generated by one element over $\bbF_p$. Since
%%$\Gamma=C_p\times C_p$,
$\Gamma\simeq \mathbb Z/p\mathbb Z\times
\mathbb Z/p\mathbb Z$, we may use the
K\"{u}nneth formula to show that
$$\bbF_pI_{\Gamma}/\bbF_pI_{\Gamma}^2=
{\rm H}_1(\Gamma,\bbF_p)\cong {\rm
H}_1(\mathbb Z/p\mathbb Z,
%%C_p,
\bbF_p)^2\cong \bbF_p^2.$$
Alternatively, for any $p$-group $H$, we may
show that
$$\bbF_pI_H/\bbF_pI_H^2\longrightarrow H/H^p[H,H],
\ \ \overline{h-1}\mapsto \overline{h}$$ is a
group isomorphism so that, in our case, we
again have
$$\bbF_pI_{\Gamma}/\bbF_pI_{\Gamma}^2\cong \bbF_p^2.$$
Then the above discussion shows that
$\bbF_pI_{\Gamma}$ is not a cyclic
$\bbF_p[\Gamma]$ module so that there is no
surjective map from $\bbF_p[\Gamma]$ to
$\bbF_pI_{\Gamma}$. This implies that
$\bbF_pL$ is not permutation and hence $L$ is
not permutation projective as a $\bbZ
[\Gamma]$-module. This implies in turn that
$\Lambda_{2p}$ is not quasi-permutation as a
$\Gamma$-lattice or as an
%%$\Sym_6$
$\Sym_p\times\Sym_p$-lattice. \quad $\square$
\renewcommand{\qed}{}
\end{proof}

\begin{remark}
Note that this argument fails for $p=2$.
Indeed, we showed that ${\rm rank}\,L = p^2 +
2$ and if $\bbF_pL$ were permutation, it would
have three transitive components. For $p > 2$,
we used these facts to conclude that $\bbF_p L
= \bbF_p[\Gamma] \oplus \bbF_p^2$. For $p=2$,
this is not so; here $\bbF_2 L$ may have three
permutation components, each of rank $2$.
Indeed, if $\Gamma=\langle
g,h\rangle\simeq\bbZ/2\bbZ\times \bbZ/2\bbZ$,
then one can define a surjective $\bbF_2
[\Gamma]$-homomorphism
$$\bbF_2[\Gamma/\langle g\rangle]\oplus
\bbF_2[\Gamma/\langle h\rangle] \oplus
\bbF_2[\Gamma/\langle gh\rangle]\to
\bbF_2I_{\Gamma} \oplus \bbF_2$$ by sending
the generator of the first component to
$(1+g,0)$, the generator of the second
component to $(1+h,0)$, and that of the third
component to $(0,1)$.

In fact, by Proposition~\ref{qndrest}, we see
that $Q_4(2)\vert_{\Gamma} \sim
\Lambda_4\vert_{\Gamma}$. Since $Q_4(2)$ is
the character lattice of the Cayley group
$\SL_4/\m_2\simeq \SO_6$, by Theorem
\ref{thm1} it must be quasi-permutation as an
$\Sym_4$-lattice and hence as a
$\Gamma$-lattice. Alternatively, one could
show directly that $Q_4(2)$ is a
sign-permutation $\Sym_4$-lattice and hence it is
quasi-permutation.
\end{remark}

\subsection { Completion of the
proof of Proposition~\ref{prop.A_n}.}
\label{subsection7.3} It now suffices to prove
the following proposition to complete the
proof of Proposition~\ref{prop.A_n}:

\begin{prop} Suppose $n/d$ is divisible by
a prime $p$.
\begin{enumerate}
\item[{\rm (a)}] If $p > 2$, then the
$\Sym_n$-lattice $Q_n(d)$ is not
quasi-permutation.

\item[{\rm (b)}] If $n>p^2$, then the
$\Sym_n$-lattice $Q_n(d)$ is not
quasi-permutation.
\end{enumerate}
\end{prop}

Indeed, by part
%%(b)
(a), the $\Sym_n$-lattice $Q_n(d)$ is not
%%$\Sym_n$
quasi-permutation if the prime factorization
of
%%$\frac{n}{d}$
$n/d$ includes a prime larger than 2. On the
other hand, if $n/d=2^k$,
%%$\frac{n}{d}=2^k$
then, by part
%%(c)
(b), the $\Sym_n$-lattice $Q_n(d)$ is not
quasi-permutation, for any $(n,d) \ne (4,2)$,
and Proposition~\ref{prop.A_n} follows.

\begin{proof}
(a) Proposition~\ref{qndrest} shows that
$Q_n(d)\vert_{Y_{p}}$ is equivalent to
$\Lambda_{2p}\vert_{\Sym_p\times \Sym_p}$
which is not quasi-permutation by
Proposition~\ref{lambda2p}. Thus $Q_n(d)$ is
not quasi-permutation as a $Y_p$-lattice and
hence as an $\Sym_n$-lattice as well.

(b) We have $n=tp$ with $t> p$. Following the
proof of Proposition 4.1(i) in \cite{ll}, we
define a subgroup $\Gamma \simeq
\bbZ/p\bbZ\times \bbZ/p\bbZ$ of $\Sym_n$ as
follows. Arrange the numbers from $1$ to $n$
into a rectangular table with $p$ columns and
$t$ rows, so that the first row is $1, \dots,
p$, the second row is $p+1, \dots, 2p$, etc.
Let $\sigma_i$ be the $p$-cycle that
cyclically permutes the $i$th row and leaves
elements of all other rows fixed. Note that
$\sigma_1, \dots, \sigma_t$ are commuting
$p$-cycles; explicitly
$$\sigma_i=\bigl((i-1)p+1,
(i-1)p+2,\dots,ip\bigr). $$ We now set $\Gamma
:= \langle\alpha,\beta\rangle$, where
\[ \textstyle
\text{$\alpha:=\prod_{i = 1}^{t-1}\sigma_i
\quad$ and\quad $\beta:=\prod_{i =
1}^{p-1}\sigma_i^{-i}\cdot\prod_{i=p+1}^t\sigma_i$.}
\]
The subgroup $\Gamma$ has orbits
$\O_i=\{(i-1)p+1,(i-1)p+2,\ldots,ip\}$,
$i=1,\dots,t$, all of length $p$ and every
cyclic subgroup $C$ of $\Gamma$ has fixed
points. This means that by Lemma~\ref{Ancoh}
$$
\text{${\rm H}^1(\Gamma,\mathbb
Z\A_{n-1})\simeq\bbZ/p\bbZ$ \hskip 2.5mm but
\hskip 2.5mm ${\rm H}^1(C,\mathbb
Z\A_{n-1})=0$.  }
$$
Also by Lemma~\ref{H1triv}, we find that
$$ \text{${\rm H}^1(C,Q_n\left(n/d \right))=0$.} $$
Then Lemma~\ref{connecthom} and the fact that
$p$ divides $n/d$ show that $
\bbZ/d\bbZ\stackrel{\partial}{\to} {\rm
H}^1(\Gamma ,\mathbb Z\A_{n-1})$ is the zero
map. The following commutative diagram
$$
\xymatrix{ \bbZ/d\bbZ\ar[r]^{0}\ar[d]^{\Res}
&{\rm H}^1(\Gamma,\mathbb Z\A_{n-1})
\ar[r]\ar[d]^{\Res} & {\rm
H}^1\left(\Gamma,Q_n\left(
%%\frac{n}{d}
n/d\right)\right) \ar[d]^{\Res}
 \\
\prod_{a\in \Gamma}\bbZ/d\bbZ\ar[r]^{0 \quad
\quad \quad} &\prod_{a\in \Gamma} {\rm
H}^1(\gen{a},\mathbb Z\A_{n-1})=0\ar[r] &
\prod_{a\in \Gamma} {\rm
H}^1\left(\gen{a},Q_n\left(
%%\frac{n}{d}
n/d\right)\right)=0 }
$$
 shows that
$$\bbZ/p\bbZ\simeq
{\Sha}^1(\Gamma,\mathbb Z\A_{n-1})
\leqslant{\Sha}^1\!\left(\Gamma, Q_n\left(
%%\frac{n}{d}
n/d\right)\right). $$

%%But then
Now if $M$ were a flasque lattice with
$\rho(Q_n(d))$ = the stable equivalence class of
$M$, then $M^*$ is a coflasque lattice
satisfying
$$0\lra M^*\lra P\lra Q_n\left(
n/d\right)\lra 0,$$ so that by
Lemma~\ref{shalem}(a),
 ${\Sha}^2(\Gamma, M^*)\simeq
 {\Sha}^1(\Gamma, Q_n( n/d))
\ne 0$.  Lemma~\ref{shalem}(c) now shows that
$M^*$ cannot be a direct summand of a
quasi-permutation lattice and hence is not stably
permutation. This implies that $M$ cannot be
stably permutation and so $Q_n(d)$ cannot be
quasi-permutation. \quad $\square$
\renewcommand{\qed}{} \end{proof}

\section{ \bf Stably Cayley groups of type \boldmath$\D_n$}
\label{sect.intD_n} % Section 8

\subsection {\bf Root system of type
\boldmath$\D_n$.} \label{subsection8.1} Let
$\ve_1,\ldots, \ve_n$ be the same as in
Subsection \ref{epsilon}. The root system of
type $\D_n$ is the set
$$
\D_n=\{\pm \ve_i\pm \ve_j \mid 1\leqslant i<
j\leqslant n\}.
$$
It has a base $ \alpha_1, \dots, \alpha_n $,
where
%%$\alpha_i=\ve_i-\ve_{i+1}$ for $i=1,\dots, n-1$
$\alpha_1, \ldots, \alpha_{n-1}$ are given by
\eqref{Aalpha} and $\alpha_n=\ve_{n-1}+\ve_n$.
The fundamental dominant weights of $\D_n$
with respect to this base are $\vp_i = \ve_1 +
\dots + \ve_i$ for $ i=1,\dots, n-2$,
$$\textstyle
\text{$\vp_{n-1}=\frac{1}{2}
\sum_{i=1}^{n-1}\ve_i- \frac{1}{2}\ve_n \quad$
and $\quad
\vp_n=\frac{1}{2}\sum_{i=1}^{n-1}\ve_i+
\frac{1}{2}\ve_n$.}
$$
The Weyl group $W(\D_n)$ of $\D_n$ is
$(\bbZ/2\bbZ)^{n-1} \ltimes \Sym_n$, where
$(\bbZ/2 \bbZ)^{n-1}$ consists of all even
numbers of sign changes on
$\{\ve_1,\dots,\ve_n\}$ and $\Sym_n$ acts via
\eqref{action}. The root and weight
$W(\D_n)$-lattices of $\D_n$ are, respectively,
$\mathbb Z\D_n$ and $\Lambda(\D_{n}):=\bbZ
\vp_1\oplus\ldots\oplus \bbZ\vp_n$.

\subsection { Lattices
\boldmath${\rm Y}_{2m}$ and \boldmath${\rm
Z}_{2m}$.} \label{subsection8.2} As we
explained in Section \ref{stablycayley}, we
are interested in the case where $n$ is even,
$n = 2m$, $m>2$. There are precisely the
following three
%%intermediate sublattices,
lattices between $\Lambda(\D_{2m})$ and $\bbZ
\D_{2m}$:
%%\[ \begin{array}{l}
\begin{gather}\label{XYZ}
{\rm X}_{2m}:=\bbZ \D_{2m}+\bbZ \vp_1,\quad
%%\\
{\rm Y}_{2m}:=\bbZ \D_{2m}+\bbZ\vp_{2m-1},
%%,\quad
%%\\
%%\quad
\quad \text{and}\quad
%%\\
{\rm  Z}_{2m}:=\bbZ \D_{2m}+\bbZ\vp_{2m}.
\end{gather}
%%\end{array} \]
%%These are, respectively,
 The
character lattice of $\Spin^{1/2}_{4m}$ (see
Section \ref{stablycayley}) is isomorphic to
either of the lattices ${\rm Y}_{2m}$ and
${\rm Z}_{2m}$ while ${\rm X}_{2m}$ is
isomorphic to the character lattice of
$\SO_{4m}$. Note that $\ve_1, \dots, \ve_n$ is
the sign-permutation basis for $X_{2m}$; this
is consistent with the fact that $\SO_{4m}$ is
Cayley; see~Theorem~\ref{thm1}(a). Also note
that
 $$\textstyle\{\frac{1}{2}(\ve_1+\ve_2+\ve_3-\ve_4),
 \frac{1}{2}(\ve_1+\ve_2-\ve_3+\ve_4),
\frac{1}{2}(\ve_1-\ve_2+\ve_3+\ve_4),
\frac{1}{2}(-\ve_1+\ve_2+\ve_3+\ve_4)\}$$ is
the sign-permutation basis for ${\rm Y}_4$,
and
$$\textstyle\{\frac{1}{2}(\ve_1+\ve_2+\ve_3+\ve_4),
\frac{1}{2}(\ve_1+\ve_2-\ve_3-\ve_4),
\frac{1}{2}(\ve_1-\ve_2+\ve_3-\ve_4),
\frac{1}{2}(-\ve_1+\ve_2+\ve_3-\ve_4)\}$$ is
that for ${\rm Z}_4$; this is consistent with
the fact that $\Spin^{1/2}_{8}$ is Cayley  (see
Section \ref{stablycayley}).
%%Example \ref{spin}.

Our goal is to prove
Proposition~\ref{prop.D_n}. In view of the
aforesaid, this is equivalent to proving the
following.

\begin{prop} The $W(\D_{2m})$-lattices
${\rm Y}_{2m}$ and ${\rm  Z}_{2m}$ are not
%%$W(\D_{2m})$
quasi-permutation for any $m > 2$.
\end{prop}

\begin{proof} For the subgroup $\Sym_{2m}$ of
$W(\D_{2m})$ acting by \eqref{action}, we
consider the $\Sym_{2m}$-lattices ${\rm
Y}_{2m}|_{\Sym_{2m}}$ and ${\rm
Z}_{2m}|_{\Sym_{2m}}$
 %%We will restrict ${\rm Y}_{2m}$ and ${\rm  Z}_{2m}$ to $\Sym_{2m}
%%\subset W(\D_{2m})$
and compare
%%the resulting
%%$\Sym_{2m}$-lattice
them to the $\Sym_{2m}$-lattice $Q_{2m}(m)$
defined by \eqref{e.q_n(d)} and \eqref{ZA},
\begin{equation}\label{Q}
\textstyle Q_{2m}(m)=\mathbb
\bbZ\alpha_1+\ldots+\bbZ\alpha_{2m-1}+
\bbZ\beta, \quad \text{where}\quad
\beta:=m\ve_1-\frac{1}{2}\sum_{i=1}^{2m}\ve_i,
\end{equation}
that is isomorphic to the character lattice of
$\SL_{2m}/ \m_{m}$; see
Subsection~\ref{epsilon}.

First we observe that
$$\textstyle
\alpha_1,\dots, \alpha_{2m - 2}, \gamma,
\ve_{2m - 2}+ \ve_{2m-1}, \quad
\text{where}\quad
\gamma:=\frac{1}{2}\sum_{i=1}^m \ve_i-
\frac{1}{2}\sum_{i=m+1}^{2m}\ve_i,
 $$ is a
%%$\bbZ$-
basis for ${\rm Y}_{2m}$ if $m$ is odd and
for ${\rm Z}_{2m}$ if $m$ is even. Since
$\alpha_1,\ldots,\alpha_{2m-2},
\ve_{2m-2}+\ve_{2m-1}$ is a basis for
$\bbZ\D_{2m-1}$, \eqref{XYZ} implies that
proving this claim is equivalent to proving
the equality
\begin{equation}\label{equality1}
\bbZ \D_{2m-1}+\bbZ\gamma=\begin{cases} \bbZ
\D_{2m}+\bbZ
%%\omega
\vp_{2m-1} &\text{if $m$ is
%%\mbox{
odd},\\
%%$$
%%$$\bbZ D_{2m-1}+\bbZ\gamma=
\bbZ \D_{2m}+\bbZ
%%\omega
\vp_{2m} &\text{if $m$
%%\mbox{
is even}.
\end{cases}
%%$$
\end{equation}
Note that
\begin{gather*}
\label{equality2}
\begin{gathered}\textstyle
\vp_{2m-1}-\gamma=\sum_{i=m+1}^{2m-1} \ve_i\in
\bbZ \D_{2m-1} \quad\text{if $m$
is odd},\\
\textstyle
\vp_{2m}+\gamma=\sum_{i=1}^{m}\ve_i\in \bbZ
\D_{2m-1}\quad\text{if $m$ is even.}
\end{gathered}
\end{gather*}
Therefore proving \eqref{equality1} is
equivalent to
%%the
proving the inclusion
$$
\bbZ\D_{2m}\subseteq \bbZ\D_{2m-1}+\bbZ\gamma,
$$
which in turn is equivalent to proving the
inclusions
\begin{equation*}\label{pmepsilon}
\ve_{2m-1}\pm\ve_{2m}\in\bbZ\D_{2m-1}+\bbZ\gamma.
\end{equation*}
Finally, the last inclusions indeed hold as we
have
\begin{gather*}\textstyle
2\gamma+(\ve_{2m-1}+\ve_{2m})
=\sum_{i=1}^m\ve_i-\sum_{i=m+1}^{2m-2}
\ve_i\in \bbZ \D_{2m-1},\\
\textstyle 2\gamma-(\ve_{2m-1}-
\ve_{2m})=\sum_{i=1}^{m-1}(\ve_i-
\ve_{m+i})+(\ve_m-\ve_{2m-1}) \in \bbZ
\D_{2m-1}.
\end{gather*}
Thus the claim is proved.

Furthermore, the easily checked equalities
\begin{gather*}\textstyle
\beta=\gamma+\sum_{i=1}^{m-1}(m-i)\alpha_i,\\
\textstyle
\alpha_{2m-1}=2\gamma-\sum_{i=1}^{m}i\alpha_i
- \sum_{i=1}^{m-2}(m-i)\alpha_{m+i},
\end{gather*}
and \eqref{Q} imply that
%%we claim that
%%$$\textstyle
$ \alpha_1,\dots,\alpha_{2m-2}, \gamma $ is a
$\bbZ$-basis for $Q_{2m}(m)$.

We thus obtain the following exact sequences
of $\Sym_{2m}$-lattices:
$$ 0\lra Q_{2m}(m)\lra {{\rm Y}_{2m}}\vert_{\Sym_{2m}}
\lra \bbZ \lra 0
$$ if $m$ is odd and
$$0\lra Q_{2m}(m)\lra {{\rm  Z}_{2m}}\vert_{\Sym_{2m}}
\lra \bbZ\lra 0 $$ if $m$ is even. Here the
$\Sym_{2m}$-lattice $\bbZ$ is generated by
$\ve_{2m-2} + \ve_{2m-1}$ modulo $Q_{2m}(m)$.
We claim that the $\Sym_{2m}$-action on this
lattice is trivial. Indeed, on the one hand,
the al\-ter\-nating subgroup of $\Sym_{2m}$
has to act on this lattice trivially because
it has no non-trivial one-dimensional
representations. On the other hand, as $m>2$,
the transposition $(1,2)$
%%clearly
acts trivially on $\ve_{2m-2} + \ve_{2m-1}$.
Since the alternating subgroup and the
transposition $(1,2)$ generate $\Sym_{2m}$,
this proves the claim.

The above exact sequences thus tell us that
${{\rm Y}_{2m}} \vert_{\Sym_{2m}}\sim
Q_{2m}(m)$ if $m$ is odd and ${{\rm Z}_{2m}}
\vert_{\Sym_{2m}}\sim Q_{2m}(m)$ if $m$ is
even. By Proposition~\ref{prop.A_n}, the
$\Sym_{2m}$-lattice $Q_{2m}(m)$ is not
quasi-permutation for any $m > 2$. Thus for $m
> 2$, the $W(\D_{2m})$-lattice ${{\rm Y}_{2m}}$ is
not quasi-permutation if $m$ is odd, and the
$W(\D_{2m})$-lattice ${\rm  Z}_{2m}$ is not
quasi-permutation if $m$ is even, as their
restrictions to $\Sym_{2m}$ are not
quasi-permutation. Since $Y_{2m} \simeq
Z_{2m}$ as $W(\D_{2m})$-lattices, this
completes the proof. \quad $\square$
\renewcommand{\qed}{} \end{proof}

\section{\bf Which stably Cayley groups are Cayley?}
\label{sect.sl_3} % Section 9

In this section we will prove
Theorem~\ref{thm3}. The groups $G = \SO_n$,
$\Symp_{2n}$, and $\PGL_n$ are shown to be
Cayley in Examples~\ref{so}
and~\ref{ex.pgl_n}. It thus remains to
consider $\SL_3$ and ${\bf G}_2$.

\subsection { The group $\SL_3$.}
\label{subssl3} % Subsection 9.1

\begin{prop} \label{prop.sl_3}
The group $\SL_3$ is Cayley.
\end{prop}

The proof below is based on analysis of the
explicit formulas in \cite[4.9]{voskresenskii}
and the geometric ideas of the proof of
Proposition \ref{prop.sl_3} given in
\cite{popov-luna}. We present it in a form
that will also help us prove that ${\bf
G}_2\times {\bf G}_m^2$ is Cayley; see
Proposition~\ref{prop.g_2} below. On the other
hand, the spirit of the arguments in
\cite{popov-luna} is close to that in
\cite{iskovskih1}. Since \cite{iskovskih1} is
the main ingredient we will use in showing
that ${\bf G}_2$ is not Cayley, see
Lemma~\ref{lem.iskovskih} and
Proposition~\ref{prop.iskovskih} below, we
will give an outline of the proof of
Proposition~\ref{prop.sl_3} from
\cite{popov-luna} in the Appendix.

\begin{proof} The Weyl group $W$
of $\SL_3$ is $\Sym_3$. Consider the following
subalgebra $D$ of ${\rm Mat}_{3\times 3}$:
\begin{equation}\label{D3}
D:=\{\diag(a_1,a_2,a_3)\in {\rm Mat}_{3\times
3}\mid a_i\in k\}
\end{equation}
and the action of $\Sym_3$ on $D$ given by
\begin{equation}\label{sigma-action}
\sigma(\diag(a_1, a_2,
a_3)):=\diag(a_{\sigma(1)}, a_{\sigma(2)},
a_{\sigma(3)})\quad\text{where}\hskip 2mm
\sigma\in\Sym_3.
\end{equation}
The $\Sym_3$-stable subvarieties
\begin{gather}\label{tori}
T =\{X\in D\mid \det X=1\}
 \quad \text{and}\quad
\lt =\{Y\in D\mid \tr Y=0\}
\end{gather}
are, respectively, the maximal torus of $\SL_3$
and its Lie algebra, considered as
$W$-varieties. By the Corollary
%%~\ref{cor.torus}
of Lemma~\ref{cor.caley}, it suffices to show
that $T$ and $\lt$ are birationally
$\Sym_3$-isomorphic.

Let $D\setminus\{0\}\to \bbP(D)$,
$X\mapsto[X]$, be the natural projection.
Denote by $\bbP^2_{\Sym_3\text{-natural}}$ and
$\bbP^2_{\Sym_3\text{-twisted}}$ the
projective plane $\bbP(D)$ endowed, 
respectively, with the na\-tu\-ral and
``twisted'' rational actions of $\Sym_3$ given
by
$$
\sigma([X]):=[\sigma(X)]\hskip
4mm\text{and}\hskip
3mm\sigma([X]):=[\sigma(X)^{\sign\sigma}],
\hskip 3mm\text{where} \hskip 2mm
\sigma\in\Sym_3, X\in D.
$$

Let $\pi: \SL_3\to\PGL_3$ be the natural
projection. Since $d_e\pi$ is an isomorphism
between the Lie algebras of $\SL_3$ and
$\PGL_3$ and since $\PGL_3$ is a Cayley group, see
Example \ref{ex.pgl_n}, the Corollary of
Lemma~\ref{cor.caley} tells us that $\lt$ is
birationally $\Sym_3$-isomorphic to the
maximal torus $\pi(T)$ of $\PGL_3$. In turn,
we have the following birational
$\Sym_3$-isomorphisms of $\Sym_3$-varieties:
\begin{gather*}
\pi(T)\stackrel{\simeq}{\dasharrow}
\bbP_{\Sym_3\text{-natural}}^2, \quad
\pi(X)\mapsto
[X],\\
\textstyle \bbP_{\Sym_3\text{-twisted}}^2
\stackrel{\simeq}{\dasharrow} T, \quad
[\diag(a_1, a_2, a_3)]\mapsto {\rm
diag}\bigl({a_2}/{a_3}, {a_3}/{a_1},
{a_1}/{a_2} \bigr).
\end{gather*}

Thus we only need to show that
$\bbP_{\Sym_3\text{-natural}}^2$ and
$\bbP_{\Sym_3\text{-twisted}}^2$ are
birationally $\Sym_3$-isomorphic. We shall
establish this in three steps.

{\it Step $1$.} Consider the action of
$\Sym_3$ on $\lt\times\lt$ given by
\begin{equation}\label{S3twisted}
\sigma(Y, Z):=
\begin{cases}
 \bigl(\sigma(Y), \sigma(Z)\bigr)
 &\text{if $\sigma$ is even,} \\
 \bigl(\sigma(Z), \sigma(Y)\bigr)
 &\text{if $\sigma$ is odd,}
\end{cases}
\hskip 4mm\text{where}\hskip 2mm
\sigma\in\Sym_3,\ Y, Z\in\lt.
\end{equation}
It determines the action of $\Sym_3$ on the
surface $\bbP(\lt)\times \bbP(\lt)$. Denote
resulting $\Sym_3$-surface by
$(\bbP(\lt)\times
\bbP(\lt))_{\Sym_3\text{-twisted}}$.

We claim that the $\Sym_3$-varieties
$\bbP^2_{\Sym_3\text{-twisted}}$ and
$(\bbP(\lt)\times \bbP
(\lt))_{\Sym_3\text{-twisted}}$ are
birationally $\Sym_3$-isomorphic. Indeed, it
is immediately seen that the rational map
\begin{equation*}\textstyle
\varphi:
\bbP^2_{\Sym_3\text{-twisted}}\dasharrow
(\bbP(\lt)\times \bbP
(\lt))_{\Sym_3\text{-twisted}},\quad
[X]\mapsto \bigl(\bigl[X-\frac{\tr
(X)}{3}I_3\bigr], \bigl[X^{-1}-\frac{\tr
(X^{-1})}{3}I_3\bigr]\bigr),
\end{equation*}
is $\Sym_3$-equivariant and we shall now
construct a rational map inverse to $\varphi$.
Note that for $Y, Z\in \lt$ in general
position,  $Y$, $Z$, $I_3$ form a basis of the
vector space $D$.
%%$k^3$.
Thus there are unique $\alpha, \beta, \gamma
\in k$ such that
\[ \alpha Z + \beta Y + \gamma I = - YZ. \]
Note that $\alpha$, $\beta$, and $\gamma$ are,
in fact, bihomogeneous rational functions of
$Y$ and $Z$ of bidegree $(1, 0)$, $(0, 1)$, and
$(1, 1)$, respectively. We now
%%define
consider the map
%%$h$
\begin{equation*}\label{psi}
\psi \colon (\bbP(\lt)\times \bbP
(\lt))_{\Sym_3\text{-twisted}}\dasharrow
\bbP^2_{\Sym_3\text{-twisted}}, \quad ([Y],
[Z]) \mapsto [Y + \alpha I_3].
\end{equation*}
%%\[ ([Y], [Z]) \mapsto [Y + \alpha I_3]. \]
To compute $\psi \circ \varphi$, note that
 if $Y = X - \frac{\tr(X)}{3} I_3$ and
$Z = X^{-1} - \frac{\tr(X^{-1})}{3} I_3$, then
expanding
\begin{equation*} \textstyle
I_3 = \bigl(Y + \frac{\tr(X)}{3} I_3\bigr)
\bigl(Z + \frac{\tr(X^{-1})}{3} I_3\bigr),
\end{equation*}
we see that $\alpha = \frac{\tr(X)}{3}$ and
thus $\psi([Y], [Z]) = [X]$. Thus
$\psi\circ\varphi={\rm id}$, and hence
$\varphi$ is a birational
$\Sym_3$-isomorphism.

\smallskip

{\it Step $2$.} We now consider the linear
action of $\Sym_3$ on $\lt \otimes \lt$
%%given
determined by the action \eqref{S3twisted} and
the corresponding action of $\Sym_3$ on
$\bbP(\lt \otimes \lt)$.
 Then the
Segre embedding
%%use the Segre
%%embedding
\begin{equation*}
%%\iota:
(\bbP(\lt) \times
\bbP(\lt))_{\Sym_3\text{-twisted}}
\hookrightarrow \bbP(\lt \otimes \lt)
\end{equation*} is $\Sym_3$-equivariant.
Its image is a quadric $Q$ in $\bbP(\lt
\otimes \lt)$ described as follows. Choose a
basis $D_1 := {\rm diag}(1, \zeta, \zeta^2)$,
$D_2 := {\rm diag}(1, \zeta^2, \zeta)$ of
$\lt$, where $\zeta$ is a primitive cube root
of unity. Set $D_{ij} = D_i \otimes D_j$.
%% for $i, j = 1, 2$.
Then
%%quadric $Q$ is given by
\begin{equation}\label{quadric} Q = \{ (\alpha_{11}: \alpha_{12}:
\alpha_{21}: \alpha_{22}) \mid \alpha_{11}
\alpha_{22} = \alpha_{12} \alpha_{21} \},
\end{equation} where $ (\alpha_{11}:
\alpha_{12}: \alpha_{21}: \alpha_{22})$ is the
point of $\bbP(\lt \otimes \lt)$ corresponding
to $\alpha_{11} D_{11} + \alpha_{12} D_{12} +
\alpha_{21} D_{21} + \alpha_{22} D_{22}\in\lt
\otimes \lt$.

\smallskip

{\it Step $3$.}  Decomposing $\lt \otimes \lt$
as a sum of
%%irreducible
$\Sym_3$-submodules,
%%representations,
we obtain
\begin{equation}\label{S3decomposition} \lt
\otimes \lt = V_1 \oplus V_2 \oplus V_3,
\end{equation} where $V_1 =kD_{11}+ kD_{22}
%%\Span_k (e_{11}, e_{22})
$ is a simple
%%a
%%faithful
2-dimensional submodule and $V_2 =kD_{12}$,
$V_3 =kD_{21}$ are trivial 1-dimensional
submodules. Since the $\Sym_3$-fixed point
$(0:0:1:0)\in \bbP(\lt \otimes \lt)$
corresponding to $V_3$ lies on $Q$, the
stereographic projection $Q \dasharrow
\bbP(V_1 \oplus V_2)$ from this point is a
birational $\Sym_3$-isomorphism.

Finally, the $\Sym_3$-module $D$ is isomorphic
to $V_1\oplus V_2$. Hence $\bbP(V_1\oplus
V_2)$ and $ \bbP^2_{\Sym_3\text{-natural}}$
are $\Sym_3$-isomorphic.

To sum up, we have established the existence
of the following birational
$\Sym_3$-isomorphisms:
\[
\xymatrix{ \bbP^2_{\Sym_3\text{-twisted}}
\ar@{-->}[r]^{\hskip -10mm\stackrel{\rm Step
1}{\simeq}}& (\bbP(\lt) \times
\bbP(\lt))_{\Sym_3\text{-twisted}}\ar@{-->}[r]^{\hskip
13mm\stackrel{\rm Step 2}{\simeq}}& Q
\ar@{-->}[r]^{\hskip -5mm\stackrel{\rm Step
3}{\simeq}}& \bbP^2_{\Sym_3\text{-natural}}. }
\]
This completes the proof of
Proposition~\ref{prop.sl_3}. \quad $\square$
\renewcommand{\qed}{} \end{proof}

\subsection{The group ${\bf
G}_2$.} \label{sect.g_2} % Subsection 9.2

The Weyl group of ${\bf G}_2$ is the dihedral
group $\Sym_3 \times \Sym_2$ of order $12$.
The maximal torus of ${\bf G}_2$ and its Lie
algebra are $\Sym_3 \times \Sym_2$-isomorphic,
respectively, to $T$ and $\lt$ given by
\eqref{tori}, where the action of the first
factor of $\Sym_3\times \Sym_2$ is defined, as
in the case of $\SL_3$, by
\eqref{sigma-action}, and that of the
non-trivial element $\theta$ of the second
factor by
\begin{equation}\label{theta}
\theta(X):=X^{-1}\ \text{for} \ X\in
T\quad\text{and}\quad \theta(Y):=-Y\
\text{for}\ Y\in\lt.
\end{equation}

We begin with the following surprising recent
result due to {\sc
Iskov\-skikh},~\cite{iskovskih1}.

\begin{lem} \label{lem.iskovskih}
The $\Sym_3\times\Sym_2$-varieties $T$ and
$\lt$ are not birationally
$\Sym_3\times\Sym_2$-isomorphic.
\end{lem}

\noindent{\em Proof outline.} Since $T$ and
$\lt$ are rational surfaces, the theory of
rational $G$-surfaces, due to {\sc
Manin}~\cite{manin} and {\sc
Iskovskikh}~\cite{iskovskih2},
\cite{iskovskih3}, can be applied; this is
precisely what is done in~\cite{iskovskih1}
(see also~\cite{iskovskih5}). Minimal rational
$\Sym_3\times\Sym_2$-surfaces are known, and
any equivariant birational isomorphism between
two such surfaces can be written as a
composition of so-called ``elementary links",
which are completely enumerated
in~\cite{iskovskih3}. The argument
in~\cite{iskovskih1} and~\cite{iskovskih5}
amounts to constructing suitable minimal
models for $T$ and $\lt$ and explicitly
checking that it is impossible to get from one
to the other by a sequence of elementary
links. \quad $\square$

\begin{prop} \label{prop.iskovskih} ${\bf G}_2$ is not a Cayley group.
\end{prop}

\begin{proof} By the Corollary of Lemma~\ref{cor.caley},
this follows from Lemma~\ref{lem.iskovskih}.
\quad $\square$
\renewcommand{\qed}{} \end{proof}

The following result illustrates how delicate
the matter is.

\begin{prop} \label{prop.g_2} ${\bf G}_2
\times {\bf G}_m^2$ is a Cayley group.
\end{prop}
\begin{proof}
By the Corollary of Lemma~\ref{cor.caley}, it
suffices to show that $T \times \bbA^2$ and
$\lt \times \bbA^2$ are birationally
$\Sym_3\times \Sym_2$-isomorphic, where in
both cases $\Sym_3\times \Sym_2$ acts via the
first factor. We shall define a birational
$\Sym_3\times \Sym_2$-isomorphism between them
in three steps.

\smallskip

{\it Step $1$.} Let $(\lt\times
\lt)_{\Sym_3\times\Sym_2\text{-twisted}}$ be
the variety $\lt\times \lt$ endowed with the
following $\Sym_3\times\Sym_2$-action:
\begin{equation} \label{e.W-action1}
(\sigma, \varepsilon) (Y, Z)\!:=\!
\begin{cases}
\sign(\sigma) (\sigma(Y), \sigma(Z))&
 \text{if \ $\sign(\sigma) =\sign(\varepsilon)$, } \\
\sign(\sigma) (\sigma(Z), \sigma(Y)) &
\text{otherwise,}
\end{cases}
\end{equation}
for any $(\sigma,
\varepsilon)\in\Sym_3\times\Sym_2$ and $Y, Z
\in \lt$.  This action descends to
$\bbP(\lt)\times \bbP(\lt)$; denote the
resulting $\Sym_3\times\Sym_2$-variety by
$(\bbP(\lt) \times
\bbP(\lt))_{\Sym_3\times\Sym_2
\text{-twisted}}$. We claim that
$(\lt\times\lt)_{\Sym_3\times\Sym_2
\text{-twisted}}$ is birationally isomorphic
to $(\bbP(\lt) \times
\bbP(\lt))_{\Sym_3\times\Sym_2
\text{-twisted}}\times \bbA^2$ as an
$\Sym_3\times\Sym_2$-variety. Here
$\Sym_3\times\Sym_2$ acts trivially on
$\bbA^2$.

To prove the claim, let $\lt '$ be $\lt$ blown
up at the origin. The
$\Sym_3\times\Sym_2$-action
\eqref{e.W-action1} on $\lt \times \lt$ lifts
to $\lt ' \times \lt'$; we shall denote the
resulting $\Sym_3\times\Sym_2$-variety by
$(\lt ' \times\lt ')_{\Sym_3\times\Sym_2
\text{-twisted}}$. The natural projection $\lt
\dasharrow \bbP(\lt)$ (which is only a
rational map, not defined at the origin)
lifts to a regular map $\lt ' \to \bbP(\lt)$.
Moreover, the natural projection
\begin{equation*}
(\lt ' \times\lt ')_{\Sym_3\times\Sym_2
\text{-twisted}}\lra( \bbP(\lt)\times
\bbP(\lt))_{\Sym_3\times\Sym_2
\text{-twisted}}
\end{equation*}
is an algebraic vector
$\Sym_3\times\Sym_2$-bundle of rank $2$. Since
$\Sym_3\times\Sym_2$ acts on $(\bbP(\lt)\times
\bbP(\lt))_{\Sym_3\times\Sym_2
\text{-twisted}}$ faithfully,
Lemma~\ref{lem.no-name}(b) shows that $(\lt '
\times \lt ')_{\Sym_3\times\Sym_2
\text{-twisted}}$ is birationally isomorphic,
as an $\Sym_3 \times \Sym_2$-variety, to
$(\bbP(\lt)\times
\bbP(\lt))_{\Sym_3\times\Sym_2
\text{-twisted}}\times \bbA^2$ (where
$\Sym_3\times\Sym_2$ acts via the first
factor, as above). Since $(\lt  \times \lt
)_{\Sym_3\times\Sym_2 \text{-twisted}}$ and
$(\lt ' \times \lt ')_{\Sym_3\times\Sym_2
\text{-twisted}}$ are birationally
$\Sym_3\times\Sym_2$-isomorphic, this proves
the claim.

\smallskip

{\it Step $2$.} Let
$\bbP^2_{\Sym_3\times\Sym_2 \text{-twisted}}$
be the projective plane $\bbP(D)$ endowed with
the action of $\Sym_3\times\Sym_2$ given by
$$ (\sigma, \varepsilon)([X]):=
[\sigma(X)^{\sign \sigma\,\sign\varepsilon}],
\hskip 4mm\text{where}\hskip 2mm
(\sigma,\varepsilon)\in\Sym_3\times\Sym_2,\
X\in D.  $$ Then the rational maps
\begin{equation*}
\textstyle \bbP^2_{\Sym_3\times\Sym_2
\text{-twisted}} \dasharrow T,\quad [{\rm
diag}(a_1, a_2, a_3)]\mapsto {\rm diag} \bigl(
{a_2}/{a_3}, {a_3}/{a_1}, {a_1}/{a_2} \bigr),
\end{equation*}
and
\begin{eqnarray*}\textstyle \bbP^2_{\Sym_3\times\Sym_2
\text{-twisted}}&\dasharrow& ( \bbP(\lt)\times
\bbP(\lt))_{\Sym_3\times\Sym_2
\text{-twisted}},\\  &[X]&\mapsto
\bigl(\bigl[X-\frac{\tr (X)}{3}I_3\bigr],
\bigl[X^{-1}-\frac{\tr
(X^{-1})}{3}I_3\bigr]\bigr),
\end{eqnarray*}
are birational
$\Sym_3\times\Sym_2$-isomorphisms\,---\,the
arguments are similar to those in the proof of
Proposition~\ref{prop.sl_3}.

\smallskip

{\it Step $3$.} The definition of $(\lt \times
\lt)_{\Sym_3\times\Sym_2 \text{-twisted}}$ in
Step 1 shows that the map
\begin{equation*}
(\lt \times \lt)_{\Sym_3\times\Sym_2
\text{-twisted}} \lra \lt, \quad (t_1, t_2)
\mapsto t_1 - t_2,
\end{equation*}
is $\Sym_3\times\Sym_2$-equivariant. Hence,
this map may be viewed as an algebraic
$\Sym_3\times\Sym_2$-vector
%%$\Sym_3\times\Sym_2$-
bundle of rank $2$. Since
$\Sym_3\times\Sym_2$ acts on $\lt$ faithfully,
applying Lemma~\ref{lem.no-name}(b) once
again, we conclude that $(\lt \times
\lt)_{\Sym_3\times\Sym_2 \text{-twisted}}$ is
birationally $\Sym_3\times\Sym_2$-isomorphic
to $\lt \times \bbA^2$, where
$\Sym_3\times\Sym_2$ acts via the first
factor.

To sum up, we have established the existence
of the following birational
$\Sym_3\times\Sym_2$-iso\-mor\-phisms:
\begin{equation*}
\xymatrix {{T \times \bbA^2}
\ar@{-->}[r]^{\hskip -19mm\stackrel{\rm Step
2}{\simeq}} & (\bbP(\lt) \times
\bbP(\lt))_{\Sym_3\times\Sym_2
\text{-twisted}} \times
\bbA^2\ar@{-->}[r]^{\hskip 10mm\stackrel{\rm
Step 1}{\simeq}}& (\lt \times
\lt)_{\Sym_3\times\Sym_2 \text{-twisted}}
\ar@{-->}[r]^{\hskip 10mm\stackrel{\rm Step
3}{\simeq}}& \lt \times \bbA^2.}
\end{equation*}
This completes the proof of
Proposition~\ref{prop.g_2}. \quad $\square$
\renewcommand{\qed}{}
\end{proof}

\begin{remark}
We do not know whether or not ${\bf G}_2
\times {\bf G}_m$ is a Cayley group.
\end{remark}

%%%%%%%%%%%%%%%%%%%%%%%%%%%%%%%%%%%%%%%%%

\section{\bf Generalization}
\label{generalization} % Section 10

The notions of Cayley map and Cayley group
naturally lead to generalizations which will
be considered in this section.

\subsection { Generalized Cayley maps.}
\label{subsection10.1} Let $G$ be a connected
linear al\-geb\-raic group and let $\g$ be its
Lie algebra. We consider $G$ and $\g$ as
$G$-varieties with respect to the conjugating
and adjoint actions, respectively, and denote by
${\rm Rat}_G(G, \g)$ the set of all rational
$G$-maps $G\dasharrow \g$ endowed with the
natural structure of a vector space over
$k(G)^G$. Set ${\rm Mor}_G(G,
\g):=\{\varphi\in {\rm Rat}_G(G,
\g)\mid\text{$\varphi$ is a morphism}\}$.

\begin{defn} \label{generalized}
An element $\varphi\in {\rm Rat}_G(G, \g)$
(respectively, $\varphi\in {\rm Mor}_G(G,
\g)$) is called a {\it generalized Cayley map}
(respectively, {\it generalized Cayley
morphism}) {\it of} $G$ if $\varphi$ is a
dominant map.
\end{defn}

We are now ready to state the main result of
this subsection.

\begin{thm} \label{thm.gen-morph}
Every connected linear algebraic group admits
a generalized Cayley morphism.
\end{thm}

Our proof of Theorem~\ref{thm.gen-morph} will
proceed in three steps. First we will
construct a ge\-neralized Cayley morphism for
every reductive group (Corollary to
Lemma~\ref{fixed}), and then a generalized Cayley
map for an arbitrary linear algebraic group
(Proposition~\ref{prop.cayley-map}), then a
generalized Cayley morphism for
%%an map for
an arbitrary linear algebraic group.

Our construction in the case of reductive
groups relies on the following known fact; see
\cite[Lemme III.1]{lu3} and 
cf.\,\cite[6.3]{popov-vinberg}.

\begin{lem}\label{fixed}
Assume that the group $G$ is reductive. Let
$X$ be an affine algebraic variety endowed
with an algebraic action of $G$ and let $x\in
X$ be a non-singular fixed point of $G$. Let
${\rm T}_x$ be the tangent space of $X$ at $x$
endowed with the natural action of $G$. Then
there is a $G$-morphism $\varepsilon: X\to
{\rm T}_{x}$ \'etale at $x$ {\rm(}\hskip -.5mm
hence dominant{\rm)} and such that
$\varepsilon(x)=0$.
\end{lem}

\begin{proof}
We may assume without loss of generality that
$X$ is a $G$-stable subvariety of a finite-dimensional 
algebraic $G$-module $V$; see
\cite[Theorem 1.5]{popov-vinberg}. Since $x$
is a fixed point of $G$, we can replace $X$ by
its image under the parallel translation
$v\mapsto v-x$ and assume that $x=0$. The
tangent space ${\rm T}_x$ is identified with a
submodule of $V$. Since $G$ is reductive, the
$G$-module $V$ is semisimple. Hence $V={\rm
T}_x\oplus M$ for some submodule $M$. Now we
can take $\varepsilon=\pi|_{X}$, where
$\pi:V\to{\rm T}_x$ is the projection parallel
to $M$.\quad $\square$
\renewcommand{\qed}{}
\end{proof}

Taking $X=G$ with the conjugating action and
$x=e$, we obtain the following.

\medskip

\noindent {\bf Corollary.} {\it Assume that
$G$ is reductive. Then there is a generalized
Cayley morphism $\varphi$ of $G$ \'etale at
$e$ and such that $\varphi (e)=0$.}

\medskip

The following special case of this
construction was considered by {\sc Kostant}
and {\sc Michor}, \cite{komi}.

\begin{example}\label{km} Assume that $G$ is reductive.
Consider an algebraic homomorphism $\nu: G\to
{\GL}(S)$, where $S$ is a finite-dimensional
vector space over~$k$. Then the $k$-vector
space $V:={\rm End}(S)$ has a natural
$G$-module structure defined by $g(h):=
\nu(g)h\nu (g)^{-1}$ for every $g \in G$ and
$h \in V$. If $\nu$ is injective, identify $G$
with the image of $\iota\circ\nu$, where
$\iota:{\GL}(S)\hookrightarrow V$ is the
natural embedding. Then $G$ is a $G$-stable
subvariety of $V$ and the restriction to
$\g={\rm T}_e$ of the $G$-invariant inner
product $(x, y)\mapsto\tr xy$ on $V$ is
non-degenerate. This yields the $G$-module
decomposition $V=\g\oplus\g^{\perp}$, where
$\g^{\perp}$ is the orthogonal complement to
$\g$ with respect to $(\ {,}\ )$. The
restriction to $G$ of the projection $V\to\g$
parallel to $\g^{\perp}$ is a generalized
Cayley morphism $\varphi: G\to\g$ \'etale at
$e$ such that $\varphi (e)=0$. \quad $\square$
\end{example}

\begin{prop} \label{prop.cayley-map} Every
connected linear algebraic group $G$ admits a
generalized Cayley map.
\end{prop}

\begin{proof}
 We use the notation
of Proposition \ref{prop.levi} and its proof.
The group $W_{L, T}$ is finite, hence
reductive, and $e\in T$ is its fixed point.
Therefore Lemma \ref{fixed} implies that there
is a dominant $W_{L, T}$-morphism
$\varepsilon\!: T\to\lt$. The arguments in the
proof of part (a) of
Proposi\-tion~\ref{prop.levi} show that
$\varepsilon$ is $N$-equivariant. Consider an
$N$-isomorphism \eqref{tau}. Then
$$
\varepsilon\times\tau: C=T\times
U\longrightarrow \lt\oplus \lu=\lc
$$
is a dominant $N$-morphism. Hence by
Lemma~\ref{properties}, there is a dominant
$G$-morphism
$$\theta:
G\times^{N}\!C\longrightarrow G\times^N\!\lc$$
such that $\theta|_C= \varepsilon\times\tau$.
Now, since, by Lemma \ref{lem.cartan}, the
$G$-morphisms $\gamma_C$ and $\gamma_{\lc}$
given by \eqref{cartan-morphisms} are
birational $G$-isomorphisms,
$\gamma_{\lc}\circ\theta \circ\gamma^{-1}_C\in
{\rm Rat}_G(G,\g)$ is a generalized Cayley
map. \quad $\square$
\renewcommand{\qed}{}
\end{proof}

% Which groups admit generalized Cayley
% morphisms? By the Corollary of Lemma
% \ref{fixed}, reductive groups have this
% property.  The following shows that this
% property is shared by a much wider class
% of groups (in particular, by all groups
% whose radical is nilpotent).

Our next task is to deduce
Theorem~\ref{thm.gen-morph} from
Proposition~\ref{prop.cayley-map}. Our
argument will rely on the following simple
lemma.

\begin{lem} \label{lem.semiinv}
%%Let G be a connected algebraic
%%group. Then every
Every semi-invariant for the conjugating
action of G on itself is, in fact, an
invariant.
\end{lem}

\begin{proof} Suppose $t \in k[G]$ is a
semi-invariant.  That is, there exists an
algebraic character $\chi \colon G \to {\bf
G}_m$ such that $t(g h g^{-1}) = \chi(g) t(h)$
for every $g, h \in G$. We may assume $t$ is
not identically zero.
%%We may assume
%%$t$ is not identically zero.
%%Our goal is to show that
%%$\chi$ is the trivial character of $G$.
Setting $g = h$ in the above formula, we
obtain
\[ \text{$t(g) =
%%t(g g g^{-1}) =
\chi(g) t(g)$ for every $g \in G$.
%%in $k$,
} \]
%%i.e., for every $g \in G$, $\chi(g) = 1$ or $t(g) = 0$.
Since $G$ is connected and $t$ is not
identically zero, this implies that
%%either (i) $t(g) = 0$ for
%%every $g \in G$ or (ii)
$\chi(g) = 1$ for every $g \in G$, i.e., $t\in
k[G]^G$. \quad $\square$
\renewcommand{\qed}{}
\end{proof}

Theorem~\ref{thm.gen-morph} is now an
immediate consequence of
Proposition~\ref{prop.cayley-map} and
Proposition~\ref{prop.Rati} below.

\begin{prop}\label{prop.Rati}
Let $\varphi\!\in\! {\rm Rat}_G(G, \g)$. Then
there is $f\!\in\! k[G]^G$ such~that
\begin{enumerate}
\item[\rm (i)] $\{g\in G\mid f(g)=0\}$ is the
indeterminacy locus of $\varphi$, \item[\rm
(ii)] $f\varphi\in {\rm Mor}_G(G, \g)$.
\end{enumerate}
Moreover, if $\varphi$ is a generalized Cayley
map of $G$, then {\rm (ii)} may be replaced by
\begin{enumerate}
\item[$\rm (ii)'$] $f\varphi$ is a generalized
Cayley morphism $G\to\g$.
\end{enumerate}
\end{prop}

\begin{proof}
We may assume that $\varphi$ is not a
morphism. Then the indeterminacy locus of
$\varphi$ is an unmixed closed subset $X$ of
$G$ of codimension~$1$. Since, by
\cite[Theorem 6]{popov-pic}, the Picard group
of the underlying variety of $G$ is finite,
this implies that there is $t\in k[G]$ such
that $\{g\in G\mid t(g)=~0\}=X$. As $\varphi$
is $G$-equivariant, $X$ is $G$-stable. Hence,
by \cite[Theorem~3.1]{popov-vinberg}, $t$ is a
semi-invariant of $G$ and therefore $t\in
k[G]^G$ by Lemma~\ref{lem.semiinv}.
Consequently the function $f=t^m$ satisfies
(i) and (ii) for a sufficiently large positive
integer $m$.  The second assertion of the
proposition follows from Lemma \ref{dominant}
below. \quad $\square$
\renewcommand{\qed}{}
\end{proof}

\begin{lem}\label{dominant}
Let $\psi: X\dasharrow V$ be a dominant
rational map, where $X$ is an irreducible
al\-geb\-raic variety, $V$ a vector space over
$k$, and $\dim X=\dim V$. Then for every
non-zero func\-tion $t\in k(X)$, at least one
of the maps $\alpha:=t\psi$ and
$\beta:=t^2\psi$ is dominant.
\end{lem}

\begin{proof}
Put $h_i:=\psi^*(x_i)\in k(X)$, where
$x_1,\ldots, x_n$ are the coordinate functions
on $V$ with respect to some basis. Then
$K:=\psi^* \bigl(k(V)\bigr)=k(h_1,\ldots,
h_n)$, $K_1:=\alpha^*\bigl(k
(\overline{\alpha(X)})\bigr)=k(th_1,\ldots,
th_n)$ and $K_2:=
\beta^*\bigl(k(\overline{\beta(X)})\bigr)=
k(t^2h_1,\ldots, t^2h_n)$, where the  bar denotes
the closure in $V$. All three fields contain
the subfield $K_0:= k(\ldots,h_i/h_j,\ldots)$.
We have ${\rm trdeg}_kK=n$. Therefore ${\rm
trdeg}_kK_0=n-1$.

Assume the contrary: neither $t\psi$ nor
$t^2\psi$ is dominant. Then ${\rm
trdeg}_kK_1={\rm trdeg}_kK_2=n-1$. Since
$K_1=K_0(th_i)$ and $K_2=K_0(t^2h_i)$ for any
$i$, this implies that both $th_i$ and
$t^2h_i$ are algebraic over $K_0$. Hence
$h_i=(th_i)^2/t^2h_i$ is algebraic over $K_0$.
Thus $K$ is algebraic over $K_0$. Hence ${\rm
trdeg}_kK={\rm trdeg}_kK_0=n-1$, a
contradiction. \quad $\square$
\renewcommand{\qed}{}
\end{proof}

\subsection { The Cayley degree.}
\label{subsection10.2} Note that every
generalized Cayley map $\varphi: G\dasharrow
\g$ has finite degree, i.e., $\deg
\varphi:=[k(G):\varphi^*(k(\g))]<\infty$. By
Definition \ref{Cayley}, Cayley maps are
exactly  generalized Cayley maps of degree
$1$. This naturally leads to the following
definition of a ``measure of non-Cayleyness''
of $G$.

\begin{defn}\label{cayley-deg}
The {\it Cayley degree} of $G$ is the number
${\rm Cay}(G):= \underset{\varphi}{\min}\deg
\varphi$, where $\varphi$ runs through all
generalized Cayley maps of $G$.
\end{defn}

% If $G$ is defined over a subfield $K$ of $k$,
% then the {\it Cayley $K$-degree} ${\rm Cay}_K(G)$
% of $G$ is defined in a similar manner by taking
% the minimum only over those generalized Cayley maps
% $\varphi$ that are defined over $K$.
%% range over all generalized Cayley maps of $G$ defined over $K$.

Clearly $G$ is a Cayley group
% (respectively, Cayley $K$-group)
if and only if ${\rm Cay}(G)=1$.
% (respectively, $ {\rm Cay}_K(G)=1$).
Theorem~\ref{thm3} may thus be interpreted as
a classification of connected simple algebraic
groups of Cayley degree $1$ and, consequently,
as a first step towards the solution of the
following general problem:

\begin{problem} \label{lastproblem} {\it Find
the Cayley degrees of connected simple
algebraic groups.}
%  and, more ge\-nerally, the Cayley
% $K$-degrees of connected simple algebraic $K$-groups.}
\end{problem}

For example, composing the natural projection
$\Spin_n \to \SO_n$ with the classical Cayley
map $\SO_n \stackrel{\simeq}{\dasharrow}
{\mathfrak so}_n$ yields a generalized Cayley
map $\Spin_n \to \SO_n \dasharrow {\mathfrak
{so}}_n = {\mathfrak {spin}}_n$ of degree $2$.
Combining this with Theorem~\ref{thm2}, we
conclude that
\[ {\rm Cay}({\bf Spin}_n)=
\begin{cases} 2&\text{for $n\geqslant 6$},\\
1&\text{for $n\leqslant 5$}.
\end{cases} \]
Other examples can be found in~\cite[Section
10]{lpr}. Note that
Definition~\ref{cayley-deg} and
Problem~\ref{lastproblem} have natural
analogues in the case where $G$ is defined
over a subfield $K$ of $k$ (here we consider
only generalized Cayley maps $\varphi$ defined
over $K$).
%%We intend to address
%%Problem~\ref{lastproblem} and its variants in
%%a separate publication.

\section*{\bf Appendix.
Alternative Proof of Proposition {
\ref{prop.sl_3}}: An outline}
\renewcommand{\theequation}{A\arabic{equation}}
\setcounter{equation}{0}

\smallskip

{\it Step} 1. Consider $D$, see \eqref{D3}, as
an open subset of $ \bbP^3$ given by $x_0\neq
0$, and extend the $\Sym_3$-action
\eqref{sigma-action} up to $\bbP^3$ by
\begin{equation*}
\sigma(a_0: a_1: a_2: a_3)= (a_0:
a_{\sigma(1)}: a_{\sigma(2)}: a_{\sigma(3)}),
\quad\text{where}\hskip 2mm \sigma\in \Sym_3.
\end{equation*}

The closure $X$ of $T$ in $\bbP^3$, see
\eqref{tori}, is the rational cubic surface
given by $ x_1x_2x_3-x_0^3=0$. It has exactly
three fixed points
$$
a_i:=(1:\varepsilon^i:\varepsilon^i:\varepsilon^i),\
i=1, 2, 3, \quad \varepsilon^3=1,\
\varepsilon\neq 1,
$$
and three singular (double) points
$$
s_1=(0: 1: 0: 0),\ s_2=(0: 0: 1: 0), \ s_3=(0:
0: 0: 1).
$$
The hyperplane section of $X$ given by $x_0=0$
is $H:=l_1+l_2+l_3$, where $l_i$ is the line
given by $x_0=x_i=0$.

\begin{center}
\leavevmode \epsfxsize =2.8cm
\epsffile{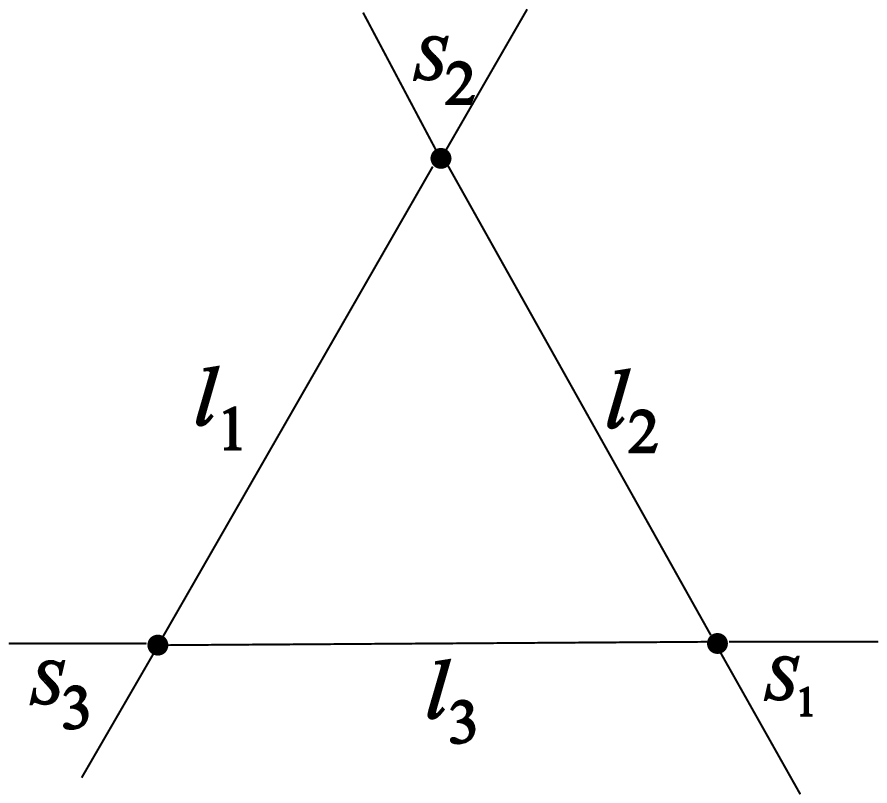}
\end{center}
Since $H$ is $\Sym_3$-invariant, the
$\Sym_3$-action on $X$ lifts to the surface
$\widetilde X$ obtained from $X$ by the
simultaneous blowing up $\mu: \widetilde X\to
X$ of $s_1$, $s_2$,~$s_3$. The surface
$\widetilde X$ is smooth and $T$ is its open
$\Sym_3$-stable subset. %%

\smallskip

{\it Step} 2. We have
$\mu^*(H)=\sum_it_i+\sum_{ij}m_{ij}$ where
$t_i$ is the proper inverse image of $l_i$ and
$\mu^{-1}(s_i)=m_{ij}\cup m_{ir}$, $\{i, j,
r\}=\{1, 2, 3\}$. The curves $t_i$, $m_{ij}$
are isomorphic to $\bbP^1$ and form a $9$-gon
as shown on the figure below. Their
intersections are transversal and the
self-intersection indices are $(t_i, t_i)=-1$,
$(m_{ij}, m_{ij})=-2$.

\begin{center}
\leavevmode \epsfxsize =4.3cm
\epsffile{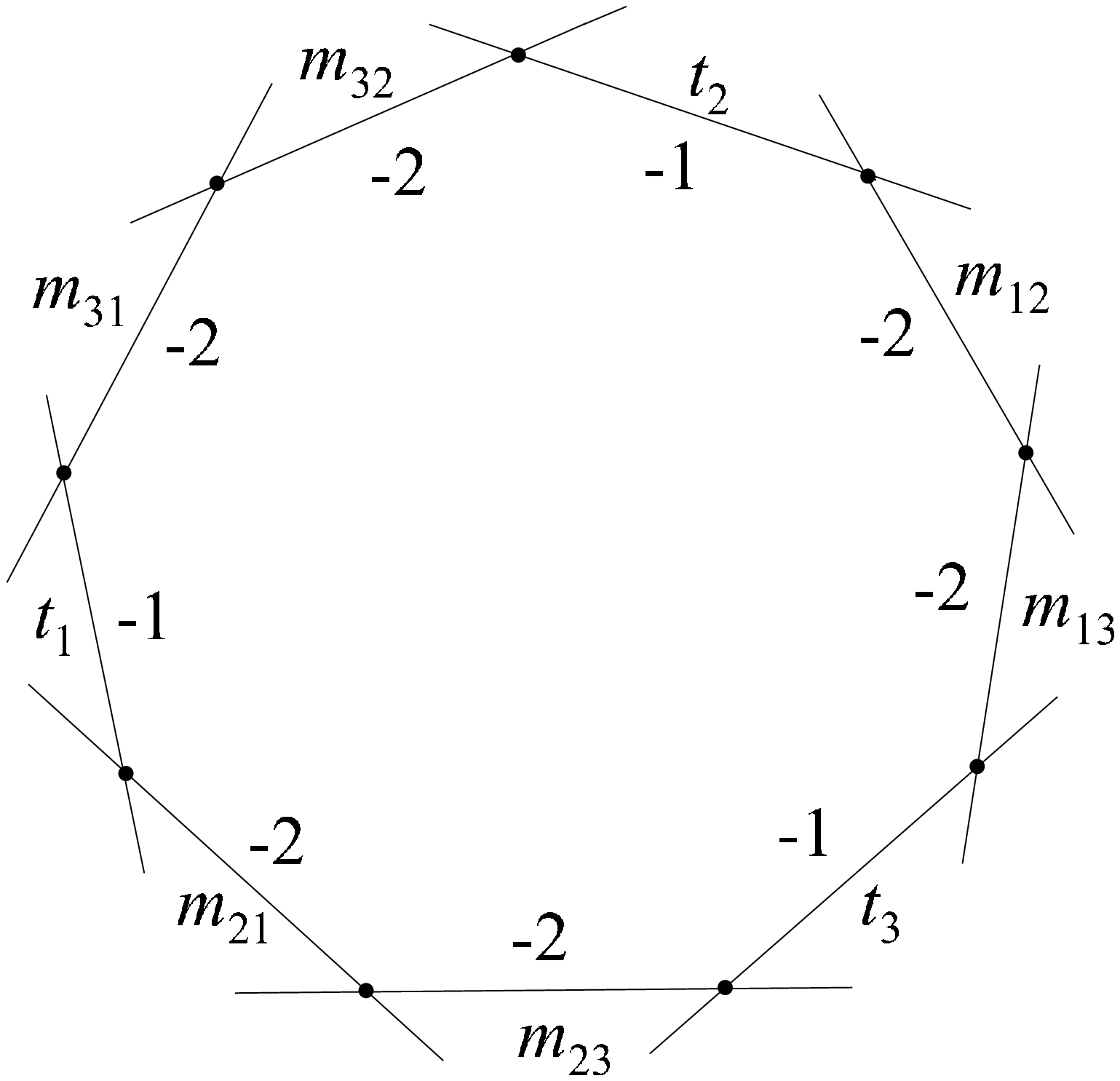}
\end{center}

Computing the canonical classes gives $K_X=-H$
and $K_{\widetilde X}=\mu^*(-H)$. Hence
\begin{equation}\label{KY2}
(K_{\widetilde X}, K_{\widetilde X})= (-H, -H)
={\rm deg}\,X=3.
\end{equation}

{\it Step} 3. By the Castelnuovo criterion,
the curves $t_i$ are exceptional, so they can
be simultaneously blown down: $\nu\!:
\widetilde X\to Y$. The surface $Y$ is smooth,
and the $\Sym_3$-invariance of $t_1+t_2+t_3$
implies that the action of $\Sym_3$ on
$\widetilde X$ descends to $Y$. We can
consider $T$ as an open $\Sym_3$-stable subset
of $Y$.

It follows from \eqref{KY2} that
\begin{equation}\label{KZ6}
(K_Y, K_Y)=6, \end{equation} and ${\rm
Pic}\;T=0$ implies that $({\rm
Pic}\;Y)^{\Sym_3}$ is generated by
\begin{equation*}\textstyle
P:=\nu_*(\sum_{ij}m_{ij}).\end{equation*}
Hence $K_Y=nP$ for some non-zero integer $n$.
Rationality of $Y$ implies $n<0$. If $C$ is a
positive divisor on $Y$, then
$\sum_{\sigma\in{\rm S}_3}\sigma(C)=cP$ for
some positive integer $c$. Using \eqref{KZ6},
we then obtain
\begin{eqnarray*}\textstyle(-K_Y, C)
&=&
%%\frac{1}{6}
\bigl(-K_Y, \sum_{\sigma\in{\rm
S}_3}\sigma(C)\bigr)/6\\&=&-
%%\frac{cn}{6}
cn(P, P)/6= -
%%\frac{c}{6n}
c(K_Y, K_Y)/6n= -
%%\frac{c}{n}
c/n>0.\end{eqnarray*} By the Nakai--Moishezon
criterion, this implies that $-K_Y$ is ample,
i.e., $Y$ is a Del Pezzo surface. From
\eqref{KZ6} it then follows  (see,
e.g.,\,\cite[\S\,24]{manin2}) that $|-K_Y|$
defines an embedding of $Y$ into $\mathbb P^6$
equivariant with respect to a certain action
of ${\rm S}_3$ on $\mathbb P^6$. We keep the
notation $Y$ for its image.

\smallskip

{\it Step} 4. Consider on $Y$ the linear
system $|R|$ of all hyperplane sections in
$\bbP^6$ containing the fixed point $a_1\in
T\subseteq Y$ and which is singular at $a_1$. These are
precisely sections by hyperplanes tangent to
$Y$ at $a_1$, so
\begin{equation}\label{R4}
\dim |{R}|=4.
\end{equation}

 The system $|R|$ is an $\Sym_3$-stable
subsystem of $|{-K_{Y}}|$. By Bertini's
theorem, its general element $R$ is an
irreducible curve. We have
\begin{equation}\label{R5}
\textstyle p_{\rm a}(R)=1+
%%\frac{
(R, (R+K_Y))/2=
%%}{2}+1=
1+
%%\frac{
(R, (R-R))/2
%%}{2}+1
=1.
\end{equation}
On the other hand, $p_{\rm
a}(R)=g+\sum_x\delta_x$, where $g$ is the
geometric genus of the normalization of $R$,
the sum is taken over all singular points $x$
of $R$, and $\delta_x>0$. This and \eqref{R5}
imply that $R$ is a rational curve whose
singular locus is the double point $a_1$.

 The system $|R|$ has no fixed components.
Indeed, if $F$ were such a component, then
$\dim {\rm H}^0(Y,\mathcal O(F))=1$ and, by
the Riemann--Roch theorem,
\begin{equation}\label{R6}
\dim {\rm H}^0(Y,\mathcal O(K_Y-F))\geqslant
\bigl((F, F)-(F, K_Y)\bigr)/2.
\end{equation}
Let $-K_Y=F+E$. Since $F>0$ and $E>0$, the
left-hand side of \eqref{R6} is zero, whence
$0\geqslant (F, F)+ (F, E)/2$. Since $(F,
E)\geqslant 0$, this yields $0\geqslant (F,
F)$. But $F=m P$ for some non-zero integer $m$.
Therefore $$0\geqslant (F, F)=m^2(P, P)=
6m^2/n^2>0,$$  a contradiction.

From \eqref{R4} we deduce that $a_1$ is a
unique base point of $|R|$.

\smallskip

{\it Step} 5. Let $\gamma\!:\widetilde Y\to Y$
be the blowing up of $a_1$. The action of
$\Sym_3$ lifts to $\widetilde Y$. The proper
inverse image $|{\widetilde R}|$ of $|{R}|$ is
a $4$-dimensional $\Sym_3$-stable linear
system on $\widetilde Y$. It has no base
points and separates points of an open subset
of $\widetilde Y$. Hence $|{\widetilde R}|$
defines an $\Sym_3$-equivariant morphism
$\widetilde Y\to \bbP^3$ with respect to a
certain $\Sym_3$-action on $\bbP^3$. Let $Z$
be its image. This morphism then yields an
$\Sym_3$-equivariant birational isomorphism
$\psi: \widetilde Y\rightarrow Z$.

Let $l=\gamma^{-1}(a_1)$ and let $\widetilde
R$ be the proper inverse image of $R$. Then
$(l, l)=-1$ and, since $a_1$ is a double point
of $R$, we have $\gamma^*(R)=\widetilde R+2l$
and $(l, \widetilde R)=2$. This yields $$6=(R,
R)= (\widetilde R, \widetilde R)+4(l,
\widetilde R)+4(l, l)= (\widetilde R,
\widetilde R)+4,$$ so $(\widetilde R,
\widetilde R)=2$. Since ${\rm
deg}\,Z=(\widetilde R, \widetilde R)$, this
means that $Z$ is an $\Sym_3$-stable quadric
in $\bbP^3$.

\smallskip

{\it Step} 6. Since the point
$a'_2:=\psi\circ\gamma^{-1}(a_2)\in Z$ is
fixed by $\Sym_3$, it follows from the
complete reducibility of representations of
reductive groups that there is an
$\Sym_3$-stable plane $L\simeq \bbP^2$ in
$\bbP^3$ not passing through $a'_2$. Consider
the stereographic projection $ \pi\!:
Z\dasharrow L$ from $a'_2$; it is birational
and $\Sym_3$-equivariant. The map $\pi$ is
defined at $\psi\circ\gamma^{-1}(a_3)$ and
$a'_3:=\pi\circ\psi\circ\gamma^{-1}(a_3)\in L$
is a fixed point of $\Sym_3$. Using the
complete reducibility argument again, we
conclude that there is an $\Sym_3$-stable line
$l\subset L$ such that $a'_3\in L\setminus l$.
Thus we obtain a faithful linear action of
$\Sym_3$ on $\bbA^2\simeq L\setminus l$. But
there is  a unique $2$-dimensional faithful
linear representation of $\Sym_3$, namely that
on~$\lt$ given by \eqref{sigma-action},
\eqref{tori}.

\smallskip
In summary, we have constructed the following
chain of birational equivariant maps of
$\Sym_3$-varieties:
\[ \lt \hookrightarrow L
\overset{\pi}\dashleftarrow Z
\overset{\psi}\longleftarrow \widetilde{Y}
\overset{\gamma}\longrightarrow Y
\overset{\nu}\longleftarrow \widetilde{X}
\overset{\mu}\longrightarrow X \hookleftarrow
T.
\]
This shows that $T$ and $\lt$ are birationally
isomorphic as $\Sym_3$-varieties, thus
completing the proof of
Proposition~\ref{prop.sl_3}. \quad $\square$

\section*{\bf Acknowledgements} We are
grateful to {\sc G.  Berhuy}, {\sc V. A.
Iskovskikh}, and {\sc D.  Luna} for
stimulating discussions related to the subject
matter of this paper.

\providecommand{\bysame}{\leavevmode\hbox
to3em{\hrulefill}\thinspace}

\end{document}